\documentclass[a4paper,12pt]{article}

\usepackage{latexsym}
\usepackage{amssymb}
\usepackage{theorem}
\usepackage{amsmath}
\usepackage{amscd}
\usepackage{graphicx}
\usepackage{xcolor}
\usepackage{url}
\usepackage{enumitem}
\usepackage{tikz}
\usetikzlibrary{arrows.meta, calc}
\usepackage{diagbox} % 用于对角线单元格
\usepackage{xcolor}  % 用于颜色
\usepackage{array}
\usepackage[
colorlinks=false,
pdfborder={0 0 1},
linkbordercolor=red,
citebordercolor=green
]{hyperref}
\usetikzlibrary{decorations.pathmorphing}
\newtheorem{theorem}{Theorem}[section]
\newtheorem{corollary}[theorem]{Corollary}
\newtheorem{lemma}[theorem]{Lemma}
\newtheorem{example}[theorem]{Example}
\newtheorem{proposition}[theorem]{Proposition}
\newtheorem{remark}[theorem]{Remark}
\newtheorem{definition}[theorem]{Definition}

\newcommand{\demo}{\par\noindent{\it Proof. \/}\ }
\newcommand{\enD}{\hfill $\Box$\vspace{3truemm} \par}
\newcommand{\R}{\mathbb{R}}

\newcommand{\be}{\mbox{\boldmath $e$}}

\newcommand{\bv}{\mbox{\boldmath $v$}}
\newcommand{\bx}{\mbox{\boldmath $x$}}
\newcommand{\by}{\mbox{\boldmath $y$}}

\newcommand{\bgamma}{\mbox{\boldmath $\gamma$}}
\newcommand{\bnu}{\mbox{\boldmath $\nu$}}
\newcommand{\bmu}{\mbox{\boldmath $\mu$}}

\newcommand{\tikznode}[2]{%
\tikz[remember picture,baseline=(#1.base)] \node (#1) [inner sep=0pt] {#2};%
}
\begin{document}

\title{Evolutes, involutes and parallels  of Legendre curves in hyperbolic and de Sitter 2-spaces via rotational transformations}

\author{Nozomi Nakatsuyama, Masatomo Takahashi and Anjie Zhou}

\date{\today}

\maketitle
\begin{abstract} 
We develop a unified approach to constructing the moving frames of the evolute, involute and parallel of a Legendre curve via rotational transformations (rotation, boost and null rotation). The definitions of the evolute, involute and parallel of a Legendre curve in hyperbolic and de Sitter 2-spaces are systematically presented. Furthermore, the connections and dualities among these curves are discussed. 
\end{abstract}
\renewcommand{\thefootnote}{\fnsymbol{footnote}}
{\hypersetup{pdfborder={0 0 0}}
	}

%%%%%%%%%%%% Section 1 %%%%%%%%%%%%%
\section{Introduction}
The Lorentz group $O(1,2)$ is the linear isometry group of Lorentz-Minkowski 3-space and serves as the basic symmetry group underlying special relativity. The   transformations in its identity component $O^{+\uparrow}(1,2)$ are classified, up to conjugacy, into three types of rotational transformations according to whether the axis is timelike, spacelike or lightlike.
A \textit{rotation} (or, \textit{elliptic transformation}) fixes a timelike axis and rotates the orthogonal spacelike plane through an angle, exactly as in Euclidean geometry. 
A \textit{boost} (or, \textit{hyperbolic transformation}) in the direction of a unit spacelike vector $\bv$ is a Lorentz transformation that fixes pointwise the  subspace orthogonal to both the timelike direction and  $\bv$, mixing the timelike direction with  $\bv$ through a rapidity parameter.
A \textit{null rotation} (or, \textit{parabolic transformation}) fixes pointwise a lightlike direction and has no analogue in Euclidean geometry or Galilean relativity (cf. \cite{CL2021}). 
The hyperbolic and de Sitter 2-spaces are orbits of $ O^{+\uparrow}(1,2) $ acting on  Lorentz-Minkowski 3-space, corresponding respectively to the unit timelike and spacelike pseudospheres. The duality between these two spaces has also been studied extensively (cf. \cite{CIT2024,CPT2020,CT2016,I2009}).
The extension of classical evolutes, involutes and parallels in the hyperbolic and de Sitter 2-spaces has been developed in  \cite{CIT2024,CT2016,IPST2004,IT2007,XC2026,ZP2026}. In our previous work \cite{NTZ2026}, we introduced the notions of horocyclic parallel and involute via Bernoulli and Riccati equations, respectively.
\par 
In the present paper, we revisit and extend these constructions from a different viewpoint. We show that the moving frames of evolutes, involutes  and parallels arise uniformly as  the image of the moving frame of the original curve in hyperbolic 2-space under a suitable Lorentz rotation, as shown in section \ref{S3}. This rotational viewpoint yields results that are hard to obtain by a direct computation. For instance, it gives the explicit rotation relating the hyperbolic horocyclic and hyperbolic evolutes in hyperbolic 2-space (Propositions \ref{prop1} and  \ref{prop2}), from which the relations between their moving frames, curvatures  and singularities can also be derived (Theorem \ref{th2} and Corollary \ref{cor1}).   Based on  de Sitter horospherical geometry  in \cite{HI2015,K2009,K2010},
sections \ref{S4} and \ref{S5} carry out the construction systematically in de Sitter 2-space. 
Evolutes and involutes have long been closely related  and their correspondence has been extensively studied in various ambient spaces (cf. \cite{FT2015,LP2022,NTZ2026,SLP2022,ZP2026}). Using rotation matrices, we generalize the  correspondence and further establish the relations between their moving frames in hyperbolic and de Sitter 2-spaces (Theorem \ref{th1}).
\par
We shall assume throughout the whole paper that all maps and manifolds are $C^{\infty}$ unless the contrary is explicitly stated.

%%%%%%%%%%% Acknowledgement %%%%%%%%%%%%
\bigskip
\noindent
{\bf Acknowledgement}. 
The first author is supported by JST SPRING (Grant Number JPMJSP2153).
The second author is partially supported by JSPS KAKENHI (Grant Number JP 24K06728).
The third author is grateful for the financial support provided by the CSC-MuroranIT Scholarship (Grant  Number 202506620020).

%%%%%%%%%%% Section 2 %%%%%%%%%%%%%
\section{Preliminaries}\label{S2}

Let $\R_1^3 $ denote \textit{Lorentz}-\textit{Minkowski} $3$-\textit{space} with the pseudo scalar product $\langle \bx,\by\rangle=-x_1y_1+x_2y_2+x_3y_3$ and the pseudo vector product 
$$
\bx \wedge \by =\det\begin{pmatrix}-\be_1&\be_2 &\be_3\\
x_1&x_2& x_3\\
y_1&y_2&y_3 
\end{pmatrix},
$$ 
where $\bx=(x_1,x_2,x_3)$, $\by=(y_1,y_2,y_3)$ and $\left\{\be_1,\be_2, \be_3\right\}$ is the canonical basis of $\R_1^3$. 
A vector $\bx \in\R_1^3 \setminus \left\{0\right\} $ is \textit{spacelike}, \textit{lightlike} or \textit{timelike} if $\langle \bx,\bx \rangle > 0$, $\langle \bx,\bx \rangle=0$ or $\langle \bx,\bx \rangle<0$, respectively. 
We define the \textit{hyperbolic} 2-\textit{space} as
$
H^2=\left\{ \bx\in\R_1^3 \mid \langle \bx,\bx\rangle=-1\right\}
$ and the \textit{de Sitter} 2-\textit{space} as
$
S_1^2=\left\{ \bx\in\R_1^3 \mid \langle \bx,\bx\rangle=1\right\}.
$
It can be verified that
$\langle \bx,\by_1 \wedge \by_2  \rangle=\det(\bx,\by_1,\by_2),$ implying that $\by_1 \wedge \by_2$ and $\by_i$ are pseudo-orthogonal $ (i=1,2).$ 
We denote
$
\Delta_1=\left\{ (\bnu_1,\bnu_2) \in H^2\times S_1^2\mid \langle \bnu_1,\bnu_2 \rangle =0 \right\}  
$, $
\overline{\Delta}_1=\left\{ (\bnu_1,\bnu_2) \in S_1^2 \times H^2 \mid \langle \bnu_1,\bnu_2 \rangle =0 \right\}  
$  and $
\Delta_5=\left\{ (\bnu_1,\bnu_2) \in  S_1^2\times S_1^2\mid \langle \bnu_1,\bnu_2 \rangle =0 \right\}  
$.

In this paper, we denote $I$ as interval  in $\R$.

%%%%%
\begin{definition}[\cite{CT2016}] \rm
Let $\bgamma_h : I \rightarrow H^2$ be a smooth curve.
We call $(\bgamma_h,\bgamma_h^d): I \rightarrow \Delta_1$ a \textit{spacelike Legendre curve} in $\Delta_1$ if there exists a smooth mapping $\bgamma_h^d : I \rightarrow S^2_1$ such that	$\langle \dot{\bgamma}_h(t),\bgamma_h^d(t) \rangle=0$ for all $t\in I$. Then we call  $\bgamma_h :I\rightarrow H^2$   a \textit{spacelike frontal} in $H^2$.
\end{definition}
%%%%%

Let $\bgamma_h^s(t) = \bgamma_h(t) \wedge \bgamma_h^d(t)  \in S_1^2$. Then $\{\bgamma_h(t),\bgamma_h^d(t),\bgamma_h^s(t)\}$ is a moving frame along $\bgamma_h(t)$. 
We have the following hyperbolic Legendre Frenet-type formula:
\begin{equation} \notag
\begin{pmatrix}
\dot{\bgamma}_h(t) \\
\dot{\bgamma}_h^d(t) \\
\dot{\bgamma}_h^s(t) 
\end{pmatrix}
=
\begin{pmatrix}
0&0&m_h(t) \\
0&0&n_h (t) \\
m_h(t) &-n_h(t) &0
\end{pmatrix}
\begin{pmatrix}
\bgamma_h(t) \\
\bgamma_h^d(t) \\
\bgamma_h^s(t)
\end{pmatrix}.
\end{equation}
The mapping $(m_h,n_h)$ is called the   \textit{spacelike hyperbolic Legendre curvature} of   $(\bgamma_h,\bgamma_h^d)$.

%%%%%
\begin{definition}[\cite{CT2016}] \rm
Suppose that $\bgamma_d : I \to S_1^2$ is a spacelike curve at regular points $t \in I$, namely, $\dot{\bgamma}_d(t)$ is a spacelike vector at the regular points.
We call $(\bgamma_d,\bgamma_d^h): I \rightarrow \overline{\Delta}_1$ a \textit{spacelike Legendre curve} in $\overline{\Delta}_1$ if there exists a smooth mapping $\bgamma_d^h : I \to H^2$ such that	$\langle \dot{\bgamma}_d(t),\bgamma_d^h(t) \rangle=0$ for all $t\in I$. Then we call   $\bgamma _d:I\rightarrow S_1^2$   a \textit{spacelike frontal} in $S_1^2$.
\end{definition}
%%%%%

Let  $\bgamma_d^s(t) = \bgamma_d(t) \wedge \bgamma_d^h(t) \in S_1^2$. Then $\{\bgamma_d(t),\bgamma_d^h(t),\bgamma_d^s (t)\}$ is a moving frame along $\bgamma_d(t)$.  We have the following spacelike de Sitter Legendre Frenet-type formula:
\[
\begin{pmatrix} 
\dot{\bgamma}_d(t) \\ 
\dot{\bgamma}_d^h(t) \\ 
\dot{\bgamma}_d^s(t) 
\end{pmatrix} = 
\begin{pmatrix} 
0 & 0 & m_d (t) \\ 
0 & 0 & n _d(t) \\ 
-m_d (t) & n_d (t) & 0 
\end{pmatrix} 
\begin{pmatrix} 
\bgamma_d(t) \\ 
\bgamma_d^h(t) \\ 
\bgamma_d^s(t) 
\end{pmatrix}.
\]
We call   $(m_d, n_d)$ the \textit{spacelike de Sitter Legendre curvature} of   $(\bgamma_d,\bgamma_d^h)$. 

%%%%%
\begin{definition}[\cite{CT2016}] \rm
Suppose that $\bgamma_T : I \to S_1^2$ is a timelike curve at regular points $t \in I$, namely, $\dot{\bgamma}_T(t)$ is a timelike vector at the regular points.
We call $(\bgamma_T,\bgamma_T^d): I \rightarrow \Delta_5$ a \textit{timelike Legendre curve} in $\Delta_5$ if there exists a smooth mapping $\bgamma_T^d : I \to S_1^2$ such that	$\langle \dot{\bgamma}_T(t),\bgamma_T^d(t) \rangle=0$ for all $t\in I$. Then we call  $\bgamma _T:I\rightarrow S_1^2$   a \textit{timelike frontal} in $S_1^2$.
\end{definition}
%%%%%

Let  $\bgamma_T^h(t) = \bgamma_T(t) \wedge \bgamma_T^d(t) \in H^2$. Then $\{\bgamma_T(t),\bgamma_T^d(t),\bgamma_T^h (t)\}$ is a moving frame along $\bgamma_T(t)$.  We have the following   timelike de Sitter Legendre Frenet-type formula:
\[
\begin{pmatrix} 
\dot{\bgamma}_T(t) \\ 
\dot{\bgamma}_T^d(t) \\ 
\dot{\bgamma}_T^h(t) 
\end{pmatrix} = 
\begin{pmatrix} 
0 & 0 & m_T (t) \\ 
0 & 0 & n _T(t) \\ 
m _T(t) & n _T(t) & 0 
\end{pmatrix} 
\begin{pmatrix} 
\bgamma_T(t) \\ 
\bgamma_T^d(t) \\ 
\bgamma_T^h(t) 
\end{pmatrix}.
\]
We call   $(m_T, n_T)$ the \textit{timelike de Sitter Legendre curvature} of   $(\bgamma_T, \bgamma_T^d)$.

The \textit{Lorentz group} is defined as
$$O(1,2) =\left\{A\in M_3(\R)\;\middle|\; A^TGA = G, \; G =
\begin{pmatrix}
-1 & 0 & 0  \\
0 & 1 & 0  \\
0 & 0 & 1 
\end{pmatrix}\right\} $$
and the \textit{special Lorentz group}   is
$ SO(1,2) = \{ A \in O(1,2) \mid \det A  = 1 \}.$  Here, $A^T$ denotes the transpose of $A$. The \textit{special orthochronous Lorentz group} is given by$$O^{+\uparrow}(1,2)=\left\{ A=(a_{ij})_{1\leq i,j\leq3} \in SO(1,2) \mid a_{11} > 0 \right\},$$and it preserves time orientation. Transformations in $O^{+\uparrow}(1,2)$ can be classified, up to conjugacy, into three types: \textit{rotations}, \textit{boosts}, and \textit{null rotations}, representing Lorentz rotations around timelike, spacelike  and lightlike directions, respectively (cf. \cite{CL2021}). Next, we examine their matrix expressions under specific bases.

%%%%%
\begin{lemma}[\cite{CL2021}]\rm
Let  $\{\be_1,\be_2, \be_3\}$  be the canonical basis   of $\R_1^3$.
\begin{enumerate}	
\item[$(1)$] Let $(\bx_1,\bx_2,\bx_3)$ be a pseudo-orthonormal basis of $\R_1^3$ such that $(\bx_1,\bx_2,\bx_3 )=(\be_1,\be_2, \be_3)$. \\
{\rm(i)}    The   matrix of  {rotation}      with     angle $\phi  \in \mathbb{R}$   is  
\[
R (\phi ) = \begin{pmatrix} 1 & 0 & 0 \\ 0 & \cos\phi  & -\sin\phi \\ 0 & \sin \phi  & \cos\phi  \end{pmatrix}.
\]
{\rm(ii)}    
The   matrices of   {boosts}   in the $\bx_2$-direction and $\bx_3$-direction   with   rapidity $\phi  \in \mathbb{R}$    are respectively  
$$B_1(\phi ) = \begin{pmatrix} \cosh\phi  & \sinh\phi  & 0  \\ \sinh\phi   & \cosh\phi & 0\\ 0  & 0& 1\end{pmatrix},\quad  
B_2(\phi ) = \begin{pmatrix} \cosh\phi  & 0 & \sinh\phi  \\ 0 & 1 & 0 \\ \sinh\phi  & 0 & \cosh\phi  \end{pmatrix}.
$$
{\rm(iii)}  The   matrices of   {null rotations}   around the $(\bx_1\pm\bx_2)$-direction and $(\bx_1\pm\bx_3)$-direction   with     angle $\phi  \in \mathbb{R}$    are respectively  
$$
N^\pm_1(\phi ) = \begin{pmatrix} 
1 + \frac{\phi^2 }{2} & \mp \frac{\phi^2 }{2} & \phi  \\ 
\pm \frac{\phi^2 }{2} & 1 - \frac{\phi^2 }{2} & \pm \phi  \\ 
\phi & \mp \phi  & 1 
\end{pmatrix}, \quad  N^\pm_2(\phi )   = \begin{pmatrix}
1 + \frac{\phi^2 }{2} &  \phi  & \mp \frac{\phi^2}{2} \\
\phi  & 1 & \mp \phi  \\
\pm \frac{\phi^2 }{2} & \pm \phi  & 1 - \frac{\phi^2 }{2}
\end{pmatrix}.$$
\item[$(2)$] Let $(\bx_1,\bx_2,\bx_3)$ be a pseudo-orthonormal basis of $\R_1^3$ such that  $(\bx_1,\bx_2,\bx_3 )=(\be_3,\be_2, \be_1) $. \\
{\rm(i)}    The   matrix of {rotation}     with     angle $\phi  \in \mathbb{R}$   is  
\[
\overline{R}(\phi ) = \begin{pmatrix}  \cos\phi  & -\sin\phi & 0\\  \sin \phi  & \cos\phi & 0 \\  0 & 0  &1   \end{pmatrix}.
\]
{\rm(ii)}    
The    matrix  of   {boost}   in the $\bx_1$-direction    with      rapidity $\phi  \in \mathbb{R}$    is 
$$ 
\overline{B}_1(\phi ) = \begin{pmatrix} \cosh\phi  & 0 & \sinh\phi  \\ 0 & 1 & 0 \\ \sinh\phi  & 0 & \cosh\phi  \end{pmatrix}.
$$
{\rm(iii)}  The    matrix of  {null rotation}   around the  $(\bx_1\pm\bx_3)$-direction   with      angle $\phi  \in \mathbb{R}$    is 
$$
\overline{N}^\pm_1(\phi )   = \begin{pmatrix}
1 - \frac{\phi^2}{2} &  \phi & \pm \frac{\phi^2}{2} \\
-\phi & 1 & \pm \phi \\
\mp \frac{\phi^2}{2} & \pm \phi & 1+ \frac{\phi^2 }{2}
\end{pmatrix}.$$
\end{enumerate}
\end{lemma}
%%%%%

We say that a smooth curve $\bgamma_h:(I,t_0)\rightarrow H^2$ at $t_0$ has a $(i,j)$-cusp, where $ (i, j)= (2, 3),(2, 5)$, $(3, 4),(3, 5) $  if $\bgamma$ is $ \mathcal{A} $-equivalent to the germ  $t \mapsto(t^i,t^j)$ at the origin. 
%%%%%
\begin{lemma}[\cite{NTZ2026}]\label{lemma1}
	Let $(\bgamma_h,\bgamma_h^d):I\rightarrow \Delta_1$ be a spacelike Legendre curve with curvature $(m_h,n_h)$.  Then we have the following.
	\begin{enumerate}%[label=\textup{(\arabic*)}]
		\item[$(1)$]  $\bgamma_h$ is singular at $t_0$  if and only if $m_h(t_0)=0$.
		\item[$(2)$]  $\bgamma_h$ has a $(2,3)$-cusp at $t_0$ if and only if $m_h(t_0)=0$, $ n_h(t_0)\neq0$, $\dot{m}_h(t_0)\neq0$.
		\item[$(3)$] $\bgamma_h$ has a $(3,4)$-cusp at $t_0$ if and only if  $m_h(t_0)=\dot{m}_h(t_0)=0$, $ n_h(t_0)\neq0$, $\ddot{m}_h(t_0)\neq0$.
		\item[$(4)$] $\bgamma_h$ has a $(2,5)$-cusp at $t_0$ if and only if  $m_h(t_0)=n_h(t_0)=0$, $\dot{m}_h(t_0)\neq0$, $\dot{m}_h(t_0)\ddot{n}_h(t_0)-\ddot{m}_h(t_0)\dot{n}_h(t_0)\neq0$.
		\item[$(5)$] $\bgamma_h$ has a $(3,5)$-cusp at $t_0$ if and only if  $m_h(t_0)=n_h(t_0)=\dot{m}_h(t_0)=0$, $\ddot{m}_h(t_0)\neq0$, $\dot{n}_h(t_0)\neq0$.
	\end{enumerate} 
\end{lemma}
%%%%%
%%%%%%%%%%% Section 3 %%%%%%%%%%%%%
\section{Evolutes, involutes and parallels of spacelike Legendre curves in $H^2$}\label{S3}
\subsection{Evolutes of spacelike Legendre curves}

%%%%%
\begin{definition}[Hyperbolic and de Sitter horocyclic evolutes \cite{CIT2024}]\rm
Let $(\bgamma_h,\bgamma_h^d):I\rightarrow\Delta_1$ be a spacelike Legendre curve  with curvature $(m_h,n_h)$. Suppose that there exists a smooth function $f:I\rightarrow\R  $ such that $m_h(t)+f(t)n_h(t)=0$ for all $t\in I$. 
Then the \textit{hyperbolic horocyclic evolute}   $\mathcal{HE}^\pm(\bgamma_h,\bgamma_h^d):I\rightarrow H^2$  of $(\bgamma_h,\bgamma_h^d)$   is defined as  
\begin{equation}\notag
\mathcal{HE}^\pm(\bgamma_h,\bgamma_h^d)(t)=\bgamma_h(t)+f(t)\bgamma_h^d(t)+\frac{f^2(t)}{2}(\bgamma_h(t)\pm\bgamma_h^s(t)).
\end{equation} 
The \textit{de Sitter horocyclic evolute}   $\mathcal{DE}^\pm(\bgamma_h,\bgamma_h^d):I\rightarrow S_1^2$    of $(\bgamma_h,\bgamma_h^d)$   is defined as  
\begin{equation}\notag
\mathcal{DE}^\pm(\bgamma_h,\bgamma_h^d)(t)=\mp\frac{f^2(t)}{2}\bgamma_h(t)\mp f(t)\bgamma_h^d(t)+ \left(1-\frac{f^2(t)}{2}\right)\bgamma_h^s(t).
\end{equation} 
\end{definition}
%%%%%

%%%%%
\begin{proposition}[\cite{CIT2024,FT2015}]
We denote  $\mathrm{Reg}(\bgamma_h^d)= \{t\in I\mid n_h(t)\neq0 \}$ as the set of regular points of  $\bgamma_h^d$. Suppose that there exists a smooth function $f:I\rightarrow\R  $ such that $f(t)=-m_h(t)/n_h(t)$ on $\mathrm{Reg}(\bgamma_h^d)$. Then   $f$ is
unique if and only if $\mathrm{Reg}(\bgamma_h^d) $ is   dense  in $ I $. 
\end{proposition}
%%%%%

The uniqueness condition of $f$ is a topological condition. In
this paper, we assume that $\mathrm{Reg}(\bgamma_h^d)= \{t\in I\mid n_h(t)\neq0 \}$  is   dense  in $ I $.

%%%%
\begin{proposition}[\cite{CIT2024}] 
Suppose that there exists   $f:I\rightarrow\R  $ such that $m_h(t)+f(t)n_h(t)=0$ for all $t\in I$. Then
\begin{enumerate} 
\item[$(1)$] $(\mathcal{HE}^\pm(\bgamma_h,\bgamma_h^d), \mathcal{DE}^\pm(\bgamma_h,\bgamma_h^d)  ):I\rightarrow \Delta_1$ is a spacelike Legendre curve with curvature 
$$
\left(m^\pm_\mathcal{HE}(t), n^\pm_\mathcal{HE}(t)\right)= \left( - \dot{f}(t)\pm ({f^2(t)n_h(t)}/{2}),\pm \dot{f}(t)- ({f^2(t)n_h(t)}/{2})+n_h(t)\right).
$$ 

\item[$(2)$] $(\mathcal{DE}^\pm(\bgamma_h,\bgamma_h^d), \mathcal{HE}^\pm(\bgamma_h,\bgamma_h^d) ):I\rightarrow \overline{\Delta}_1$ is a spacelike Legendre curve with curvature 
$$
\left(m^\pm_\mathcal{DE}(t), n^\pm_\mathcal{DE}(t)\right)= \left( \mp \dot{f}(t)+({f^2(t)n_h(t)}/{2})-n_h(t),  \dot{f}(t)\mp ({f^2(t)n_h(t)}/{2}) \right).
$$ 
\end{enumerate}
\end{proposition}
%%%%

%%%%
\begin{proposition} \label{prop1}
Suppose that there exists   $f:I\rightarrow\R  $ such that $m_h(t)+f(t)n_h(t)=0$ for all $t\in I$. Let  $\bmu^\pm_\mathcal{E}(\bgamma_h,\bgamma_h^d) =\mathcal{HE}^\pm(\bgamma_h,\bgamma_h^d) \wedge \mathcal{DE}^\pm(\bgamma_h,\bgamma_h^d) $ and denote the moving frame of $\mathcal{HE}^\pm(\bgamma_h,\bgamma_h^d)$ by $(\mathcal{HE}^\pm, \mathcal{DE}^\pm,\bmu^\pm_\mathcal{E} )(\bgamma_h,\bgamma_h^d) $. Then
$$ 
(\mathcal{HE}^\pm, \mathcal{DE}^\pm,\bmu^\pm_\mathcal{E}  )^T(\bgamma_h,\bgamma_h^d) =N_1^{\mp}(-f )R(- \pi/2) (\bgamma_h, \bgamma_h^d,\bgamma_h^s   ) ^T .$$
That is,  $( \mathcal{HE}^\pm(\bgamma_h,\bgamma_h^d)(t), \mathcal{DE}^\pm(\bgamma_h,\bgamma_h^d)(t),\bmu^\pm_\mathcal{E}(\bgamma_h,\bgamma_h^d)(t) )$ is obtained by applying a rotation  with   angle $ -\pi/2$ to $(\bgamma_h(t), \bgamma_h^d(t),\bgamma_h^s(t) )$, followed by a null rotation    around the $(\bgamma_h\mp \bgamma_h^s)$-direction with   angle $-f(t)$  for all $t\in I$. 
\end{proposition}
%%%%

%%%%%
\begin{definition}[Hyperbolic and de Sitter evolutes]\rm
Let $(\bgamma_h,\bgamma_h^d):I\rightarrow\Delta_1$ be a spacelike Legendre curve  with curvature $(m_h,n_h)$. 
Suppose that  
there exists a smooth function $\varphi:I\rightarrow\R$ such that $m_h(t)\cosh\varphi(t)+ n_h(t)\sinh\varphi(t)=0$ for all $t\in I$. Then 
the \textit{hyperbolic evolute} $ HE(\bgamma_h,\bgamma_h^d):I\rightarrow H^2$   of  $(\bgamma_h,\bgamma_h^d)$   is defined as  
\begin{equation}\notag
HE(\bgamma_h,\bgamma_h^d)(t)=\cosh\varphi(t)\bgamma_h(t)+\sinh\varphi(t)\bgamma_h^d(t).
\end{equation} 
Suppose    there exists a smooth function $\psi:I\rightarrow\R$ such that $m_h(t)\sinh\psi(t)+ n_h(t)\cosh\psi(t)=0$ for all $t\in I$. Then 
the \textit{de Sitter evolute} $ DE(\bgamma_h,\bgamma_h^d):I\rightarrow S_1^2$   of $(\bgamma_h,\bgamma_h^d)$  is defined as  
\begin{equation}\notag
DE(\bgamma_h,\bgamma_h^d)(t)=\sinh\psi(t)\bgamma_h(t)+\cosh\psi(t)\bgamma_h^d(t).
\end{equation} 
\end{definition}
%%%%%

If  $\varphi:I\rightarrow\R$ exists, then $\varphi$ is unique under the assumption that $\mathrm{Reg}(\bgamma_h^d)=\{t\in I\mid n_h(t)\neq0\}$  is   dense  in $ I $. In this paper, we also assume that   $\mathrm{Reg}(\bgamma_h)=\{t\in I\mid m_h(t)\neq0\}$  is   dense  in $ I $. Then, if  $\psi:I\rightarrow\R$ exists,  $\psi$ is also unique.

%%%%%
\begin{remark}\rm
If $m_h^2(t)=n_h^2(t)$ for all $t\in I$, then $\bgamma_h$ is a part of horocycle (cf. \cite{NTZ2026}).  In this case,
there exists no smooth function $\varphi:I\rightarrow\R$ such that for all $t\in I$, either $m_h(t)\cosh\varphi(t)+ n_h(t)\sinh\varphi(t)=0$ or $m_h(t)\sinh\varphi(t)+ n_h(t)\cosh\varphi(t)=0$ holds. Consequently,  the  hyperbolic and de Sitter evolutes of $(\bgamma_h,\bgamma_h^d)$   do  not exist.	 
\end{remark}
%%%%%

%%%%
\begin{proposition} 
Let $(\bgamma_h,\bgamma_h^d):I\rightarrow\Delta_1$ be a spacelike Legendre curve  with curvature $(m_h,n_h)$. 
\begin{enumerate} 
\item[$(1)$] 
Suppose	that  there exists   $\varphi:I\rightarrow\R$ such that $m_h(t)\cosh\varphi(t)+ n_h(t)\sinh\varphi(t)=0$ for all $t\in I$.   Then  
$(HE(\bgamma_h,\bgamma_h^d), \bnu_{HE}(\bgamma_h,\bgamma_h^d)  ):I\rightarrow \Delta_1$ is a spacelike Legendre curve with curvature 
$ 
(m _{HE}(t), n _{HE}(t))= ( -  \dot{\varphi}(t), m_h(t)\sinh\varphi(t)+ n_h(t)\cosh\varphi(t)),
$ 
where
$  
\bnu_{HE}(\bgamma_h,\bgamma_h^d)(t) = \bgamma_h^s(t). 
$ 
\item[$(2)$] Suppose	that  there exists   $\psi:I\rightarrow\R$ such that $m_h(t)\sinh\psi(t)+ n_h(t)\cosh\psi(t)=0$ for all $t\in I$.   Then  
$(DE(\bgamma_h,\bgamma_h^d), \bnu_{DE} (\bgamma_h,\bgamma_h^d)):I\rightarrow \Delta_5$ is a timelike Legendre curve with curvature 
$ 
(m _{DE}(t), n_{DE}(t))=   ( -  \dot{\psi}(t),  -m_h(t)\cosh\psi(t)- n_h(t)\sinh\psi(t) ), 
$  
where\\
$ \bnu_{DE}(\bgamma_h,\bgamma_h^d)(t) = \bgamma_h^s(t). $ 
\end{enumerate}
\end{proposition}
%%%%

%%%%
\begin{proposition}  \label{prop2}
Let $(\bgamma_h,\bgamma_h^d):I\rightarrow\Delta_1$ be a spacelike Legendre curve  with curvature $(m_h,n_h) $. 	 
Suppose	that there exists $\varphi:I\rightarrow\R$ such that $m_h(t)\cosh\varphi(t)+ n_h(t)\sinh\varphi(t)=0$ for all $t\in I$. Let $\bmu_{HE}(\bgamma_h,\bgamma_h^d)= HE(\bgamma_h,\bgamma_h^d)\wedge \bnu_{HE}(\bgamma_h,\bgamma_h^d)$ and denote the moving frame of $HE(\bgamma_h,\bgamma_h^d)$ by $ (HE, \bnu_{HE},\bmu_{HE}  )(\bgamma_h,\bgamma_h^d) $.  Then  
$$ (HE, \bnu_{HE},\bmu_{HE}  )^T(\bgamma_h,\bgamma_h^d)  =B_2(-\varphi)R( -\pi/2) ( \bgamma_h,\bgamma_h^d,\bgamma_h^s  )^T . $$
That is,  $ (HE(\bgamma_h,\bgamma_h^d)(t), \bnu_{HE}(\bgamma_h,\bgamma_h^d)(t),\bmu_{HE}(\bgamma_h,\bgamma_h^d)(t) )$ is obtained by applying a rotation with   angle $ -\pi/2$ to  $ (\bgamma_h(t),\bgamma_h^d(t),\bgamma_h^s(t) )$, followed by a boost in the $(-\bgamma_h^d)$-direction with rapidity $-\varphi(t)$  for all $t\in I$.
\end{proposition}
%%%%

Combining Propositions \ref{prop1} and  \ref{prop2}, we have the relations between the hyperbolic horocyclic and hyperbolic evolutes of $(\bgamma_h,\bgamma_h^d)$.

%%%%
\begin{theorem} \label{th2} 
Let $(\bgamma_h,\bgamma_h^d):I\rightarrow\Delta_1$ be a spacelike Legendre curve  with curvature $(m_h,n_h)$. 
\begin{enumerate}	
\item[$(1)$]   Suppose	that   there exists   $\varphi:I\rightarrow\R$ such that $m_h(t)\cosh\varphi(t)+ n_h(t)\sinh\varphi(t)=0$ for all $t\in I$.  Then 
$$  (\mathcal{HE} ^\pm ,\mathcal{DE} ^\pm,\bmu^\pm_\mathcal{E}) ^T  (\bgamma_h,\bgamma_h^d)    =   N_1^\mp(-\tanh\varphi )B_2(\varphi )(HE,\bnu_{HE},\bmu_{HE} ) ^T  (\bgamma_h,\bgamma_h^d) .  $$ 
\item[$(2)$] Suppose that there exists   $f:I\rightarrow\R  $ such that $m_h(t)+f(t)n_h(t)=0$  for all $t\in I$.  
If $|f(t)|<1$ for all $t\in I$, then $$ ( HE  ,\bnu_{HE},\bmu_{HE})^T (\bgamma_h,\bgamma_h^d) 	 =  B_2(-\operatorname{arctanh}f )N_1^\mp(f )(\mathcal{HE} ^\pm,\mathcal{DE} ^\pm,\bmu^\pm_\mathcal{E} )^T(\bgamma_h,\bgamma_h^d)   .$$
\end{enumerate}
\end{theorem}
%%%%

By Lemma \ref{lemma1}, we can discuss the singularities of hyperbolic horocyclic and hyperbolic evolutes.
%%%%%
\begin{proposition}[\cite{NTZ2026}]
Let $(\bgamma_h,\bgamma_h^d):I\rightarrow\Delta_1$ be a spacelike Legendre curve  with curvature $(m_h,n_h)$.	   Suppose that there exists   $f:I\rightarrow\R  $ such that $m_h(t)+f(t)n_h(t)=0$ for all $t\in I$.  For any $t_0\in I$, the following assertions hold.
\begin{enumerate} 
\item[$(1)$]  $\mathcal{HE}^\pm(\bgamma_h,\bgamma_h^d) $ is singular at $t_0$ if and only if  $ - 2\dot{f}(t_0)\pm f^2(t_0)n_h(t_0)=0 $.
\item[$(2)$]  $\mathcal{HE}^\pm(\bgamma_h,\bgamma_h^d) $ has a $(2,3)$-cusp at $t_0$ if and only if  $n_h(t_0)\neq0$,  $  -2\dot{f}(t_0)\pm f^2(t_0)n_h(t_0)=0 $ and $ -2\ddot{f}(t_0)\pm 2f(t_0)\dot{f}(t_0)n_h(t_0)\pm f^2(t_0)\dot{n}_h(t_0)\neq0$.
\item[$(3)$] $\mathcal{HE}^\pm(\bgamma_h,\bgamma_h^d) $ has a  $(3,4)$-cusp at $t_0$ if and only if  $n_h(t_0)\neq0$, $-2\dot{f}(t_0)\pm f^2(t_0)n_h(t_0)=-2\ddot{f}(t_0)\pm2f(t_0)\dot{f}(t_0)n_h(t_0)\pm f^2(t_0)\dot{n}_h(t_0)=0$     and  $ \pm2\dot{f^2}(t_0)n_h(t_0)\pm 2f(t_0)\ddot{f}(t_0)n_h(t_0)\pm 4f(t_0)\dot{f}(t_0)\dot{n}_h(t_0) \pm f^2(t_0)\ddot{n}_h(t_0)-2\dddot{f}(t_0)\neq0$.
\item[$(4)$] $\mathcal{HE}^\pm(\bgamma_h,\bgamma_h^d) $ has a  $(2,5)$-cusp at $t_0$ if and only if  $n_h(t_0)=0,$ $\dot{f}(t_0)=0$, $ -2\ddot{f}(t_0)\pm f^2(t_0)\dot{n}_h(t_0)\neq0 $ and $ -\ddot{f}(t_0)\ddot{n}_h(t_0)+\dddot{f}(t_0)\dot{n}_h(t_0)\neq0 $.
\item[$(5)$] $\mathcal{HE}^\pm(\bgamma_h,\bgamma_h^d) $ has a  $(3,5)$-cusp at $t_0$ if and only if   $n_h(t_0)=\dot{f}(t_0)=  -2\ddot{f}(t_0)\pm f^2(t_0)\dot{n}_h(t_0)=0 $, $\dot{n}_h(t_0)\neq0$ and $-2\dddot{f}(t_0)\pm f^2(t_0)\ddot{n}_h(t_0)\neq0$. 
\end{enumerate}
\end{proposition}
%%%%%

%%%%%
\begin{proposition}
Let $(\bgamma_h,\bgamma_h^d):I\rightarrow\Delta_1$ be a spacelike Legendre curve  with curvature $(m_h,n_h)$.	Suppose  there exists   $\varphi:I\rightarrow\R$ such that $m_h(t)\cosh\varphi(t)+ n_h(t)\sinh\varphi(t)=0$ for all $t\in I$.       For any $t_0\in I$, the following assertions hold.
\begin{enumerate} 
\item[$(1)$]  $ HE(\bgamma_h,\bgamma_h^d) $ is singular at $t_0$ if and only if $ \dot{\varphi}(t_0)=0 $.
\item[$(2)$]  $ HE(\bgamma_h,\bgamma_h^d) $ has a $(2,3)$-cusp at $t_0$ if and only if  $ \dot{\varphi}(t_0)=0 $, $ n_h(t_0)\neq0 $ and $ \ddot{\varphi}(t_0)\neq0 $.
\item[$(3)$] $ HE(\bgamma_h,\bgamma_h^d) $ has a  $(3,4)$-cusp at $t_0$ if and only if $ \dot{\varphi}(t_0)=\ddot{\varphi}(t_0)=0 $, $ n_h(t_0)\neq0 $ and $\dddot{\varphi}(t_0)\neq0 $.
\item[$(4)$] $ HE(\bgamma_h,\bgamma_h^d) $ has a  $(2,5)$-cusp at $t_0$ if and only if  $ \dot{\varphi}(t_0)=n_h(t_0)=0 $,  $\ddot{\varphi}(t_0)\neq0 $ and   $\ddot{\varphi}(t_0)\ddot{n}_h(t_0)-\dddot{\varphi}(t_0)\dot{n}_h(t_0)\neq0 $.
\item[$(5)$] $ HE(\bgamma_h,\bgamma_h^d) $ has a  $(3,5)$-cusp at $t_0$ if and only if  $ \dot{\varphi}(t_0)=\ddot{\varphi}(t_0)=n_h(t_0)=0 $,  $\dddot{\varphi}(t_0)\neq0 $ and $\dot{n}_h(t_0)\neq0 $.
\end{enumerate}
\end{proposition}
%%%%%

%%%%%
\begin{corollary}\label{cor1}
If  there exist   $f,\varphi:I\rightarrow\R  $ such that $m_h(t)+f(t)n_h(t)=0$ and $m_h(t)\cosh\varphi(t)+ n_h(t)\sinh\varphi(t)=0$ hold for all $t\in I$, then $f(t)= \tanh\varphi(t)$ for all $t\in I$ by the assumption that $\mathrm{Reg}(\bgamma_h^d)= \{t\in I\mid n_h(t)\neq0 \}$  is   dense  in $ I $. Furthermore,  for any $t_0\in I$, the following assertions hold.
\begin{enumerate} 
\item[$(1)$] Suppose that $m_h(t_0)=0 $. Then $t_0$ is a singular point of $\mathcal{HE}^\pm(\bgamma_h,\bgamma_h^d)$ if and only if $t_0$ is a singular point of $ HE(\bgamma_h,\bgamma_h^d)$.
Especially,  \\
{\rm (i)}   $t_0$ is a  $(2,3)$-cusp of $\mathcal{HE}^\pm(\bgamma_h,\bgamma_h^d)$ if and only if $t_0$ is a $(2,3)$-cusp of $ HE(\bgamma_h,\bgamma_h^d)$.\\
{\rm (ii)}   $t_0$ is a  $(3,4)$-cusp of $\mathcal{HE}^\pm(\bgamma_h,\bgamma_h^d)$ if and only if $t_0$ is a $(3,4)$-cusp of $ HE(\bgamma_h,\bgamma_h^d)$.
\item[$(2)$] Suppose  that $m_h(t_0)\neq0 $. Then,\\
{\rm (i)} If $t_0$ is a singular point of $\mathcal{HE}^\pm(\bgamma_h,\bgamma_h^d)$, then   $t_0$ is a regular point of $ HE(\bgamma_h,\bgamma_h^d)$.\\
{\rm (ii)} If $t_0$ is a singular point of $ HE(\bgamma_h,\bgamma_h^d)$, then   $t_0$ is a regular point of $\mathcal{HE}^\pm(\bgamma_h,\bgamma_h^d)$.
\item[$(3)$] Suppose that $\dot{m}_h(t_0)=0 $. Then,\\
{\rm (i)}    $t_0$ is a  $(2,5)$-cusp of $\mathcal{HE}^\pm(\bgamma_h,\bgamma_h^d)$ if and only if $t_0$ is a $(2,5)$-cusp of $ HE(\bgamma_h,\bgamma_h^d)$.\\
{\rm (ii)}   $t_0$ is a  $(3,5)$-cusp of $\mathcal{HE}^\pm(\bgamma_h,\bgamma_h^d)$ if and only if $t_0$ is a $(3,5)$-cusp of $ HE(\bgamma_h,\bgamma_h^d)$. 
\end{enumerate}	 
\end{corollary}
%%%%%

%%%%%%
\begin{example}\rm
Let $(\bgamma_h,\bgamma_h^d):(\R,0)\rightarrow\Delta_1$   be a spacelike Legendre curve given by
$$\begin{aligned}
\bgamma_h(t) =&~
\Big( \cosh(t^3+1),  \sinh(t^3+1)\cos t^4,\sinh(t^3+1)\sin t^4 
\Big),\\
\bgamma_h^d(t) =&~
\frac{1}{\sqrt{16t^2\sinh^2(t^3+1)+9}}
\Big(
4t\sinh^2(t^3+1),
4t\sinh(t^3+1)\cosh(t^3+1)\cos t^4+3\sin t^4,\\
&4t\sinh(t^3+1)\cosh(t^3+1)\sin t^4-3\cos t^4
\Big). \end{aligned} $$
%%%%%
\begin{figure}[h] 
	\includegraphics[width = 14.5cm]{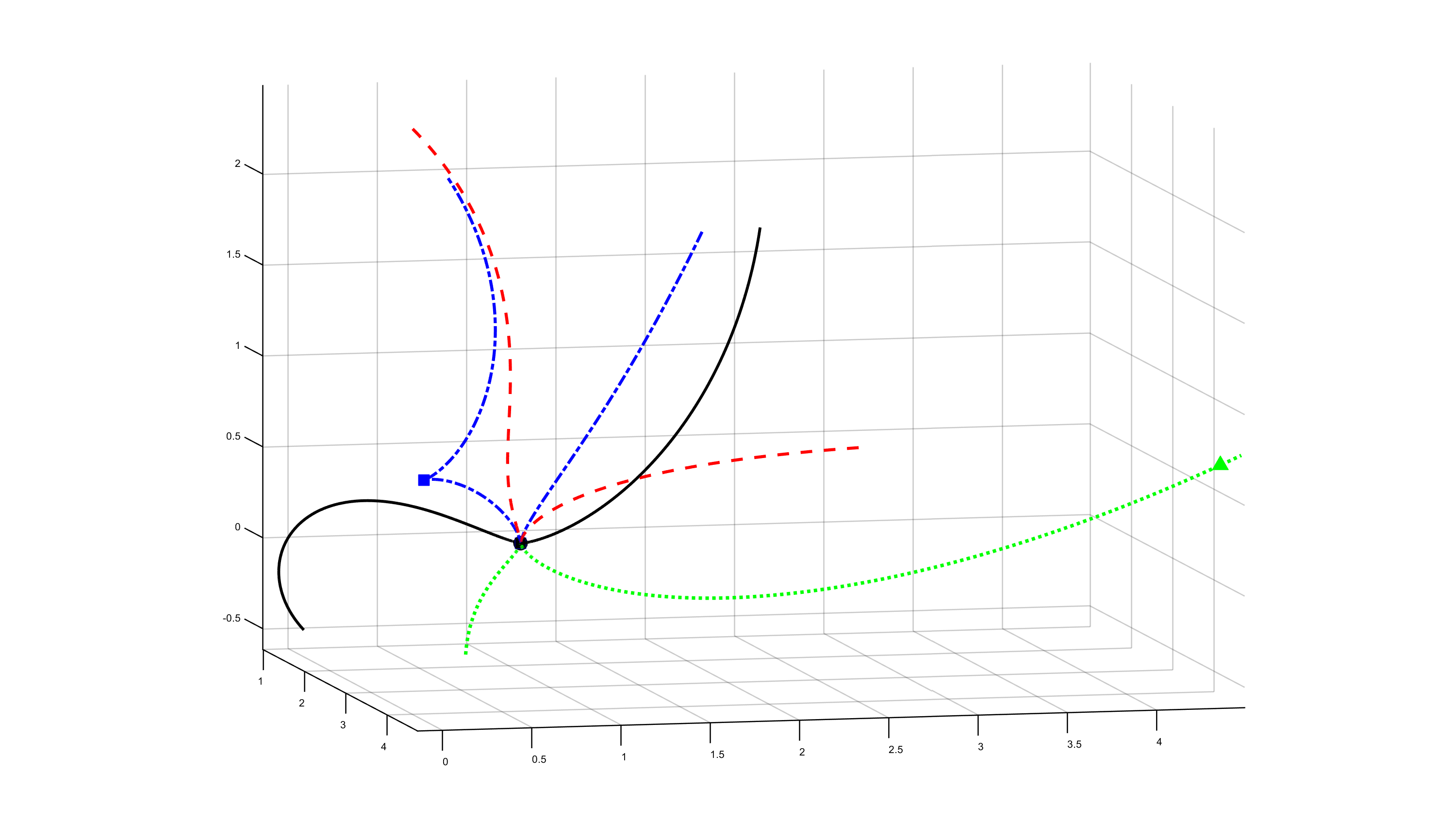}
	\centering
	\caption{The hyperbolic  horocyclic evolutes $ \mathcal{HE}^+(\bgamma_h,\bgamma_h^d) $ (red dashed curve), $\mathcal{HE}^-(\bgamma_h,\bgamma_h^d)$ (blue dash-dotted curve) and  hyperbolic   evolute $HE  (\bgamma_h,\bgamma_h^d) $ (green dotted curve) of $(\bgamma_h,\bgamma_h^d)$. The black solid curve is the spacelike frontal $\bgamma_h$.}
	\label{figure1}
\end{figure}
%%%%%
Then 
$$\begin{aligned}
\bgamma_h^s(t) =&~\bgamma_h(t)\wedge\bgamma_h^d(t)=
\frac{1}{\sqrt{16t^2\sinh^2(t^3+1)+9}}
\Big(
3\sinh(t^3+1),\\
&3\cos t^4\cosh(t^3+1)-4t\sin t^4\sinh(t^3+1),
3\sin t^4\cosh(t^3+1)+4t\cos t^4\sinh(t^3+1)
\Big)\end{aligned}$$ and the curvature $(m_h,n_h)$ of $(\bgamma_h,\bgamma_h^d) $ is
$$ \left(t^2\sqrt{16t^2\sinh^2(t^3+1)+9},\;
\frac{4\left(t^3\big(16t^2\sinh^2(t^3+1)+18\big)\cosh(t^3+1)+3\sinh(t^3+1)\right)}{16t^2\sinh^2(t^3+1)+9}\right).$$
There exists a unique smooth function $f:(\R,0)\rightarrow\R$, $$f(t)=\frac{-t^2(16t^2\sinh^2(t^3+1)+9)^{\frac{3}{2}}}{4\left(t^3\big(16t^2\sinh^2(t^3+1)+18\big)\cosh(t^3+1)+3\sinh(t^3+1)\right)}  $$ such that $f(t)n_h(t)+m_h(t)=0$ for all $t\in(\R,0)$. By Proposition \ref{prop1}, we have the hyperbolic horocyclic evolute  of $(\bgamma_h,\bgamma_h^d)$ satisfying
$$ 
(\mathcal{HE}^\pm, \mathcal{DE}^\pm,\bmu^\pm_\mathcal{E}  )^T(\bgamma_h,\bgamma_h^d) =N_1^{\mp}(-f )R(- \pi/2) (\bgamma_h, \bgamma_h^d,\bgamma_h^s  ) ^T . $$
It can be proved that $|f(t)|<1$ for all $t\in (\R,0)$. Then  by Theorem \ref{th2}, the hyperbolic  evolute  of $(\bgamma_h,\bgamma_h^d)$ satisfies $$ ( HE  ,\bnu_{HE},\bmu_{HE})^T  (\bgamma_h,\bgamma_h^d)	 =  B_2(-\operatorname{arctanh}f )N_1^\mp(f )(\mathcal{HE} ^\pm  ,\mathcal{DE} ^\pm  ,\bmu^\pm_\mathcal{E} )^T (\bgamma_h,\bgamma_h^d) .$$

Since $m_h(0)=\dot{m}_h(0)=0$, $\ddot{m}_h(0)\neq0$ and $n_h(0)\neq0$, $\bgamma_h$ has a $(3,4)$-cusp at $t_1=0$ by Lemma  \ref{lemma1}. It follows that $\mathcal{HE}^\pm(\bgamma_h,\bgamma_h^d)$ has a  $(2,3)$-cusp at $t_1=0$. Furthermore,  $t_1=0$ is also a $(2,3)$-cusp of $HE  (\bgamma_h,\bgamma_h^d) $ by Corollary \ref{cor1} (1)(i) (black circle point in Figure \ref{figure1}). 

$\mathcal{HE}^-(\bgamma_h,\bgamma_h^d)$ has a singular point at $t_2$, where $t_2\approx 0.72$ (blue square point in Figure \ref{figure1}). Since  $m_h(t_2)\neq0$,  $t_2$ is a regular point of $HE  (\bgamma_h,\bgamma_h^d) $ by Corollary \ref{cor1} (2)(i) (green triangle point in Figure \ref{figure1}). 
\end{example}
%%%%%%

\subsection{Involutes of spacelike Legendre curves}

%%%%%
\begin{definition}[Hyperbolic  horocyclic involute  \cite{NTZ2026}]\rm
Let $(\bgamma_h,\bgamma_h^d):I\rightarrow\Delta_1$ be a spacelike Legendre curve  with curvature $(m_h,n_h)$. The \textit{hyperbolic horocyclic involute}  $ \mathcal{HI}_{s_\pm}^\pm(\bgamma_h,\bgamma_h^d):I\rightarrow H^2$ of $(\bgamma_h,\bgamma_h^d)$  is defined as 
$$
\mathcal{HI}_{s_\pm}^\pm(\bgamma_h,\bgamma_h^d)(t)=\bgamma_h(t)+s_\pm(t)\bgamma_h^s(t)+\frac{s_\pm^2(t)}{2}(\bgamma_h(t)\pm\bgamma_h^d(t)), 
$$
where $s_\pm(t)$ is a solution of the Riccati equation
\begin{equation}\label{Riccati1}
\frac{{\rm d}s_\pm}{{\rm d}t}(t)= \frac{m_h(t)\pm n_h(t)}{2}s_\pm^2(t)-m_h(t).\end{equation}
\end{definition}
%%%%%

%%%%%
\begin{remark} 
\rm The equation  \eqref{Riccati1}  admits  a unique local solution. That is, for any $t_0\in I$, there exists an interval $\overline{I} \subset I$ containing $t_0$
on which a solution exists and is unique. If a global solution does not exist, we examine the behavior of curves on such an  interval   $\overline{I}$. A similar property holds for equations \eqref{Bernoulli1}, \eqref{Riccati2} and \eqref{Bernoulli2}.
\end{remark}
%%%%%

%%%%
\begin{proposition}[\cite{NTZ2026}] 
Let $\mathcal{HI}_{s_\pm}^\pm(\bgamma_h,\bgamma_h^d)$   be a  hyperbolic horocyclic involute  of $(\bgamma_h,\bgamma_h^d)$. Then  
$(\mathcal{HI}_{s_\pm}^\pm(\bgamma_h,\bgamma_h^d), \bnu_\mathcal{HI}^\pm(\bgamma_h,\bgamma_h^d) ):I\rightarrow \Delta_1$ is a spacelike Legendre curve with curvature  
$ (m^\pm_\mathcal{HI}(t), n^\pm_\mathcal{HI}(t))= (- s_\pm(t)(n_h(t)\pm m_h(t)),n_h(t)\pm m_h(t)), 
$  where $\bnu_\mathcal{HI}^\pm(\bgamma_h,\bgamma_h^d)(t)=-s_\pm(t)\bgamma_h(t)-\bgamma_h^s(t)\mp s_\pm(t)\bgamma_h^d(t). $
\end{proposition}
%%%%

%%%%
\begin{proposition} \label{prop3}
Let  $  \bmu^\pm_{\mathcal{HI}}(\bgamma_h,\bgamma_h^d)=\mathcal{HI}_{s_\pm}^\pm(\bgamma_h,\bgamma_h^d) \wedge \bnu_\mathcal{HI}^\pm(\bgamma_h,\bgamma_h^d)  $  and denote the moving frame of $\mathcal{HI}_{s_\pm}^\pm(\bgamma_h,\bgamma_h^d)$ by $( \mathcal{HI}_{s_\pm}^\pm,\bnu_\mathcal{HI}^\pm ,   \bmu^\pm_\mathcal{HI}  )(\bgamma_h,\bgamma_h^d)  $.
Then  
$$ ( \mathcal{HI}_{s_\pm}^\pm,\bnu_\mathcal{HI}^\pm,   \bmu^\pm_\mathcal{HI}  )^T(\bgamma_h,\bgamma_h^d) =N_2^\mp(-s_\pm )R(\pi/2) ( \bgamma_h,   \bgamma_h^d,\bgamma_h^s )^T . $$
That is,  $ ( \mathcal{HI}_{s_\pm}^\pm(\bgamma_h,\bgamma_h^d)(t),\bnu_\mathcal{HI}^\pm(\bgamma_h,\bgamma_h^d)(t),   \bmu^\pm_\mathcal{HI}(\bgamma_h,\bgamma_h^d)(t)  )$ is obtained by applying a rotation with angle $  \pi/2$ to $ (\bgamma_h(t), \bgamma_h^d(t),\bgamma_h^s(t)  )$, followed by a null rotation   around the $(\bgamma_h\mp \bgamma_h^d)$-direction with angle $-s_\pm(t)$ for all $t\in I$. 
\end{proposition}
%%%%

%%%%%
\begin{definition}[Hyperbolic and de Sitter involutes \cite{ZP2026}]\rm
Let $(\bgamma_h,\bgamma_h^d):I\rightarrow\Delta_1$ be a spacelike Legendre curve  with curvature $(m_h,n_h)$ and  $t_0\in I$.
Then the \textit{hyperbolic involute}  $ HI_{t_0}(\bgamma_h,\bgamma_h^d):I\rightarrow H^2$  of    $(\bgamma_h,\bgamma_h^d)$  is defined as  
\begin{equation}\notag
HI_{t_0}(\bgamma_h,\bgamma_h^d)(t)=  \cosh\left(\int_{t_0}^t{m_h(t)}{\rm d}t\right)\bgamma_h(t)-\sinh\left(\int_{t_0}^t{m_h(t)}{\rm d}t\right)\bgamma_h^s(t).
\end{equation} 
The \textit{de Sitter involute} $ DI_{t_0}(\bgamma_h,\bgamma_h^d):I\rightarrow S_1^2$   of $(\bgamma_h,\bgamma_h^d)$  is defined as  
\begin{equation}\notag
DI_{t_0}(\bgamma_h,\bgamma_h^d)(t)=  \sinh\left(\int_{t_0}^t{m_h(t)}{\rm d}t\right)\bgamma_h(t)-\cosh\left(\int_{t_0}^t{m_h(t)}{\rm d}t\right)\bgamma_h^s(t).
\end{equation}  
\end{definition}
%%%%%

%%%%%
\begin{proposition}[\cite{ZP2026}]  
Let $ HI_{t_0}(\bgamma_h,\bgamma_h^d) $ and $DI_{t_0}(\bgamma_h,\bgamma_h^d)$ be  hyperbolic and de Sitter involutes of  $(\bgamma_h,\bgamma_h^d)$, respectively. Then
\begin{enumerate} 
\item[$(1)$]   
$ ( HI_{t_0}(\bgamma_h,\bgamma_h^d), DI_{t_0}(\bgamma_h,\bgamma_h^d)  ):I\rightarrow \Delta_1$ is a spacelike Legendre curve with curvature \\
$ (m _{HI}(t), n_{HI}(t) )=  \left(     n_h(t)\sinh (\int_{t_0}^t{m_h(t)}{\rm d}t ), n_h(t)\cosh(\int_{t_0}^t{m_h(t)}{\rm d}t) \right). 
$  
\item[$(2)$]    
$ ( DI_{t_0}(\bgamma_h,\bgamma_h^d), HI_{t_0}(\bgamma_h,\bgamma_h^d)   ):I\rightarrow \overline{\Delta}_1$ is a spacelike Legendre curve with curvature \\
$ (m _{DI}(t), n_{DI}(t) )= \left(      -n_h(t)\cosh (\int_{t_0}^t{m_h(t)}{\rm d}t ), -n_h(t)\sinh (\int_{t_0}^t{m_h(t)}{\rm d}t )\right).
$  
\end{enumerate}
\end{proposition}
%%%%

%%%%
\begin{proposition}  \label{prop4}
  Let $\bmu_I(\bgamma_h,\bgamma_h^d)=HI_{t_0}(\bgamma_h,\bgamma_h^d)\wedge DI_{t_0}(\bgamma_h,\bgamma_h^d)$ and denote the moving frame of $HI_{t_0}(\bgamma_h,\bgamma_h^d)$ by $ 
( HI_{t_0} , DI_{t_0},\bmu_I )  (\bgamma_h,\bgamma_h^d)     $.	Then  
$$ 
( HI_{t_0} , DI_{t_0},\bmu_I)^T (\bgamma_h,\bgamma_h^d)    =	B_1\left( \int_{t_0}^t{m_h }{\rm d}t\right)R(\pi/2) \left(   \bgamma_h,\bgamma_h^d, \bgamma_h^s  \right)^T. $$
That is,	$ ( HI_{t_0}(\bgamma_h,\bgamma_h^d)(t), DI_{t_0}(\bgamma_h,\bgamma_h^d)(t),\bmu_I(\bgamma_h,\bgamma_h^d)  (t)  )$  is obtained by applying a rotation  with angle $ \pi/2$ to $ (\bgamma_h(t), \bgamma_h^d(t),\bgamma_h^s(t)  )$, followed by a boost   in the $(-\bgamma_h^s)$-direction with rapidity $ \int_{t_0}^t{m_h(t)}{\rm d}t$ for all $t\in I$.
\end{proposition}
%%%%

Combining Propositions \ref{prop3} and   \ref{prop4}, we have the relations between  the hyperbolic horocyclic   and hyperbolic  involutes of $(\bgamma_h,\bgamma_h^d)$.

%%%%
\begin{theorem} 
Let $(\bgamma_h,\bgamma_h^d):I\rightarrow\Delta_1$ be a spacelike Legendre curve  with curvature $(m_h,n_h)$.  	     $  \mathcal{HI}_{s_\pm}^\pm(\bgamma_h,\bgamma_h^d)   $ and
$ HI_{t_0}(\bgamma_h,\bgamma_h^d)   $  are hyperbolic horocyclic and hyperbolic involutes of $(\bgamma_h,\bgamma_h^d)$, respectively. Then  
$$\begin{aligned}  ( \mathcal{HI}_{s_\pm}^\pm,\bnu_\mathcal{HI}^\pm,\bmu^\pm_\mathcal{HI}    )^T(\bgamma_h,\bgamma_h^d) &=  N_2^\mp (-s_\pm)  B_1\left(-\int_{t_0}^t{m_h }{\rm d}t\right)\left(  HI_{t_0}, DI_{t_0},\bmu_I  \right)^T   (\bgamma_h,\bgamma_h^d),\\
(  HI_{t_0}, DI_{t_0},\bmu_I  )^T (\bgamma_h,\bgamma_h^d)  &=  B_1\left( \int_{t_0}^t{m_h }{\rm d}t\right) N_2^\mp( s_\pm)( \mathcal{HI}_{s_\pm}^\pm,\bnu_\mathcal{HI}^\pm,\bmu^\pm_\mathcal{HI}    )^T (\bgamma_h,\bgamma_h^d)  .\end{aligned}$$ 
\end{theorem}
%%%%

We  discuss the singularities of hyperbolic horocyclic and hyperbolic involutes.
%%%%%
\begin{proposition}[\cite{NTZ2026}]  
Let $(\bgamma_h,\bgamma_h^d):I\rightarrow\Delta_1$ be a spacelike Legendre curve  with curvature $(m_h,n_h)$.	   $\mathcal{HI}_{s_\pm}^\pm(\bgamma_h,\bgamma_h^d) $ is a hyperbolic horocyclic involute of $ (\bgamma_h,\bgamma_h^d)$ and $t_1\in I$.
\begin{enumerate} 
\item[$(1)$] $\mathcal{HI}_{s_\pm}^\pm(\bgamma_h,\bgamma_h^d)$ is singular at $t_1$ if and only if $ s_\pm(t_1)(n_h(t_1)\pm m_h(t_1))=0$.
\item[$(2)$] $\mathcal{HI}_{s_\pm}^\pm(\bgamma_h,\bgamma_h^d)$ has a $(2,3)$-cusp at $t_1$ if and only if $ s_\pm(t_1)=0$, $m_h(t_1)\neq0$, $  n_h(t_1)\pm m_h(t_1) \neq0$.
\item[$(3)$] $\mathcal{HI}_{s_\pm}^\pm(\bgamma_h,\bgamma_h^d)$ has a $(3,4)$-cusp at $t_1$ if and only if  $ s_\pm(t_1)=0$, $m_h(t_1)=0$, $n_h(t_1)\neq0$ and $\dot{m}_h(t_1)\neq0$.
\item[$(4)$] $\mathcal{HI}_{s_\pm}^\pm(\bgamma_h,\bgamma_h^d)$ has a $(2,5)$-cusp at $t_1$ if and only if $n_h(t_1)\pm m_h(t_1) =0$, $\dot{n}_h(t_1)\pm \dot{m}_h(t_1) \neq0$, $ s_\pm(t_1)\neq0$ and $m_h(t_1)\neq0$.
\item[$(5)$] $\mathcal{HI}_{s_\pm}^\pm(\bgamma_h,\bgamma_h^d)$ has a $(3,5)$-cusp at $t_1$ if and only if $ s_\pm(t_1)=n_h(t_1)\pm m_h(t_1) =0$, $ \dot{n}_h(t_1)\pm \dot{m}_h(t_1) \neq0$ and $m_h(t_1)\neq0$.
\end{enumerate}
\end{proposition}
%%%%%

%%%%%
\begin{proposition}
Let $(\bgamma_h,\bgamma_h^d):I\rightarrow\Delta_1$ be a spacelike Legendre curve  with curvature $(m_h,n_h)$.       For any $t_1\in I$, the following assertions hold.
\begin{enumerate} 
\item[$(1)$]  $ HI_{t_0}(\bgamma_h,\bgamma_h^d) $ is singular at $t_1$ if and only if $ n_h(t_1)\int_{t_0}^{t_1}{m_h(t)}{\rm d}t =0$.
\item[$(2)$]  $ HI_{t_0}(\bgamma_h,\bgamma_h^d) $ has a $(2,3)$-cusp at $t_1$ if and only if  $ \int_{t_0}^{t_1}{m_h(t)}{\rm d}t =0$, $ n_h(t_1)\neq0 $,  $ m_h(t_1)\neq0 $.
\item[$(3)$] $ HI_{t_0}(\bgamma_h,\bgamma_h^d) $ has a  $(3,4)$-cusp at $t_1$ if and only if $ \int_{t_0}^{t_1}{m_h(t)}{\rm d}t =m_h(t_1)=0$, $ n_h(t_1)\neq0 $ and $\dot{m}_h(t_1)\neq0 $.
\item[$(4)$] $ HI_{t_0}(\bgamma_h,\bgamma_h^d) $ has a  $(2,5)$-cusp at $t_1$ if and only if  $  n_h(t_1)=0 $,  $\dot{n}_h(t_1)\int_{t_0}^{t_1}{m_h(t)}{\rm d}t\neq0 $ and   $m_h(t_1) \neq0 $.
\item[$(5)$] $ HI_{t_0}(\bgamma_h,\bgamma_h^d) $ has a  $(3,5)$-cusp at $t_1$ if and only if  $  n_h(t_1)=\int_{t_0}^{t_1}{m_h(t)}{\rm d}t=0 $,  $\dot{n}_h(t_1)\neq0 $ and $m_h(t_1)\neq0 $.
\end{enumerate}
\end{proposition}
%%%%%

%%%%%
\begin{corollary}\label{cor3}
Let $  \mathcal{HI}_{s_\pm}^\pm(\bgamma_h,\bgamma_h^d)   $ and
$ HI_{t_0}(\bgamma_h,\bgamma_h^d)   $  be hyperbolic horocyclic and hyperbolic involutes of $(\bgamma_h,\bgamma_h^d)$, respectively. For any $t_1\in I$,  suppose that $s_\pm(t_1)=\int_{t_0}^{t_1}{m_h(t)}{\rm d}t$.  Then
\begin{enumerate} 
\item[$(1)$]  $t_1$ is a $(3,4)$-cusp of $\mathcal{HI}_{s_\pm}^\pm(\bgamma_h,\bgamma_h^d)$	 if and only if $t_1$ is a $(3,4)$-cusp of $HI_{t_0}(\bgamma_h,\bgamma_h^d)$.
\item[$(2)$] If $t_1$ is a $(2,5)$-cusp of $\mathcal{HI}_{s_\pm}^\pm(\bgamma_h,\bgamma_h^d)$, then $t_1$ is a  regular point of $HI_{t_0}(\bgamma_h,\bgamma_h^d)$.
\item[$(3)$] If $t_1$ is a $(3,5)$-cusp of $\mathcal{HI}_{s_\pm}^\pm(\bgamma_h,\bgamma_h^d)$, then $t_1$ is a  $(2,3)$-cusp of $HI_{t_0}(\bgamma_h,\bgamma_h^d)$.
\item[$(4)$] If $t_1$ is a $(2,5)$-cusp of $HI_{t_0}(\bgamma_h,\bgamma_h^d)$, then $t_1$ is a  regular point of $\mathcal{HI}_{s_\pm}^\pm(\bgamma_h,\bgamma_h^d)$.
\item[$(5)$] If $t_1$ is a $(3,5)$-cusp of $HI_{t_0}(\bgamma_h,\bgamma_h^d)$, then $t_1$ is a  $(2,3)$-cusp of $\mathcal{HI}_{s_\pm}^\pm(\bgamma_h,\bgamma_h^d)$.	 
\end{enumerate}
\end{corollary}
%%%%%

%%%%%
\begin{example}\rm
Let $(\bgamma_h,\bgamma_h^d ):\R\rightarrow\Delta_1$   be a spacelike Legendre curve with curvature $$(m_h(t),n_h(t))=(ct^2,{(2\sinh t+ct^2-ct^2\sinh^2t)}/{\cosh^2t}),$$
where $c=-\sinh1$ and $\bgamma_h^s(t)=\bgamma_h(t) \wedge\bgamma_h^d(t)$. Consider the Riccati equation $$\frac{{\rm d}s_+}{{\rm d}t}(t)= \frac{m_h(t)+ n_h(t)}{2}s_+^2(t)-m_h(t)=\frac{\sinh t+ct^2}{\cosh^2t}s_+^2(t)-ct^2.$$ 
Since $s_+(t)=\cosh t$ is a solution of the above equation, by Propositions \ref{prop3} and \ref{prop4}, the corresponding hyperbolic horocyclic involute is  	$$\begin{aligned} 
( \mathcal{HI}_{s_+}^+ ,\bnu_{\mathcal{HI}}^+  ,   \bmu^+_{\mathcal{HI}}  )^T(\bgamma_h,\bgamma_h^d) = N_2^-(-\cosh t)R(\pi/2)( \bgamma_h  ,   \bgamma_h^d  ,\bgamma_h^s   )^T   \end{aligned}$$ 
and the   hyperbolic   involute  at $t_0=\sqrt[3]{3\coth 1 + 1}$ is $$\begin{aligned}  ( HI_{t_0}, DI_{t_0} ,\bmu_I     )^T(\bgamma_h,\bgamma_h^d) = B_1\left( c(t^3-3\coth 1 - 1)/3\right)R(\pi/2) (   \bgamma_h ,\bgamma_h^d , \bgamma_h^s   )^T.\end{aligned}$$

Since $n_h(1)+ m_h(1) =0$, $\dot{n}_h(1)+ \dot{m}_h(1)=2(\cosh 1-2\sinh 1)/\cosh ^21\neq0$, $ s_+(1)=\cosh 1\neq0$ and $m_h(1)=c=-\sinh1\neq0$,   $\mathcal{HI}_{s_+}^+(\bgamma_h,\bgamma_h^d)$ has a $(2,5)$-cusp at $t_1=1$. Moreover, noting that $s_+(1)=\cosh 1=\int_{t_0}^{1}{ct^2}{\rm d}t$, it follows from Corollary \ref{cor3} (2) that $t_1$ is a regular point of $HI_{t_0}(\bgamma_h,\bgamma_h^d)$.
\end{example}
%%%%%

\subsection{Parallels of spacelike Legendre curves  }

%%%%%
\begin{definition}[Hyperbolic horocyclic parallel  \cite{NTZ2026}]  \rm
Let $(\bgamma_h,\bgamma_h^d):I\rightarrow\Delta_1$ be a spacelike Legendre curve with curvature $(m_h,n_h)$. 
The \textit{hyperbolic horocyclic parallel} $\mathcal{HP}_{\lambda_\pm}^\pm (\bgamma_h,\bgamma_h^d):I\rightarrow H^2$   of $(\bgamma_h,\bgamma_h^d)$  is   defined as
$$
\mathcal{HP}_{\lambda_\pm}^\pm (\bgamma_h,\bgamma_h^d)(t)=\bgamma_h(t)+\lambda _\pm(t)\bgamma_h^d(t)+\frac{\lambda_\pm^2(t)}{2}(\bgamma_h(t)\pm\bgamma_h^s(t)), 
$$
where 
$ \lambda_\pm(t)  $ is a solution of the Bernoulli equation \begin{equation}\label{Bernoulli1}
\frac{{\rm d}\lambda _\pm}{{\rm d}t}(t)=\mp \left(\frac{ \lambda _\pm^2(t) n_h(t)}{2}+\lambda _\pm(t)m_h(t)\right).
\end{equation}
\end{definition}
%%%%%

%%%%%
\begin{proposition}[\cite{NTZ2026}]
Let  $\mathcal{HP}_{\lambda_\pm}^\pm (\bgamma_h,\bgamma_h^d)$   be a  hyperbolic  horocyclic parallel of $(\bgamma_h,\bgamma_h^d)$. Then  
$(\mathcal{HP}_{\lambda_\pm}^\pm (\bgamma_h,\bgamma_h^d),\bnu_\mathcal{HP}^\pm (\bgamma_h,\bgamma_h^d)):I\rightarrow\Delta_1$ is a spacelike Legendre curve with curvature
$  \left(m_\mathcal{HP}^\pm(t),n_\mathcal{HP}^\pm(t)\right)=\left(\lambda _\pm(t)n_h(t)+m_h(t),n_h(t)\right), $  where $\bnu_\mathcal{HP}^\pm (\bgamma_h,\bgamma_h^d)(t)=\lambda _\pm(t)\bgamma_h(t)+\bgamma_h^d(t)\pm \lambda _\pm(t)\bgamma_h^s(t).$  
\end{proposition}
%%%%%

%%%%
\begin{proposition} \label{prop5}
  Let  $\bmu_{\mathcal{HP}}^\pm(\bgamma_h,\bgamma_h^d)=\mathcal{HP}_{\lambda_\pm}^\pm (\bgamma_h,\bgamma_h^d)\wedge \bnu_\mathcal{HP}^\pm (\bgamma_h,\bgamma_h^d)$ and denote the moving frame of $\mathcal{HP}_{\lambda_\pm}^\pm (\bgamma_h,\bgamma_h^d)$ by $ (	\mathcal{HP}_{\lambda_\pm}^\pm, \bnu_\mathcal{HP}^\pm ,	 \bmu_{\mathcal{HP}}^\pm    )(\bgamma_h,\bgamma_h^d)$. Then   
$$
(	\mathcal{HP}_{\lambda_\pm}^\pm, \bnu_\mathcal{HP}^\pm ,	 \bmu_{\mathcal{HP}}^\pm   )^T(\bgamma_h,\bgamma_h^d)
=N_2^\mp( \lambda_\pm )(	\bgamma_h, \bgamma_h^d,\bgamma_h^s ) ^T.$$
That is,	$ ( 	\mathcal{HP}_{\lambda_\pm}^\pm (\bgamma_h,\bgamma_h^d)(t), \bnu_\mathcal{HP}^\pm (\bgamma_h,\bgamma_h^d)(t),	 \bmu_{\mathcal{HP}}^\pm(\bgamma_h,\bgamma_h^d)(t)   )$
is obtained by applying a null rotation  around the $(\bgamma_h\mp \bgamma_h^s)$-direction with   angle $ \lambda_\pm(t)$ to $ ( \bgamma_h(t), \bgamma_h^d(t),\bgamma_h^s(t)   )$ for all $t\in I$.  
\end{proposition}
%%%%

%%%%%
\begin{definition}[Hyperbolic and de Sitter parallels \cite{IT2007}]  \rm
Let $(\bgamma_h,\bgamma_h^d):I\rightarrow\Delta_1$ be a spacelike Legendre curve with curvature $(m_h,n_h)$ and $\theta\in \R$. 
The \textit{hyperbolic parallel}\\ $HP_{\theta}(\bgamma_h,\bgamma_h^d)  :I\rightarrow H^2$    of $(\bgamma_h,\bgamma_h^d)$ is   defined as
$$
HP_{\theta}(\bgamma_h,\bgamma_h^d)(t)=\cosh \theta\bgamma_h(t)+\sinh\theta\bgamma_h^d(t).$$
The \textit{de Sitter parallel} $DP_{\theta}(\bgamma_h,\bgamma_h^d)  :I\rightarrow S_1^2$  of $(\bgamma_h,\bgamma_h^d)$  is   defined as
$$
DP_{\theta}(\bgamma_h,\bgamma_h^d)(t)=\sinh \theta\bgamma_h(t)+\cosh\theta\bgamma_h^d(t).$$
\end{definition}
%%%%%

%%%%
\begin{proposition}\label{prop0} 
Let $ HP_{\theta}(\bgamma_h,\bgamma_h^d)$ and $ DP_{\theta}(\bgamma_h,\bgamma_h^d)$ be   hyperbolic and de Sitter parallels of $(\bgamma_h,\bgamma_h^d)$, respectively.  Then   \begin{enumerate} 
\item[$(1)$] 
$( HP_{\theta}(\bgamma_h,\bgamma_h^d),  DP_{\theta}(\bgamma_h,\bgamma_h^d) ):I\rightarrow \Delta_1$ is a spacelike Legendre curve with curvature\\ 
$ 
(m _{HP}(t), n_{HP}(t))= (      m_h(t)\cosh \theta+ n_h(t)\sinh\theta,   n_h(t)\cosh \theta+m_h(t)\sinh\theta ).
$

\item[$(2)$]   
$( DP_{\theta}(\bgamma_h,\bgamma_h^d), HP_{\theta}(\bgamma_h,\bgamma_h^d) ):I\rightarrow \overline{\Delta}_1$ is a spacelike Legendre curve with curvature \\
$ (m _{DP}(t), n_{DP}(t))= (    - n_h(t)\cosh \theta- m_h(t)\sinh\theta, -  m_h(t)\cosh \theta- n_h(t)\sinh\theta).
$  
\end{enumerate}
\end{proposition}
%%%%

%%%%
\begin{proposition} \label{prop6}
  Let $\bmu_P(\bgamma_h,\bgamma_h^d)=HP_{\theta}(\bgamma_h,\bgamma_h^d)\wedge DP_{\theta}(\bgamma_h,\bgamma_h^d)$ and denote the moving frame of $ HP_{\theta}(\bgamma_h,\bgamma_h^d)$ by $( 	HP_{\theta},  DP_{\theta},	\bmu_P    ) (\bgamma_h,\bgamma_h^d)$.  Then      
$$
( 	HP_{\theta},  DP_{\theta},	\bmu_P    )^T(\bgamma_h,\bgamma_h^d)
=B_1(\theta) ( \bgamma_h, \bgamma_h^d,\bgamma_h^s )^T .$$
That is,	$( 	HP_{\theta}(\bgamma_h,\bgamma_h^d)(t),  DP_{\theta}(\bgamma_h,\bgamma_h^d)(t),	\bmu_P(\bgamma_h,\bgamma_h^d)(t) )    $ is obtained by applying a  boost in the $\bgamma_h^d$-direction  with rapidity $\theta$ to $( \bgamma_h(t), \bgamma_h^d(t),\bgamma_h^s(t) )  $ for all $t\in I$.   
\end{proposition}
%%%%

Combining Propositions \ref{prop5} and  \ref{prop6}, we have  the relations between  the hyperbolic horocyclic parallels  and hyperbolic   parallels    of $(\bgamma_h,\bgamma_h^d)$.

%%%%
\begin{theorem}
Let $(\bgamma_h,\bgamma_h^d):I\rightarrow\Delta_1$ be a spacelike Legendre curve  with curvature $(m_h,n_h)$.  	      $  \mathcal{HP}_{\lambda_\pm}^\pm (\bgamma_h,\bgamma_h^d)    $ and
$ HP_\theta(\bgamma_h,\bgamma_h^d)  $  are  hyperbolic  horocyclic and  hyperbolic parallels of $(\bgamma_h,\bgamma_h^d) $, respectively. Then 
$$\begin{aligned}
(	\mathcal{HP}_{\lambda_\pm}^\pm , \bnu^\pm_{\mathcal{HP}}  ,	 \bmu^\pm_{\mathcal{HP}} )^T (\bgamma_h,\bgamma_h^d)& =N_2^\mp( \lambda_\pm )B_1(-\theta)( 	HP_{\theta} , DP_{\theta},	\bmu_P )  ^T(\bgamma_h,\bgamma_h^d) ,\\
( 	HP_{\theta} ,  DP_{\theta} ,	\bmu_P  )^T (\bgamma_h,\bgamma_h^d)  &=B_1( \theta) N_2^\mp(- \lambda_\pm ) (	\mathcal{HP}_{\lambda_\pm}^\pm  , \bnu^\pm_{\mathcal{HP}}  ,	 \bmu^\pm_{\mathcal{HP}}   )^T (\bgamma_h,\bgamma_h^d)    
.\end{aligned}$$ 
\end{theorem}

\begin{figure}[!h]	
\centering

\begin{tikzpicture}[
>=stealth, % 
node distance=3cm,
vertex/.style={circle, minimum size=6mm, inner sep=0pt, font= \small}
]
\node[vertex] (D1) at (-4.4, 2.8) {$(\mathcal{HE}^\pm, \mathcal{DE}^\pm ,\bmu^\pm_\mathcal{E}  )(\bgamma_h,\bgamma_h^d)$};
\node[vertex] (D2) at (4.4, 2.8) {$ (HE, \bnu_{HE},\bmu_{HE}   ) (\bgamma_h,\bgamma_h^d)$};
\node[vertex] (D3) at (-4.4, -2.8) {$( \mathcal{HI}_{s_\pm}^\pm, \bnu_\mathcal{HI}^\pm ,\bmu^\pm_\mathcal{HI} )(\bgamma_h,\bgamma_h^d)$};
\node[vertex] (D4) at (4.4, -2.8) {$ ( HI_{t_0}, DI_{t_0} ,\bmu_I  ) (\bgamma_h,\bgamma_h^d)   $};
\node[vertex] (D5) at (0,1.2) {$(\bgamma_h,\bgamma_h^s,-\bgamma_h^d )$};
\node[vertex] (D0) at (0,0) {$( \bgamma_h,\bgamma_h^d,\bgamma_h^s    )$};
\node[vertex] (D55) at (0,-1.2) {$( \bgamma_h,-\bgamma_h^s,\bgamma_h^d  )$};
\node[vertex] (D11) at (-4.4, 6.2) {$(  \mathcal{HP}_{\lambda_\pm}^\pm  , \bnu^\pm_{\mathcal{ HP}}  ,\bmu^\pm_{\mathcal{HP}}  )(\bgamma_h,\bgamma_h^d)   $};
\node[vertex] (D22) at (4.4, 6.2) {$ (HP_{\theta},  {DP}_{\theta},\bmu_P   ) (\bgamma_h,\bgamma_h^d)$};
\node[vertex] (D6) at (0,7.6) {$( \bgamma_h, \bgamma_h^d, \bgamma_h^s )$};

% 0 -> 5  
\draw[->] ([yshift=-23pt]D0.north) -- ([yshift= 26pt]D5.south) 
node[midway, xshift=22pt, font=\footnotesize, ] { $R(-\pi/2) $};
% 0 -> 55  
\draw[->] ([yshift=23pt]D0.south) -- ([yshift=-29pt]D55.north) 
node[midway, xshift=17pt, font=\footnotesize, ] { $R(\pi/2) $};

% 1 <-> 2  
\draw[->] ([xshift=2pt, yshift=1.5pt]D1.east) -- ([xshift=-2pt, yshift=2.5pt]D2.west) 
node[midway,   yshift=14pt,font=\footnotesize,align=center] { $ B_2(-\operatorname{arctanh}f)N_1^\mp(f)$,\\ $fn_h+m_h=0$, $|f |<1$};
\draw[->] ([xshift=-2pt, yshift=-2.5pt]D2.west) -- ([xshift=2pt, yshift=-2.5pt]D1.east) 
node[midway,   yshift=-14pt,font=\footnotesize,align=center] { $N_1^\mp(-\tanh\varphi)B_2(\varphi)$,\\ $m_h \cosh\varphi + n_h \sinh\varphi =0$};
%
% 3 <-> 4  
\draw[->] ([xshift=2pt, yshift=1.5pt]D3.east) -- ([xshift=-2pt, yshift=2.5pt]D4.west) 
node[midway,   yshift=9pt,font=\footnotesize,align=center] {$B_1( \int_{t_0}^t{m_h}{\rm d}t)N_2^\mp (s_\pm)$};
\draw[->] ([xshift=-2pt, yshift=-2.5pt]D4.west) -- ([xshift=2pt, yshift=-2.5pt]D3.east) 
node[midway,   yshift=-11pt,font=\footnotesize,align=center] {$ N_2^\mp (-s_\pm) B_1(-\int_{t_0}^t{m_h}{\rm d}t)$ };

% 5  -> 1
\draw[->] ([xshift=-2pt, yshift=3pt]D5.west) -- ([xshift=0pt, yshift=50pt]D1.south) 
node[midway,font=\footnotesize,align=center, xshift=-23pt,yshift=-8pt] {$N_1^{\mp}(-f) $,\\ $fn_h+m_h=0$};

% 5  -> 2
\draw[->] ([xshift=2pt, yshift=3pt]D5.east) --([xshift=0pt, yshift=50pt]D2.south) 
node[midway,font=\footnotesize,align=center,  xshift=29pt,yshift=-8pt] {$B_2(-\varphi) $,\\ $m_h \cosh\varphi + n_h \sinh\varphi =0$};

%	5  -> 3 
\draw[->] ([xshift=-2pt, yshift=-3pt]D55.west) -- ([xshift=5pt, yshift=-54pt]D3.north) 
node[midway,font=\footnotesize,align=center, xshift=-18pt,yshift=8pt] {$N_2^\mp (-s_\pm) $};
%	5  -> 4 
\draw[->] ([xshift=2pt, yshift=-3pt]D55.east) -- ([xshift=-5pt, yshift=-50pt]D4.north) 
node[midway,font=\footnotesize,align=center, xshift=30pt,yshift=5pt] {$B_1( \int_{t_0}^t{m_h}{\rm d}t)$};
%
%%	11 -> 1
\draw[<-> ] ([xshift=-35pt, yshift=56pt]D11.south) -- ([xshift=-35pt, yshift=-50pt]D1.north) 
node[midway,font=\footnotesize,align=center, xshift= 39pt,yshift=-2pt] {The image of\\$\mathcal{HE}^\pm(\bgamma_h,\bgamma_h^d)$ is\\ contained in \\the set of \\singular values \\of   $\mathcal{HP}_{\lambda_\pm}^\pm (\bgamma_h,\bgamma_h^d)$ };
%%	22  -> 2
\draw[<- >] ([xshift=45pt, yshift=50pt]D22.south) -- ([xshift=45pt, yshift=-52pt]D2.north) 
node[midway,font=\footnotesize,align=center, xshift=-36pt,yshift=0pt] {The image of\\ $HE(\bgamma_h,\bgamma_h^d)$ is\\ contained in \\the set of \\singular values \\of   $HP_{\theta}(\bgamma_h,\bgamma_h^d)$};

% 6  -> 11
\draw[->] ([xshift=-2pt, yshift=-2pt]D6.west) -- ([xshift=0pt, yshift=74pt]D11.south) 
node[midway,font=\footnotesize,align=center, xshift=-5pt,yshift=10pt] {$N_2^\mp ( \lambda_\pm)$ };

% 6  -> 22
\draw[->] ([xshift=1pt, yshift=-2pt]D6.east) --([xshift=0pt, yshift=66pt]D22.south) 
node[midway,font=\footnotesize,align=center,  xshift=12pt,yshift=6pt] {$B_1( \theta) $ };

% 11 <-> 22  
\draw[->] ([xshift=2pt, yshift=1.5pt]D11.east) -- ([xshift=-2pt, yshift=2.5pt]D22.west) 
node[midway,   yshift=9pt,font=\footnotesize,align=center] {$ B_1( \theta) N_2^\mp(- \lambda_\pm )$};
\draw[->] ([xshift=-2pt, yshift=-2.5pt]D22.west) -- ([xshift=2pt, yshift=-2.5pt]D11.east) 
node[midway,   yshift=-8pt,font=\footnotesize,align=center] { $N_2^\mp( \lambda_\pm )B_1(-\theta)$};
% 1 <-> 3 
\draw[<->] ([xshift=-36pt, yshift=-54pt]D3.north) -- ([xshift=-36pt, yshift=51pt]D1.south) 
node[midway,font=\footnotesize,align=center, xshift=33pt,yshift=5pt] {Corresponding};
% 2 <-> 4 
\draw[<->] ([xshift=47pt, yshift=-51pt]D4.north) -- ([xshift=47pt, yshift=48pt]D2.south) 
node[midway,font=\footnotesize,align=center, xshift= -33pt,yshift=5pt] {Corresponding};

\pgfresetboundingbox
\path (-6.5, -4.0) rectangle (6.5, 8.5); 

\end{tikzpicture}

\caption{The relations among the evolutes, involutes  and parallels of $(\bgamma_h,\bgamma_h^d)$.} 
\label{figure2} 
\end{figure}
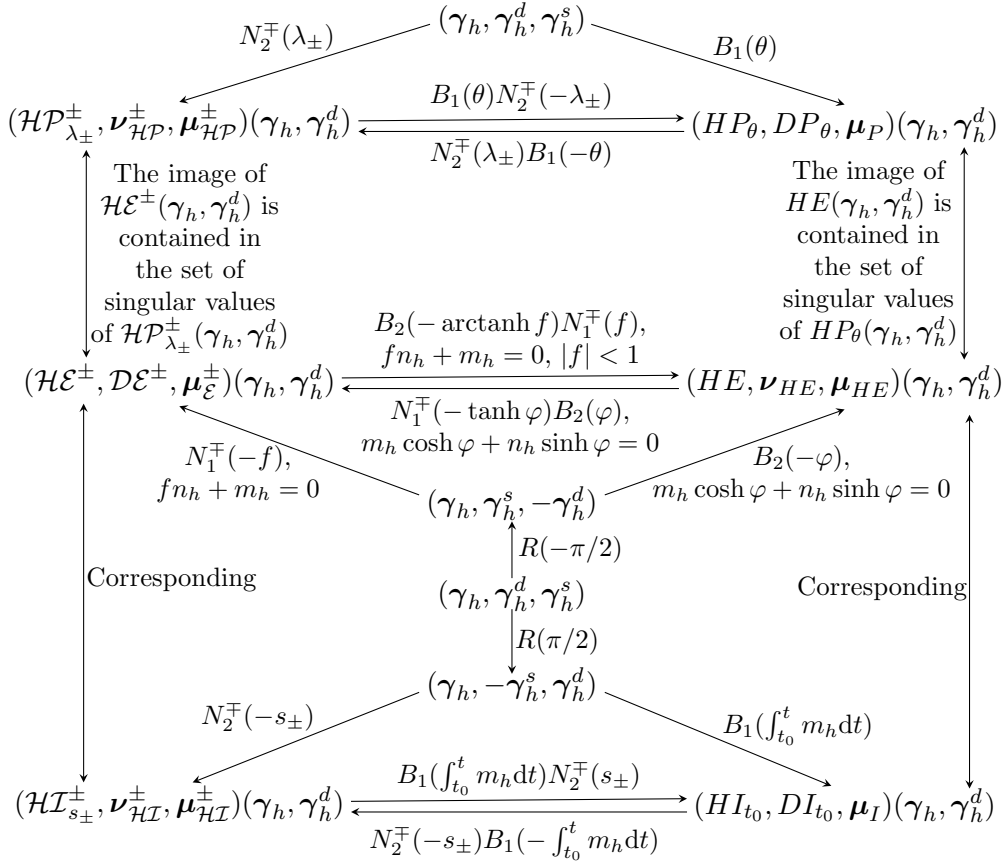

Proposition \ref{prop0} implies that the image of hyperbolic evolute $HE(\bgamma_h,\bgamma_h^d)$ (respectively, de Sitter evolute $DE(\bgamma_h,\bgamma_h^d)$) is contained in the set of singular values of the hyperbolic parallels $HP_{\theta}(\bgamma_h,\bgamma_h^d)$ (respectively, de Sitter parallels $DP_{\theta}(\bgamma_h,\bgamma_h^d)$) of $(\bgamma_h,\bgamma_h^d)$. 
In \cite{NTZ2026}, we have shown that the image of the hyperbolic horocyclic evolute $\mathcal{HE}^\pm(\bgamma_h,\bgamma_h^d)$ is contained in the set of singular values of the hyperbolic horocyclic parallels $\mathcal{HP}_{\lambda_\pm}^\pm (\bgamma_h,\bgamma_h^d)$ of $(\bgamma_h,\bgamma_h^d)$. 
Moreover, Theorem \ref{th1} provides the correspondence between the hyperbolic horocyclic  evolute $ \mathcal{HE}^\pm(\bgamma_h,\bgamma_h^d) $  (respectively, hyperbolic   evolute $HE(\bgamma_h,\bgamma_h^d)$)  and hyperbolic horocyclic involute $\mathcal{HI}_{s_\pm}^\pm(\bgamma_h,\bgamma_h^d)$ (respectively, hyperbolic   involute $HI_{t_0}(\bgamma_h,\bgamma_h^d)$) of $(\bgamma_h,\bgamma_h^d)$. 
Consequently, we eventually establish the relations among the evolutes, involutes  and parallels of $(\bgamma_h,\bgamma_h^d)$, as illustrated in Figure \ref{figure2}.

%%%%%%%%%%% Section 4 %%%%%%%%%%%%%
\section{Evolutes,  involutes and parallels  of spacelike Legendre curves in $S_1^2$}\label{S4}
\subsection{Evolutes of spacelike Legendre curves}

%%%%%
\begin{definition}[Hyperbolic and de Sitter evolutes]\rm
Let $(\bgamma_d,\bgamma_d^h):I\rightarrow\overline{\Delta}_1$ be a spacelike Legendre curve  with curvature $(m_d,n_d)$. 
Suppose that  
there exists a smooth function $\overline{\psi}:I\rightarrow\R$ such that $m_d(t)\cosh\overline{\psi}(t)+ n_d(t)\sinh\overline{\psi}(t)=0$ for all $t\in I$. Then the \textit{de Sitter evolute} $\overline{DE}(\bgamma_d,\bgamma_d^h):I\rightarrow S_1^2$ of  $(\bgamma_d,\bgamma_d^h)$   is defined as  
\begin{equation}\notag
\overline{DE}(\bgamma_d,\bgamma_d^h)(t)=\cosh\overline{\psi}(t)\bgamma_d(t)+\sinh\overline{\psi}(t)\bgamma_d^h(t).
\end{equation} 
Suppose     there exists a smooth function $\overline{\varphi}:I\rightarrow\R$ such that $m_d(t)\sinh\overline{\varphi}(t)+ n_d(t)\cosh\overline{\varphi}(t)=0$ for all $t\in I$. Then
the  \textit{hyperbolic evolute} $\overline{HE}(\bgamma_d,\bgamma_d^h):I\rightarrow H^2$  of $(\bgamma_d,\bgamma_d^h)$  is defined as  
\begin{equation}\notag
\overline{HE}(\bgamma_d,\bgamma_d^h)(t)=\sinh\overline{\varphi}(t)\bgamma_d(t)+\cosh\overline{\varphi}(t)\bgamma_d^h(t).
\end{equation} 
\end{definition}
%%%%%

In this paper, we assume    $\mathrm{Reg}(\bgamma_d)=\{t\in I\mid m_d(t)\neq0 \}$ and  $\mathrm{Reg}(\bgamma_d^h)= \{t\in I\mid n_d(t)\neq0 \}$ are   dense  in $ I $. If  $\overline{\psi},\overline{\varphi}:I\rightarrow\R$ exist, then $\overline{\psi},\overline{\varphi}$ are unique.
%%%%
\begin{proposition} 
Let $(\bgamma_d,\bgamma_d^h):I\rightarrow\overline{\Delta}_1$ be a spacelike Legendre curve  with curvature $(m_d,n_d)$. 
\begin{enumerate} 
\item[$(1)$] 
Suppose that	   there exists   $\overline{\psi}:I\rightarrow\R$ such that $m_d(t)\cosh\overline{\psi}(t)+ n_d(t)\sinh\overline{\psi}(t)=0$ for all $t\in I$.  Then  
$(\overline{DE}(\bgamma_d,\bgamma_d^h), \overline{\bnu}_{DE}(\bgamma_d,\bgamma_d^h) ):I\rightarrow \Delta_5$ is a timelike Legendre curve with curvature 
$ 
(\overline{m} _{DE}(t), \overline{n} _{DE}(t))= ( \dot{\overline{\psi}}(t),  m_d(t)\sinh\overline{\psi}(t)+n_d(t)\cosh\overline{\psi}(t)),
$  
where 
$  
\overline{\bnu}_{DE}(\bgamma_d,\bgamma_d^h)(t) = \bgamma_d^s(t). 
$ 
\item[$(2)$] Suppose that  there exists $\overline{\varphi}:I\rightarrow\R$ such that $m_d(t)\sinh\overline{\varphi}(t)+ n_d(t)\cosh\overline{\varphi}(t)=0$ for all $t\in I$.    Then  
$(\overline{HE}(\bgamma_d,\bgamma_d^h), \overline{\bnu}_{HE}(\bgamma_d,\bgamma_d^h) ):I\rightarrow \Delta_1$ is a spacelike Legendre curve with curvature 
$ (\overline{m} _{HE}(t), \overline{n}_{HE}(t))=  (   \dot{\overline{\varphi}}(t),  -m_d(t)\cosh\overline{\varphi}(t)- n_d(t)\sinh\overline{\varphi}(t) ), 
$  
where 
$ \overline{\bnu}_{HE}(\bgamma_d,\bgamma_d^h)(t) =\bgamma_d^s(t). $ 
\end{enumerate}
\end{proposition}
%%%%

%%%%%
\begin{remark}\rm
If $(\bgamma_h,\bgamma_h^d)=(\bgamma_d^h,\bgamma_d)$, then the hyperbolic and de Sitter evolutes of $(\bgamma_h,\bgamma_h^d)$ and $(\bgamma_d,\bgamma_d^h)$ coincide (cf. \cite{IT2007}). That is,  $HE(\bgamma_h,\bgamma_h^d)=\overline{HE}(\bgamma_d,\bgamma_d^h)$ and $DE(\bgamma_h,\bgamma_h^d)=\overline{DE}(\bgamma_d,\bgamma_d^h)$.
\end{remark}
%%%%%

%%%%
\begin{proposition}   
Suppose   there exists   $\overline{\psi}:I\rightarrow\R$ such that $m_d(t)\cosh\overline{\psi}(t)+ n_d(t)\sinh\overline{\psi}(t)=0$ for all $t\in I$. Let $\overline{\bmu}_{DE}(\bgamma_d,\bgamma_d^h)= \overline{DE}(\bgamma_d,\bgamma_d^h)\wedge \overline{\bnu}_{DE}(\bgamma_d,\bgamma_d^h)$ and denote the moving frame of $\overline{DE}(\bgamma_d,\bgamma_d^h)$ by $(\overline{DE}, \overline{\bnu}_{DE},\overline{\bmu}_{DE}   ) (\bgamma_d,\bgamma_d^h)$.  Then  
$$(\overline{DE}, \overline{\bnu}_{DE},\overline{\bmu}_{DE}   )^T(\bgamma_d,\bgamma_d^h) =\overline{B}_1(\overline{\psi} ) ( \bgamma_d,   \bgamma_d^s,\bgamma_d^h  )^T .$$
That is,  $(\overline{DE}(\bgamma_d,\bgamma_d^h)(t), \overline{\bnu}_{DE}(\bgamma_d,\bgamma_d^h)(t),\overline{\bmu}_{DE}(\bgamma_d,\bgamma_d^h)(t))$ is obtained by applying a  boost in the $\bgamma_d  $-direction with rapidity $ \overline{\psi}(t)$  to  $(\bgamma_d(t),   \bgamma_d^s(t),\bgamma_d^h(t))$ for all $t\in I$.
\end{proposition}
%%%%

\subsection{Involutes of spacelike Legendre curves}
\begin{definition}[De Sitter horocyclic involute]\rm
Let $(\bgamma_d,\bgamma_d^h):I\rightarrow\overline{\Delta}_1$ be a spacelike Legendre curve  with curvature $(m_d,n_d)$. The  \textit{de Sitter horocyclic involute} $\overline{\mathcal{DI}}_{\overline{s}_\pm}^\pm(\bgamma_d,\bgamma_d^h):I\rightarrow S_1^2$ of $(\bgamma_d,\bgamma_d^h)$  is defined as 
$$
\overline{\mathcal{DI}}_{\overline{s}_\pm}^\pm(\bgamma_d,\bgamma_d^h)(t)=\bgamma_d(t)+\overline{s}_\pm(t)\bgamma_d^s(t)-\frac{\overline{s}_\pm^2(t)}{2}(\bgamma_d(t)\mp\bgamma_d^h(t)),$$
where $\overline{s}_\pm(t)$ is a solution of the Riccati equation
\begin{equation}\label{Riccati2}
\frac{{\rm d}\overline{s}_\pm}{{\rm d}t}(t)=- \frac{m_d(t)\mp n_d(t)}{2}\overline{s}_\pm^2(t)-m_d(t).\end{equation}
\end{definition}
%%%%

%%%%%
\begin{proposition}
Let $	\overline{\mathcal{DI}}_{\overline{s}_\pm}^\pm(\bgamma_d,\bgamma_d^h)$   be a   de Sitter horocyclic involute  of $(\bgamma_d,\bgamma_d^h)$. Then\\
$(\overline{\mathcal{DI}}_{\overline{s}_\pm}^\pm(\bgamma_d,\bgamma_d^h), \overline{\bnu}_\mathcal{DI}^\pm  (\bgamma_d,\bgamma_d^h)):I\rightarrow \Delta_5$ is a timelike Legendre curve with curvature  
$$
\left(\overline{m}^\pm_\mathcal{DI}(t), \overline{n}^\pm_\mathcal{DI}(t)\right)= \left(  \overline{s}_\pm(t)(n_d(t)\mp m_d(t)),n_d(t)\mp m_d(t)\right),  
$$ 
where $\overline{\bnu}_\mathcal{DI}^\pm(\bgamma_d,\bgamma_d^h)(t) =-\overline{s}_\pm(t)\bgamma_d(t)\pm\overline{s}_\pm(t)\bgamma_d^h(t)+\bgamma_d^s(t)$. 	 
\end{proposition}
%%%%%
%%%%
\begin{proposition} \label{prop7}
 Let  $  \overline{\bmu}^\pm_\mathcal{DI} (\bgamma_d,\bgamma_d^h)=\overline{\mathcal{DI}}_{\overline{s}_\pm}^\pm(\bgamma_d,\bgamma_d^h)\wedge \overline{\bnu}_\mathcal{DI}^\pm (\bgamma_d,\bgamma_d^h) $  and denote the moving frame of $	\overline{\mathcal{DI}}_{\overline{s}_\pm}^\pm(\bgamma_d,\bgamma_d^h)$ by  $	 ( \overline{\mathcal{DI}}_{\overline{s}_\pm}^\pm,\overline{\bnu}_\mathcal{DI}^\pm ,   \overline{\bmu}^\pm_{\mathcal{DI} }   ) (\bgamma_d,\bgamma_d^h) $.
Then  
$$ ( \overline{\mathcal{DI}}_{\overline{s}_\pm}^\pm,\overline{\bnu}_\mathcal{DI}^\pm ,   \overline{\bmu}^\pm_\mathcal{DI}  )^T(\bgamma_d,\bgamma_d^h) =\overline{N}_1^\pm( \overline{s}_\pm ) ( \bgamma_d,   \bgamma_d^s,\bgamma_d^h )^T . $$
That is,  $ ( \overline{\mathcal{DI}}_{\overline{s}_\pm}^\pm(\bgamma_d,\bgamma_d^h)(t),\overline{\bnu}_\mathcal{DI}^\pm(\bgamma_d,\bgamma_d^h) (t),   \overline{\bmu}^\pm_\mathcal{DI}(\bgamma_d,\bgamma_d^h) (t)  )$ is obtained by applying a null rotation   around the $(\bgamma_d\pm \bgamma_d^h)$-direction with angle $  \overline{s}_\pm(t)$ to $ (\bgamma_d(t),   \bgamma_d^s(t),\bgamma_d^h(t)  )$ for all $t\in I$. 
\end{proposition}
%%%%
%%%%%
\begin{definition}[De Sitter involute]\rm
Let $(\bgamma_d,\bgamma_d^h):I\rightarrow\overline{\Delta}_1$ be a spacelike Legendre curve  with curvature $(m_d,n_d)$ and  $t_0\in I$.
The  \textit{de Sitter involute} $\overline{DI}_{t_0}(\bgamma_d,\bgamma_d^h):I\rightarrow S_1^2$   of $(\bgamma_d,\bgamma_d^h)$  is defined as  
\begin{equation}\notag
\overline{DI}_{t_0}(\bgamma_d,\bgamma_d^h)(t)=   \cos\left(\int_{t_0}^t{m_d(t)}{\rm d}t\right)\bgamma_d(t)-\sin\left(\int_{t_0}^t{m_d(t)}{\rm d}t\right)\bgamma_d^s(t).
\end{equation}  
\end{definition}
%%%%%

%%%%
\begin{proposition} 
Let $ \overline{DI}_{t_0}(\bgamma_d,\bgamma_d^h) $  be  a de Sitter involute  of  $(\bgamma_d,\bgamma_d^h)$. Then\\
$ ( \overline{DI}_{t_0}(\bgamma_d,\bgamma_d^h), \overline{\bnu}_{DI}(\bgamma_d,\bgamma_d^h)   ):I\rightarrow \Delta_5$ is a timelike Legendre curve with curvature 
$$
\left(\overline{m }_{DI}(t), \overline{n}_{DI}(t)\right)= \left(     -n_d(t)\sin\left(\int_{t_0}^t{m_d(t)}{\rm d}t\right),  n_d(t)\cos\left(\int_{t_0}^t{m_d(t)}{\rm d}t\right)\right),  
$$ 
where	 $ \overline{\bnu}_{DI}(\bgamma_d,\bgamma_d^h)(t) =\cos (\int_{t_0}^t{m_d(t)}{\rm d}t )\bgamma_d^s(t)+\sin (\int_{t_0}^t{m_d(t)}{\rm d}t )\bgamma_d(t).$ 
\end{proposition}
%%%%

%%%%
\begin{proposition}\label{prop8}   
 
Let $\overline{\bmu}_{DI}(\bgamma_d,\bgamma_d^h)= \overline{DI}_{t_0}(\bgamma_d,\bgamma_d^h)\wedge \overline{\bnu}_{DI}(\bgamma_d,\bgamma_d^h)$ and denote the moving frame of $\overline{DI}_{t_0}(\bgamma_d,\bgamma_d^h)$ by $ (\overline{DI}_{t_0}, \overline{\bnu}_{DI},\overline{\bmu}_{DI}    )(\bgamma_d,\bgamma_d^h)$.  Then  
$$\begin{aligned}\left(\overline{DI}_{t_0}, \overline{\bnu}_{DI},\overline{\bmu}_{DI}   \right)^T(\bgamma_d,\bgamma_d^h) =\overline{R}\left(\int_{t_0}^t{m_d }{\rm d}t\right) \left( \bgamma_d,   \bgamma_d^s,\bgamma_d^h  \right)^T .\end{aligned}$$
That is,  $ (\overline{DI}_{t_0}(\bgamma_d,\bgamma_d^h)(t), \overline{\bnu}_{DI}(\bgamma_d,\bgamma_d^h)(t),\overline{\bmu}_{DI}(\bgamma_d,\bgamma_d^h)(t) )$ is obtained by applying a  rotation with  angle $\int_{t_0}^t{m_d(t)}{\rm d}t$   to  $ (\bgamma_d(t),   \bgamma_d^s(t),\bgamma_d^h(t) )$ for all $t\in I$.
\end{proposition}
%%%%

Combining Propositions \ref{prop7} and  \ref{prop8}, we have the relations between  the de Sitter horocyclic and de Sitter  involutes of $(\bgamma_d,\bgamma_d^h)$.
%%%%
\begin{theorem} 
Let $(\bgamma_d,\bgamma_d^h):I\rightarrow\overline{\Delta}_1$ be a spacelike Legendre curve  with curvature $(m_d,n_d)$.  	     $ \overline{\mathcal{DI}}_{\overline{s}_\pm}^\pm(\bgamma_d,\bgamma_d^h)  $ and
$ \overline{DI}_{t_0}(\bgamma_d,\bgamma_d^h)$  are de Sitter horocyclic and de Sitter involutes of $(\bgamma_d,\bgamma_d^h)$, respectively. Then  
$$\begin{aligned} ( \overline{\mathcal{DI}}_{\overline{s}_\pm}^\pm,\overline{\bnu}_\mathcal{DI}^\pm ,   \overline{\bmu}^\pm_\mathcal{DI}     )^T  (\bgamma_d,\bgamma_d^h) &=  \overline{N}_1^\pm ( \overline{s}_\pm  )  \overline{R}\left(-\int_{t_0}^t{m _d}{\rm d}t\right)\left(  \overline{DI}_{t_0}, \overline{\bnu}_{DI},\overline{\bmu}_{DI}   \right)^T   (\bgamma_d,\bgamma_d^h)  ,\\
(  \overline{DI}_{t_0}, \overline{\bnu}_{DI},\overline{\bmu}_{DI}   )^T  	(\bgamma_d,\bgamma_d^h)   &= \overline{R}\left( \int_{t_0}^t{m _d}{\rm d}t\right) \overline{N}_1^\pm(- \overline{s}_\pm  )( \overline{\mathcal{DI}}_{\overline{s}_\pm}^\pm,\overline{\bnu}_\mathcal{DI}^\pm,   \overline{\bmu}^\pm_\mathcal{DI}      )^T   (\bgamma_d,\bgamma_d^h) .\end{aligned}$$ 
\end{theorem}
%%%%
\subsection{Parallels of spacelike Legendre curves}
%%%%%
\begin{definition}[De Sitter and  hyperbolic   parallels \cite{IT2007}]  \rm
Let $(\bgamma_d,\bgamma_d^h):I\rightarrow\overline{\Delta}_1$ be a spacelike Legendre curve  with curvature $(m_d,n_d)$ and $\theta\in\R$.  
The \textit{de Sitter parallel} $\overline{DP}_{\theta}(\bgamma_d,\bgamma_d^h):I\rightarrow S_1^2$      of $(\bgamma_d,\bgamma_d^h)$ is   defined as
$$
\overline{DP}_{\theta}(\bgamma_d,\bgamma_d^h)(t)=\sinh \theta\bgamma_d^h(t)+\cosh\theta\bgamma_d(t).$$
The \textit{hyperbolic parallel} $\overline{HP}_{\theta}(\bgamma_d,\bgamma_d^h):I\rightarrow H^2$   of $(\bgamma_d,\bgamma_d^h)$ is   defined as
$$
\overline{HP}_{\theta}(\bgamma_d,\bgamma_d^h)(t)=\cosh \theta\bgamma_d^h(t)+\sinh\theta\bgamma_d(t).$$
\end{definition}
%%%%%

%%%%
\begin{proposition} 
Let $ \overline{DP}_{\theta}(\bgamma_d,\bgamma_d^h)$  and $\overline{HP}_{\theta}(\bgamma_d,\bgamma_d^h)$ be   de Sitter and hyperbolic  parallels of $(\bgamma_d,\bgamma_d^h)$, respectively.  Then   \begin{enumerate} 
\item[$(1)$] 
$ ( \overline{DP}_{\theta}(\bgamma_d,\bgamma_d^h),  \overline{HP}_{\theta}(\bgamma_d,\bgamma_d^h)  ):I\rightarrow \overline{\Delta}_1$ is a spacelike Legendre curve with curvature 
$ (	\overline{m} _{DP}(t), 	\overline{n}_{DP}(t) )=  (     m_d(t) \cosh \theta+ n_d(t)\sinh\theta,  n_d(t) \cosh \theta+ m_d(t) \sinh\theta).
$  
\item[$(2)$]   
$ ( \overline{HP}_{\theta}(\bgamma_d,\bgamma_d^h),\overline{DP}_{\theta}(\bgamma_d,\bgamma_d^h)  ):I\rightarrow \Delta_1$ is a spacelike Legendre curve with curvature 
$ (	\overline{m }_{HP}(t), 	\overline{n}_{HP}(t) )=  (    - n_d(t)\cosh \theta- m_d(t)\sinh\theta, -  m_d(t)\cosh \theta- n_d(t)\sinh\theta ).
$  
\end{enumerate}
\end{proposition}
%%%%

%%%%
\begin{remark}\rm
If $(\bgamma_h,\bgamma_h^d)=(\bgamma_d^h,\bgamma_d)$, then the hyperbolic and de Sitter parallels of $(\bgamma_h,\bgamma_h^d)$ and $(\bgamma_d,\bgamma_d^h)$ coincide. That is,  $HP_{\theta}(\bgamma_h, \bgamma_h^d)=\overline{HP}_{\theta}(\bgamma_d,\bgamma_d^h)$ and $DP_{\theta}(\bgamma_h, \bgamma_h^d)=\overline{DP}_{\theta}(\bgamma_d,\bgamma_d^h)$.
\end{remark}
%%%%
%%%%
\begin{proposition}  
 Let $\overline{\bmu}_{P}(\bgamma_d,\bgamma_d^h)=\overline{DP}_{\theta}(\bgamma_d,\bgamma_d^h)\wedge \overline{HP}_{\theta}(\bgamma_d,\bgamma_d^h)$ and denote the moving frame of $\overline{DP}_{\theta}(\bgamma_d,\bgamma_d^h)$ by $ ( 	\overline{DP}_{\theta},\overline{\bmu}_{P},  \overline{HP}_{\theta}   )  (\bgamma_d,\bgamma_d^h)$.  Then      
$$( 	\overline{DP}_{\theta},\overline{\bmu}_{P},  \overline{HP}_{\theta}  )^T (\bgamma_d,\bgamma_d^h)
=\overline{B}_1(\theta) ( \bgamma_d,   \bgamma_d^s,\bgamma_d^h )^T .$$
That is,	$ ( 	\overline{DP}_{\theta}(\bgamma_d,\bgamma_d^h)(t), \overline{\bmu}_{P}(\bgamma_d,\bgamma_d^h)(t), \overline{HP}_{\theta}(\bgamma_d,\bgamma_d^h)(t)   )    $ is obtained by applying a  boost in the $\bgamma_d $-direction    with rapidity   $\theta$  to $ ( \bgamma_d(t),   \bgamma_d^s(t),\bgamma_d^h(t)  )  $ for all $t\in I$.   
\end{proposition}
%%%%

\begin{figure}[h] 
\centering
\begin{tikzpicture}[
>=stealth, % 
node distance=3cm,
vertex/.style={circle, minimum size=6mm, inner sep=0pt, font= \small}
]

\node[vertex] (D1) at (-3.5, 3) {$(\overline{DE}, \overline{\bnu}_{DE},\overline{\bmu}_{DE}   )(\bgamma_d,\bgamma_d^h) $};
\node[vertex] (D2) at (3.5, 3) {$ (\overline{DP}_{\theta},\overline{\bmu}_{P},  \overline{HP}_{\theta}   ) (\bgamma_d,\bgamma_d^h)  $};
\node[vertex] (D3) at (-3.8, -0.3) {$( \overline{\mathcal{DI}}_{\overline{s}_\pm}^\pm,\overline{\bnu}_\mathcal{DI}^\pm ,   \overline{\bmu}^\pm_\mathcal{DI}   )(\bgamma_d,\bgamma_d^h)$};
\node[vertex] (D4) at (3.8, -0.3) {$ (\overline{DI}_{t_0}, \overline{\bnu}_{DI},\overline{\bmu}_{DI}   ) (\bgamma_d,\bgamma_d^h)$};
\node[vertex] (D5) at (0,1.5) {$(   \bgamma_d,   \bgamma_d^s,\bgamma_d^h    )$};

% 1 <- > 2  
\draw[<-> ] ([xshift=-2pt, yshift=-2 pt]D2.west) -- ([xshift=2pt, yshift=-2 pt]D1.east) 
node[midway,   yshift=25pt,font=\footnotesize,align=center] {The image of $\overline{DE}(\bgamma_d,\bgamma_d^h)$\\ is  contained  in  the  set\\ of  singular  values \\of   $\overline{DP}_\theta(\bgamma_d,\bgamma_d^h)$};

% 3 <-> 4  
\draw[->] ([xshift=2pt, yshift=1.5pt]D3.east) -- ([xshift=-2pt, yshift=2.5pt]D4.west) 
node[midway,   yshift=10pt,font=\footnotesize,align=center] {$ \overline{R}( \int_{t_0}^t{m _d}{\rm d}t) \overline{N}_1^\pm( -\overline{s}_\pm )$};
\draw[->] ([xshift=-2pt, yshift=-2.5pt]D4.west) -- ([xshift=2pt, yshift=-2.5pt]D3.east) 
node[midway,   yshift=-12pt,font=\footnotesize,align=center] { $\overline{N}_1^\pm ( \overline{s}_\pm  )  \overline{R}(-\int_{t_0}^t{m _d}{\rm d}t)$};

% 5  -> 1
\draw[->] ([xshift=-2pt, yshift=3pt]D5.west) -- ([xshift=0pt, yshift=49pt]D1.south) 
node[midway,font=\footnotesize,align=center, xshift=-11pt,yshift=-7pt] {$\overline{B}_1(\overline{\psi})$};

% 5  -> 2
\draw[->] ([xshift=2pt, yshift=3pt]D5.east) --([xshift=0pt, yshift=52pt]D2.south) 
node[midway,font=\footnotesize,align=center,  xshift=12pt,yshift=-7pt] {$\overline{B}_1(\theta)  $};

%	5  -> 3 
\draw[->] ([xshift=-2pt, yshift=-3pt]D5.west) -- ([xshift=1pt, yshift=-55pt]D3.north) 
node[midway,font=\footnotesize,align=center, xshift=-16pt,yshift=6pt] {$\overline{N}_1^\pm( \overline{s}_\pm) $};
%	5  -> 4 
\draw[->] ([xshift=2pt, yshift=-3pt]D5.east) -- ([xshift=-5pt, yshift=-52pt]D4.north) 
node[midway,font=\footnotesize,align=center, xshift=29pt,yshift=3pt] {$\overline{R}(\int_{t_0}^t{m_d }{\rm d}t)$};

\pgfresetboundingbox
\path (-6, -1.3) rectangle (6, 5.0);

\end{tikzpicture}
\caption{The relations among the evolutes, involutes and parallels of $(\bgamma_d,\bgamma_d^h)$.}  
\label{figure3}
\end{figure}
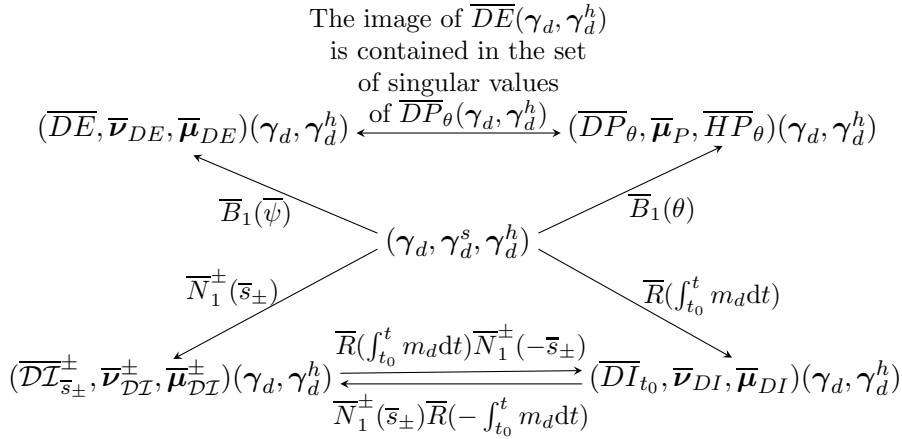

Based on the above conclusions, we present the relations among the evolutes, involutes and parallels of $(\bgamma_d,\bgamma_d^h)$, as illustrated in Figure \ref{figure3}.

\vspace{0.3cm}

%%%%%%%%%%% Section 5 %%%%%%%%%%%%%
\section{Evolutes,   involutes and parallels of timelike Legendre curves in $S_1^2$}\label{S5}
\subsection{Evolutes of timelike Legendre curves}
%%%%%
\begin{definition}[De Sitter and hyperbolic   horocyclic evolutes]\rm
Let $(\bgamma_T,\bgamma_T^d):I\rightarrow\Delta_5$ be a timelike Legendre curve  with curvature $(m_T,n_T)$. Suppose that there exists a smooth function $\widetilde{f}:I\rightarrow\R  $ such that $m_T(t)+\widetilde{f}(t)n_T(t)=0$ for all $t\in I$. 
Then the \textit{de Sitter horocyclic evolute}      $\widetilde{\mathcal{DE}}^\pm(\bgamma_T,\bgamma_T^d):I\rightarrow S_1^2$ of $(\bgamma_T,\bgamma_T^d)$ is defined as  
\begin{equation}\notag
\widetilde{\mathcal{DE}}^\pm(\bgamma_T,\bgamma_T^d)(t)=\bgamma_T(t)+\widetilde{f}(t)\bgamma_T^d(t)-\frac{\widetilde{f}^2(t)}{2}(\bgamma_T(t)\mp\bgamma_T^h(t)).
\end{equation} 
The \textit{hyperbolic horocyclic evolute}  $\widetilde{\mathcal{HE}}^\pm(\bgamma_T,\bgamma_T^d):I\rightarrow H^2$   of $(\bgamma_T,\bgamma_T^d)$   is defined as  
\begin{equation}\notag
\widetilde{\mathcal{HE}}^\pm(\bgamma_T,\bgamma_T^d)(t)=\mp\frac{\widetilde{f}^2(t)}{2}\bgamma_T(t)\pm \widetilde{f}(t)\bgamma_T^d(t)+ \left(1+\frac{\widetilde{f}^2(t)}{2}\right)\bgamma_T^h(t).
\end{equation} 
\end{definition}
%%%%%

In this paper, we assume that $\mathrm{Reg}(\bgamma_T^d)= \{t\in I\mid n_T(t)\neq0 \}$  is   dense  in $ I $. Then, if  $\widetilde{f} :I\rightarrow\R $ exists,  it is unique.

%%%%
\begin{proposition} 
Suppose that there exists  $\widetilde{f}:I\rightarrow\R  $ such that $m_T(t)+\widetilde{f}(t)n_T(t)=0$ for all $t\in I$. Then
\begin{enumerate} 
\item[$(1)$] $ ( \widetilde{\mathcal{DE}}^\pm(\bgamma_T,\bgamma_T^d),\widetilde{\mathcal{HE}}^\pm(\bgamma_T,\bgamma_T^d)   ):I\rightarrow \overline{\Delta}_1$ is a spacelike Legendre curve with curvature 
$$
\left(\widetilde{m}^\pm_\mathcal{DE}(t), \widetilde{n}^\pm_\mathcal{DE}(t)\right)= \left(   \dot{\widetilde{f}}(t)\pm ({\widetilde{f}^2(t)n_T(t)}/{2}),\pm \dot{\widetilde{f}}(t)+ ({\widetilde{f}^2(t)n_T(t)}/{2})+n_T(t)\right).
$$ 

\item[$(2)$] $ ( \widetilde{\mathcal{HE}}^\pm(\bgamma_T,\bgamma_T^d),\widetilde{\mathcal{DE}}^\pm(\bgamma_T,\bgamma_T^d)  ):I\rightarrow \Delta_1$ is a spacelike Legendre curve with curvature 
$$
\left(\widetilde{m}^\pm_\mathcal{HE}(t), \widetilde{n}^\pm_\mathcal{HE}(t)\right)= \left( \mp \dot{\widetilde{f}}(t)-({\widetilde{f}^2(t)n_T(t)}/{2})-n_T(t),  -\dot{\widetilde{f}}(t)\mp ({\widetilde{f}^2(t)n_T(t)}/{2}) \right).
$$ 
\end{enumerate}
\end{proposition}
%%%%

%%%%
\begin{proposition}\label{prop9} 
Suppose that there exists   $\widetilde{f}:I\rightarrow\R  $ such that $m_T(t)+\widetilde{f}(t)n_T(t)=0$ for all $t\in I$.  Let  $\widetilde{\bmu}^\pm_\mathcal{E}(\bgamma_T,\bgamma_T^d) =\widetilde{\mathcal{DE}}^\pm(\bgamma_T,\bgamma_T^d) \wedge  \widetilde{\mathcal{HE}}^\pm(\bgamma_T,\bgamma_T^d) $ and denote the moving frame of $\widetilde{\mathcal{DE}}^\pm(\bgamma_T,\bgamma_T^d)$ by $ (\widetilde{\mathcal{DE}}^\pm,   \widetilde{\bmu}^\pm_\mathcal{E} ,\widetilde{\mathcal{HE}}^\pm   ) (\bgamma_T,\bgamma_T^d)$. Then
$$(\widetilde{\mathcal{DE}}^\pm,   \widetilde{\bmu}^\pm_\mathcal{E} ,\widetilde{\mathcal{HE}}^\pm )^T(\bgamma_T,\bgamma_T^d) =\overline{N}_1^{\pm}( \widetilde{f} )  (\bgamma_T, \bgamma_T^d,\bgamma_T^h ) ^T .$$
That is,  $ ( \widetilde{\mathcal{DE}}^\pm(\bgamma_T,\bgamma_T^d)(t), \widetilde{\bmu}^\pm_\mathcal{E}(\bgamma_T,\bgamma_T^d)(t),\widetilde{\mathcal{HE}}^\pm(\bgamma_T,\bgamma_T^d)(t)  )$ is obtained by applying a null rotation   around the $(\bgamma_T\pm \bgamma_T^h)$-direction with angle $ \widetilde{f}(t)$ to $(\bgamma_T(t), \bgamma_T^d(t),\bgamma_T^h(t) )$   for all $t\in I$. 
\end{proposition}
%%%%% 

%%%%%
\begin{definition}[De Sitter evolute]\rm
Let $(\bgamma_T,\bgamma_T^d ):I\rightarrow\Delta_5$ be a timelike Legendre curve  with curvature  $(m_T,n_T)$. 
Suppose that
there exists a smooth function $\widetilde{\varphi}:I\rightarrow\R$ such that $m_T(t)\cos\widetilde{\varphi}(t)+ n_T(t)\sin\widetilde{\varphi}(t)=0$ for all $t\in I$. 
Then the  \textit{de Sitter evolute} $\widetilde{DE}(\bgamma_T,\bgamma_T^d):I\rightarrow S_1^2$   of $(\bgamma_T,\bgamma_T^d)$  is defined as  
\begin{equation}\notag
\widetilde{DE}(\bgamma_T,\bgamma_T^d)(t)=\cos\widetilde{\varphi}(t)\bgamma_T(t) +\sin\widetilde{\varphi}(t)\bgamma_T^d(t).
\end{equation} 

\end{definition}
%%%%%

If $m_T(t)\cos\widetilde{\varphi}(t)+n_T(t)\sin\widetilde{\varphi}(t)=0$, we have
$ m_T(t)\cos\big(\widetilde{\varphi}(t)+k\pi\big)+n_T(t)\sin\big(\widetilde{\varphi}(t)+k\pi\big)=0 $ 
for all  $ k\in\mathbb{Z} $.
Hence we consider the equivalence class
$ [\widetilde{\varphi}(t)] := \widetilde{\varphi}(t)  \pmod{\pi} $, and by abuse of notation, we still denote it by $\widetilde{\varphi}(t)$.
Recall that the set   $\mathrm{Reg}(\bgamma_T^d)= \{t\in I\mid n_T(t)\neq0 \}$  is   dense  in $ I $. Hence, if $\widetilde{\varphi}$ exists, it is unique in the sense of the equivalence class. 
%%%%
\begin{proposition} 
Suppose that
there exists  $\widetilde{\varphi}:I\rightarrow\R$ such that $m_T(t)\cos\widetilde{\varphi}(t)+ n_T(t)\sin\widetilde{\varphi}(t)=0$ for all $t\in I$.  Then
$( \widetilde{DE}(\bgamma_T,\bgamma_T^d), \widetilde{\bnu}_{DE} (\bgamma_T,\bgamma_T^d) ):I\rightarrow \overline{\Delta}_1$ is a spacelike Legendre curve with curvature 
$ (\widetilde{m }_{DE}(t), \widetilde{n}_{DE}(t) )=  (     \dot{\widetilde{\varphi}}(t),  n_T(t)\cos\widetilde{\varphi}(t)-m_T(t)\sin\widetilde{\varphi}(t) ),  
$ 
where	 $ \widetilde{\bnu}_{DE}(\bgamma_T,\bgamma_T^d)(t) = \bgamma_T^h(t).$ 
\end{proposition}
%%%%

%%%%
\begin{proposition} \label{prop10}  
 
Suppose that
there exists  $\widetilde{\varphi}:I\rightarrow\R$ such that $m_T(t)\cos\widetilde{\varphi}(t)+ n_T(t)\sin\widetilde{\varphi}(t)=0$ for all $t\in I$.  Let $\widetilde{\bmu}_{DE}(\bgamma_T,\bgamma_T^d)= \widetilde{DE}(\bgamma_T,\bgamma_T^d)\wedge \widetilde{\bnu}_{DE}(\bgamma_T,\bgamma_T^d)$ and denote the moving frame of $\widetilde{DE}(\bgamma_T,\bgamma_T^d)$ by $ ( \widetilde{DE}, \widetilde{\bmu}_{DE},\widetilde{\bnu}_{DE}   ) (\bgamma_T,\bgamma_T^d)$.  Then  
$$( \widetilde{DE}, \widetilde{\bmu}_{DE},\widetilde{\bnu}_{DE} )^T(\bgamma_T,\bgamma_T^d) =\overline{R}(-\widetilde{\varphi} ) ( \bgamma_T,   \bgamma_T^d ,\bgamma_T^h  )^T .$$
That is,  $(\widetilde{DE}(\bgamma_T,\bgamma_T^d)(t), \widetilde{\bmu}_{DE}(\bgamma_T,\bgamma_T^d)(t),\widetilde{\bnu}_{DE}(\bgamma_T,\bgamma_T^d)(t))$ is obtained by applying a  rotation with angle $-\widetilde{\varphi}(t)$   to  $ (\bgamma_T(t),   \bgamma_T^d (t),\bgamma_T^h(t) )$ for all $t\in I$.
\end{proposition}
%%%%

Combining Propositions \ref{prop9} and  \ref{prop10}, we have the relations between  the de Sitter horocyclic and  de Sitter evolutes of $(\bgamma_T,\bgamma_T^d)$.
%%%%
\begin{theorem}\label{th3} 
Let $(\bgamma_T,\bgamma_T^d ):I\rightarrow\Delta_5$ be a timelike Legendre curve  with curvature  $(m_T,n_T)$.  \begin{enumerate} 
\item[$(1)$]  	Suppose that
there exists   $\widetilde{\varphi}:I\rightarrow\R$ such that $m_T(t)\cos\widetilde{\varphi}(t)+ n_T(t)\sin\widetilde{\varphi}(t)=0$ for all $t\in I$. If $\cos\widetilde{\varphi}(t)\neq0$ for all $t\in I$,  then  
$$ ( \widetilde{\mathcal{DE}}^\pm ,   \widetilde{\bmu}^\pm_\mathcal{E},\widetilde{\mathcal{HE}}^\pm)^T(\bgamma_T,\bgamma_T^d)       =  \overline{N}_1^\pm ( \tan\widetilde{\varphi} )  \overline{R} (\widetilde{\varphi} )(\widetilde{DE} , \widetilde{\bmu}_{DE},\widetilde{\bnu}_{DE}   )^T  (\bgamma_T,\bgamma_T^d)       . $$ 
\item[$(2)$] 	Suppose that there exists   $\widetilde{f}:I\rightarrow\R  $ such that $m_T(t)+\widetilde{f}(t)n_T(t)=0$ for all $t\in I$. Then	$$  
( \widetilde{DE} , \widetilde{\bmu}_{DE},\widetilde{\bnu}_{DE} )^T  (\bgamma_T,\bgamma_T^d) 	   = \overline{R} ( -\arctan\widetilde{f} ) \overline{N}_1^\pm(  -\widetilde{f}  )( \widetilde{\mathcal{DE}}^\pm,   \widetilde{\bmu}^\pm_\mathcal{E},\widetilde{\mathcal{HE}}^\pm   )^T  (\bgamma_T,\bgamma_T^d)    . $$ 
\end{enumerate} 	
\end{theorem}
%%%%

For curves with $(i,j)$-cusps in $\R^2$ and $H^2$, the  criteria
are available  in \cite{BG1982,HS2019,NTZ2026,P2001}, where $ (i, j)= (2, 3),(2, 5)$, $(3, 4),(3, 5) $.  Similar  criteria hold for curves in $S_1^2$. Here, we omit the details and proceed   to the application.
By  direct calculations, we obtain the following results.
%%%%%
\begin{proposition} 
Let $(\bgamma_T,\bgamma_T^d):I\rightarrow\Delta_5$ be a timelike Legendre curve  with curvature $(m_T,n_T)$. Suppose that there exists   $\widetilde{f}:I\rightarrow\R  $ such that $m_T(t)+\widetilde{f}(t)n_T(t)=0$ for all $t\in I$.  For any $t_0\in I$, the following assertions hold.
\begin{enumerate} 
\item[$(1)$]  $\widetilde{\mathcal{DE}}^\pm(\bgamma_T,\bgamma_T^d) $ is singular at $t_0$ if and only if  $  2\dot{\widetilde{f}}(t_0)\pm \widetilde{f}^2(t_0)n_T(t_0)=0 $.
\item[$(2)$]  $\widetilde{\mathcal{DE}}^\pm(\bgamma_T,\bgamma_T^d) $ has a $(2,3)$-cusp at $t_0$ if and only if  $n_T(t_0)\neq0$,  $  2\dot{\widetilde{f}}(t_0)\pm \widetilde{f}^2(t_0)n_T(t_0)=0 $ and $ 2\ddot{\widetilde{f}}(t_0)\pm 2\widetilde{f}(t_0)\dot{\widetilde{f}}(t_0)n_T(t_0)\pm \widetilde{f}^2(t_0)\dot{n}_T(t_0)\neq0$.
\item[$(3)$] $\widetilde{\mathcal{DE}}^\pm(\bgamma_T,\bgamma_T^d) $ has a  $(3,4)$-cusp at $t_0$ if and only if  $n_T(t_0)\neq0$, $2\dot{\widetilde{f}}(t_0)\pm \widetilde{f}^2(t_0)n_T(t_0)=2\ddot{\widetilde{f}}(t_0)\pm2\widetilde{f}(t_0)\dot{\widetilde{f}}(t_0)n_T(t_0)\pm \widetilde{f}^2(t_0)\dot{n}_T(t_0)=0$     and  $ \pm2\dot{\widetilde{f}^2}(t_0)n_T(t_0)\pm 2\widetilde{f}(t_0)\ddot{\widetilde{f}}(t_0)n_T(t_0)\pm 4\widetilde{f}(t_0)\dot{\widetilde{f}}(t_0)\dot{n}_T(t_0) \pm \widetilde{f}^2(t_0)\ddot{n}_T(t_0)+2\dddot{\widetilde{f}}(t_0)\neq0$.

\item[$(4)$] $\widetilde{\mathcal{DE}}^\pm(\bgamma_T,\bgamma_T^d) $ has a  $(2,5)$-cusp at $t_0$ if and only if  $n_T(t_0)=0$, $\dot{\widetilde{f}}(t_0)=0$, $  2\ddot{\widetilde{f}}(t_0)\pm \widetilde{f}^2(t_0)\dot{n}_T(t_0)\neq0 $ and $ -\ddot{\widetilde{f}}(t_0)\ddot{n}_T(t_0)+\dddot{\widetilde{f}}(t_0)\dot{n}_T(t_0)\neq0 $.

\item[$(5)$] $\widetilde{\mathcal{DE}}^\pm(\bgamma_T,\bgamma_T^d) $ has a  $(3,5)$-cusp at $t_0$ if and only if   $n_T(t_0)=\dot{\widetilde{f}}(t_0)=   2\ddot{\widetilde{f}}(t_0)\pm \widetilde{f}^2(t_0)\dot{n}_T(t_0)=0 $, $\dot{n}_T(t_0)\neq0$ and $ 2\dddot{\widetilde{f}}(t_0)\pm \widetilde{f}^2(t_0)\ddot{n}_T(t_0)\neq0$. 
\end{enumerate}
\end{proposition}
%%%%%
%%%%%
\begin{proposition}
Let $(\bgamma_T,\bgamma_T^d ):I\rightarrow\Delta_5$ be a timelike Legendre curve  with curvature  $(m_T,n_T)$. Suppose that there exists   $ \widetilde{\varphi}:I\rightarrow\R  $ such that   $m_T(t)\cos\widetilde{\varphi}(t)+ n_T(t)\sin\widetilde{\varphi}(t)=0$   for all $t\in I$.      For any $t_0\in I$, the following assertions hold.
\begin{enumerate} 
\item[$(1)$]  $ \widetilde{DE}(\bgamma_T,\bgamma_T^d) $ is singular at $t_0$ if and only if $ \dot{\widetilde{\varphi}}(t_0)=0 $.

\item[$(2)$]  $ \widetilde{DE}(\bgamma_T,\bgamma_T^d) $ has a $(2,3)$-cusp at $t_0$ if and only if  $ \dot{\widetilde{\varphi}}(t_0)=0 $, $ n_T(t_0)\neq0 $ and $ \ddot{\widetilde{\varphi}}(t_0)\neq0 $.

\item[$(3)$] $ \widetilde{DE}(\bgamma_T,\bgamma_T^d) $ has a  $(3,4)$-cusp at $t_0$ if and only if $ \dot{\widetilde{\varphi}}(t_0)=\ddot{\widetilde{\varphi}}(t_0)=0 $, $ n_T(t_0)\neq0 $ and $\dddot{\widetilde{\varphi}}(t_0)\neq0 $.

\item[$(4)$] $\widetilde{DE}(\bgamma_T,\bgamma_T^d) $ has a  $(2,5)$-cusp at $t_0$ if and only if  $ \dot{\widetilde{\varphi}}(t_0)=n_T(t_0)=0 $,  $\ddot{\widetilde{\varphi}}(t_0)\neq0 $ and   $\ddot{\widetilde{\varphi}}(t_0)\ddot{n}_T(t_0)-\dddot{\widetilde{\varphi}}(t_0)\dot{n}_T(t_0)\neq0 $.

\item[$(5)$] $ \widetilde{DE}(\bgamma_T,\bgamma_T^d)$ has a  $(3,5)$-cusp at $t_0$ if and only if  $ \dot{\widetilde{\varphi}}(t_0)=\ddot{\widetilde{\varphi}}(t_0)=n_T(t_0)=0 $,  $\dddot{\widetilde{\varphi}}(t_0)\neq0 $ and $\dot{n}_T(t_0)\neq0 $.
\end{enumerate}
\end{proposition}
%%%%%
%%%%% 
\begin{corollary}\label{cor2}
Let $(\bgamma_T,\bgamma_T^d ):I\rightarrow\Delta_5$ be a timelike Legendre curve  with curvature  $(m_T,n_T)$. If there exist   $\widetilde{f},\widetilde{\varphi}:I\rightarrow\R  $ such that $m_T(t)+\widetilde{f}(t)n_T(t)=0$ and $m_T(t)\cos\widetilde{\varphi}(t)+ n_T(t)\sin\widetilde{\varphi}(t)=0$ hold for all $t\in I$,   then $\widetilde{f}(t)=\tan \widetilde{\varphi}(t)$ for all $t\in I$ by the assumption that $\mathrm{Reg}(\bgamma_T^d)=\left\{t\in I\mid n_T(t)\neq0\right\}$  is   dense  in $ I $. Furthermore,  for any $t_0\in I$, the following assertions hold.
\begin{enumerate} 
\item[$(1)$] Suppose that $m_T(t_0)=0 $. Then $t_0$ is a singular point of $\widetilde{\mathcal{DE}}^\pm(\bgamma_T,\bgamma_T^d)$ if and only if $t_0$ is a singular point of $ \widetilde{DE}(\bgamma_T,\bgamma_T^d)$.
Especially,\\
{\rm (i)}   $t_0$ is a  $(2,3)$-cusp of $\widetilde{\mathcal{DE}}^\pm(\bgamma_T,\bgamma_T^d)$ if and only if $t_0$ is a $(2,3)$-cusp of $ \widetilde{DE}(\bgamma_T,\bgamma_T^d)$.\\
{\rm (ii)}   $t_0$ is a  $(3,4)$-cusp of $\widetilde{\mathcal{DE}}^\pm(\bgamma_T,\bgamma_T^d)$ if and only if $t_0$ is a $(3,4)$-cusp of $ \widetilde{DE}(\bgamma_T,\bgamma_T^d)$. 
\item[$(2)$] Suppose  that $m_T(t_0)\neq0 $. Then,\\
{\rm (i)} If $t_0$ is a singular point of $\widetilde{\mathcal{DE}}^\pm(\bgamma_T,\bgamma_T^d)$, then   $t_0$ is a regular point of $ \widetilde{DE}(\bgamma_T,\bgamma_T^d)$.\\
{\rm (ii)} If $t_0$ is a singular point of $ \widetilde{DE}(\bgamma_T,\bgamma_T^d)$, then   $t_0$ is a regular point of $\widetilde{\mathcal{DE}}^\pm(\bgamma_T,\bgamma_T^d)$.
\item[$(3)$] Suppose that $\dot{m}_T(t_0)=0 $. Then, \\
{\rm (i)}   $t_0$ is a  $(2,5)$-cusp of $\widetilde{\mathcal{DE}}^\pm(\bgamma_T,\bgamma_T^d)$ if and only if $t_0$ is a $(2,5)$-cusp of $ \widetilde{DE}(\bgamma_T,\bgamma_T^d)$. \\ 
{\rm (ii)} $t_0$ is a  $(3,5)$-cusp of $\widetilde{\mathcal{DE}}^\pm(\bgamma_T,\bgamma_T^d)$ if and only if $t_0$ is a $(3,5)$-cusp of $ \widetilde{DE}(\bgamma_T,\bgamma_T^d)$.
\end{enumerate}	 
\end{corollary}
%%%%%

\begin{figure}[h] 
\centering
\includegraphics[width = 14cm]{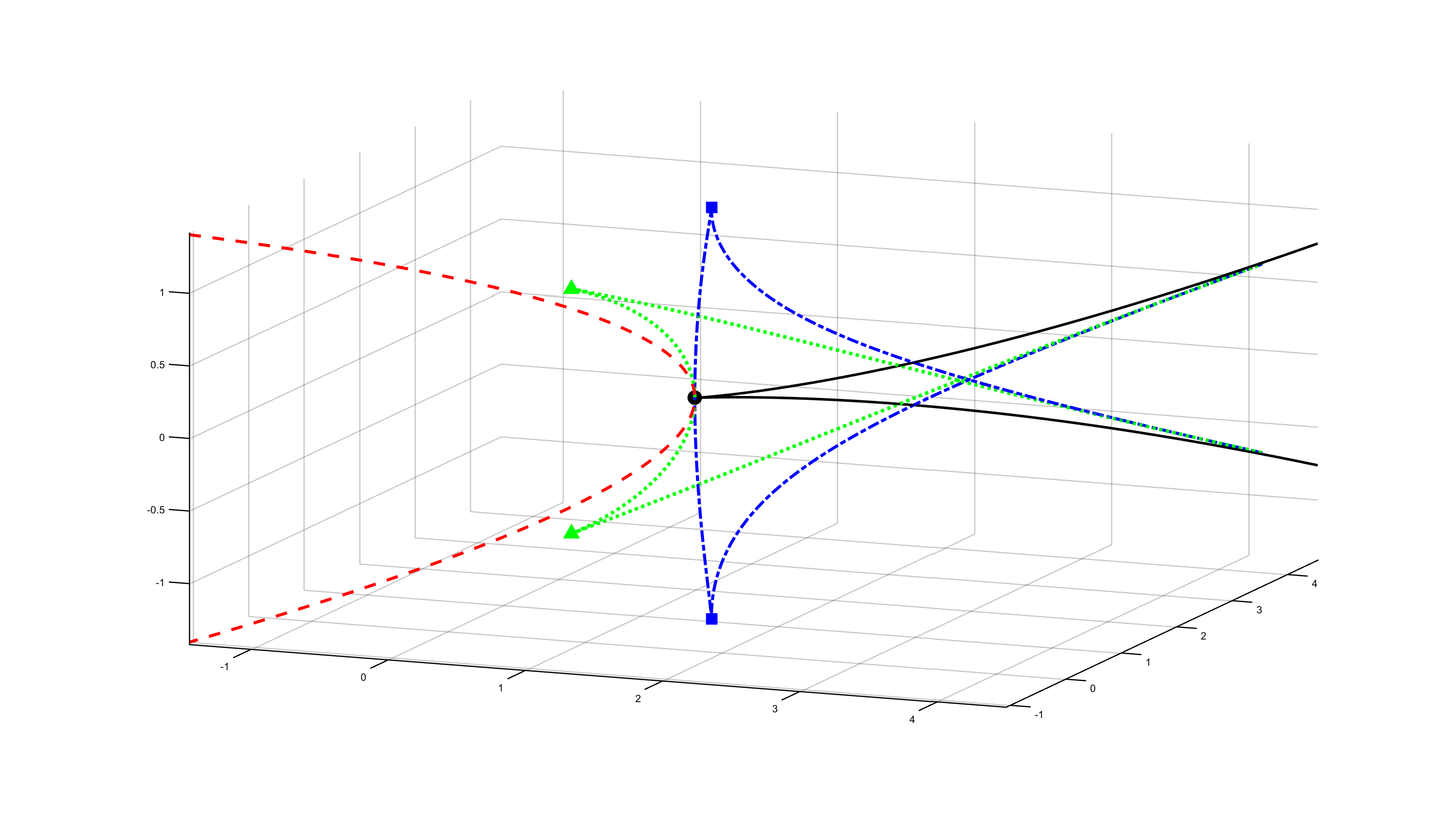}
\caption{The de Sitter  horocyclic evolutes $ \widetilde{\mathcal{DE}}^+(\bgamma_T,\bgamma_T^d) $ (red dashed curve), $\widetilde{\mathcal{DE}}^-(\bgamma_T,\bgamma_T^d)$ (blue dash-dotted curve) and  de Sitter  evolute $\widetilde{DE}  (\bgamma_T,\bgamma_T^d) $ (green dotted curve) of $(\bgamma_T,\bgamma_T^d)$. The black solid curve is the timelike frontal $\bgamma_T$.}
\label{figure4}
\end{figure}
\begin{example}\rm
Let $(\bgamma_T,\bgamma_T^d):(\R,0)\rightarrow\Delta_5$   be a timelike Legendre curve given by
$$\begin{aligned}
\bgamma_T(t) =&~
\Big( \sinh(4t^2+1),\;
\cosh(4t^2+1)\cos t^3,\;
\cosh(4t^2+1)\sin t^3
\Big),\\
\bgamma_T^d(t) =&~
\dfrac{1}{\sqrt{64 - 9t^2\cosh^2(4t^2+1)}}
\Big(
3t\cosh^2(4t^2+1),\\
&\quad 3t\sinh(4t^2+1)\cosh(4t^2+1)\cos t^3 - 8\sin t^3,\\
&\quad 3t\sinh(4t^2+1)\cosh(4t^2+1)\sin t^3 + 8\cos t^3
\Big). \end{aligned} $$
Then 
$$\begin{aligned}
\bgamma_T^h(t) =\bgamma_T(t)\wedge\bgamma_T^d(t)=&~
\dfrac{-1}{\sqrt{64 - 9t^2\cosh^2(4t^2+1)}}
\Big(
8\cosh(4t^2+1),\\
&\quad 8\sinh(4t^2+1)\cos t^3 - 3t\sin t^3\cosh(4t^2+1),\\
&\quad 8\sinh(4t^2+1)\sin t^3 + 3t\cos t^3\cosh(4t^2+1)
\Big)\end{aligned}$$ and the curvature $(m_T,n_T)$ of $(\bgamma_T,\bgamma_T^d) $ is given by
$$  
m_T(t)=-t\sqrt{64 - 9t^2\cosh^2(4t^2+1)}, $$ $$ 
n_T(t)=-\dfrac{3\left(9t^4\cosh^2(4t^2+1)\sinh(4t^2+1)
- 128t^2\sinh(4t^2+1) - 8\cosh(4t^2+1)\right)}
{9t^2\cosh^2(4t^2+1) - 64}
.$$
There exists a unique smooth function $\widetilde{f}:(\R,0)\rightarrow\R$, $$\widetilde{f}(t)=\dfrac{-t (64 - 9t^2\cosh^2(4t^2+1))^{\frac{3}{2}} }{3\left(9t^4\cosh^2(4t^2+1)\sinh(4t^2+1)
- 128t^2\sinh(4t^2+1) - 8\cosh(4t^2+1)\right)}  $$ such that $\widetilde{f}(t)n_T(t)+m_T(t)=0$ for all $t\in(\R,0)$. By Proposition \ref{prop9}, we have the de Sitter horocyclic evolute  of $(\bgamma_T,\bgamma_T^d)$ satisfying
$$\begin{aligned}
(\widetilde{\mathcal{DE}}^\pm,   \widetilde{\bmu}^\pm_\mathcal{E},\widetilde{\mathcal{HE}}^\pm )^T(\bgamma_T,\bgamma_T^d)  =\overline{N}_1^{\pm}( \widetilde{f} )  (\bgamma_T, \bgamma_T^d,\bgamma_T^h   ) ^T .\end{aligned}$$
Then  by Theorem \ref{th3} (2), the de Sitter  evolute  of $(\bgamma_T,\bgamma_T^d)$ satisfies $$( \widetilde{DE}, \widetilde{\bmu}_{DE},\widetilde{\bnu}_{DE} )^T  (\bgamma_T,\bgamma_T^d) 	   = \overline{R} ( -\arctan\widetilde{f} ) \overline{N}_1^\pm(  -\widetilde{f}  )( \widetilde{\mathcal{DE}}^\pm,   \widetilde{\bmu}^\pm_\mathcal{E},\widetilde{\mathcal{HE}}^\pm )^T  (\bgamma_T,\bgamma_T^d)     .$$

Since $m_T(0)=0$, $\dot{m}_T(0)\neq0$ and $n_T(0)\neq0$, $\bgamma_T$ has a $(2,3)$-cusp at $t_1=0$ by Lemma \ref{lemma1}. Since $2\dot{\widetilde{f}}(0)\pm \widetilde{f}^2(0)n_T(0)\neq0$,   $t_1=0$ is a regular point of $\widetilde{\mathcal{DE}}^\pm(\bgamma_T,\bgamma_T^d)$. It follows that  $t_1=0$ is also a regular point of $\widetilde{DE}  (\bgamma_T,\bgamma_T^d) $   by Corollary \ref{cor2} (1)  (black circle point in Figure \ref{figure4}). 

$\widetilde{\mathcal{DE}}^-(\bgamma_T,\bgamma_T^d)$ has  two singular points at $\pm t_2$, where $t_2\approx 0.18$ (blue square points in Figure \ref{figure4}). Since  $m_T(  t_2)\neq0$ and $m_T(-t_2)\neq0$, $\pm t_2$ are regular points of $\widetilde{DE}  (\bgamma_T,\bgamma_T^d) $ by Corollary \ref{cor2} (2)(i). 
$\widetilde{DE}  (\bgamma_T,\bgamma_T^d) $
has  two singular points at $\pm t_3$, where $t_3\approx 0.19$ (green triangle points in Figure \ref{figure4}). Since  $m_T(t_3)\neq0$ and $m_T(-t_3)\neq0$, $\pm t_3$ are regular points of $\widetilde{\mathcal{DE}}^\pm(\bgamma_T,\bgamma_T^d)$ by Corollary \ref{cor2} (2)(ii). 
\end{example}
\subsection{Involutes of timelike Legendre curves}
%%%%%
\begin{definition}[De Sitter and hyperbolic   involutes]\rm
Let $(\bgamma_T,\bgamma_T^d):I\rightarrow\Delta_5$ be a timelike Legendre curve  with curvature $(m_T,n_T)$ and  $t_0\in I$.
The \textit{de Sitter involute} $\widetilde{DI}_{t_0}(\bgamma_T,\bgamma_T^d):I\rightarrow S_1^2$ of    $(\bgamma_T,\bgamma_T^d)$  is defined as  
\begin{equation}\notag
\widetilde{DI}_{t_0}(\bgamma_T,\bgamma_T^d)(t)=  \cosh\left(\int_{t_0}^t{m_T(t)}{\rm d}t\right)\bgamma_T(t)-\sinh\left(\int_{t_0}^t{m_T(t)}{\rm d}t\right)\bgamma_T^h(t).
\end{equation} 
The \textit{hyperbolic involute}  $\widetilde{HI}_{t_0}(\bgamma_T,\bgamma_T^d):I\rightarrow H^2$  of $(\bgamma_T,\bgamma_T^d)$  is defined as  
\begin{equation}\notag
\widetilde{HI}_{t_0}(\bgamma_T,\bgamma_T^d)(t)=   \cosh\left(\int_{t_0}^t{m_T(t)}{\rm d}t\right)\bgamma_T^h(t)-\sinh\left(\int_{t_0}^t{m_T(t)}{\rm d}t\right)\bgamma_T(t).
\end{equation}  
\end{definition}
%%%%%

%%%%
\begin{proposition} 
Let   $\widetilde{DI}_{t_0}(\bgamma_T,\bgamma_T^d)$ and $ \widetilde{HI}_{t_0}(\bgamma_T,\bgamma_T^d) $ be de Sitter and hyperbolic   involutes of  $(\bgamma_T,\bgamma_T^d)$, respectively. Then
\begin{enumerate} 
\item[$(1)$]   
$( \widetilde{DI}_{t_0}(\bgamma_T,\bgamma_T^d), \widetilde{HI}_{t_0}(\bgamma_T,\bgamma_T^d)  ):I\rightarrow \overline{\Delta}_1$ is a spacelike Legendre curve with curvature 
$$
\left(\widetilde{m} _{DI}(t), \widetilde{n}_{DI}(t)\right)= \left(     -n_T(t)\sinh\left(\int_{t_0}^t{m_T(t)}{\rm d}t\right),  n_T(t)\cosh\left(\int_{t_0}^t{m_T(t)}{\rm d}t\right)\right).
$$

\item[$(2)$]    
$( \widetilde{HI}_{t_0}(\bgamma_T,\bgamma_T^d), \widetilde{DI}_{t_0}(\bgamma_T,\bgamma_T^d) ):I\rightarrow \Delta_1$ is a spacelike Legendre curve with curvature 
$$
\left(\widetilde{m} _{HI}(t), \widetilde{n}_{HI}(t)\right)= \left(     -n_T(t)\cosh\left(\int_{t_0}^t{m_T(t)}{\rm d}t\right), n_T(t)\sinh\left(\int_{t_0}^t{m_T(t)}{\rm d}t\right)\right). 
$$ 
\end{enumerate}
\end{proposition}
%%%%

%%%%
\begin{proposition}   
 Let $\widetilde{\bmu}_{I}(\bgamma_T,\bgamma_T^d)=\widetilde{DI}_{t_0}(\bgamma_T,\bgamma_T^d)\wedge\widetilde{HI}_{t_0}(\bgamma_T,\bgamma_T^d)$ and denote the moving frame of $\widetilde{DI}_{t_0}(\bgamma_T,\bgamma_T^d)$ by $ ( \widetilde{DI}_{t_0}, \widetilde{\bmu}_{I},\widetilde{HI}_{t_0}   ) (\bgamma_T,\bgamma_T^d)$.	Then  
$$\begin{aligned} 
( \widetilde{DI}_{t_0}, \widetilde{\bmu}_{I},\widetilde{HI}_{t_0}  )^T(\bgamma_T,\bgamma_T^d)=	\overline{B}_1\left(- \int_{t_0}^t{m_T }{\rm d}t\right) (  \bgamma_T,\bgamma_T^d, \bgamma_T^h  )^T.\end{aligned}$$
That is,	$ ( \widetilde{DI}_{t_0}(\bgamma_T,\bgamma_T^d)(t), \widetilde{\bmu}_{I}(\bgamma_T,\bgamma_T^d)(t),\widetilde{HI}_{t_0}(\bgamma_T,\bgamma_T^d) (t) )$  is obtained by applying a boost   in the $\bgamma_T $-direction with rapidity $ -\int_{t_0}^t{m_T(t)}{\rm d}t$ to $ (\bgamma_T(t),\bgamma_T^d(t), \bgamma_T^h(t) )$  for all $t\in I$.
\end{proposition}
%%%%
\subsection{Parallels of timelike Legendre curves}
%%%%%
\begin{definition}[De Sitter horocyclic parallel]  \rm
Let $(\bgamma_T,\bgamma_T^d):I\rightarrow\Delta_5$ be a timelike Legendre curve  with curvature $(m_T,n_T)$.
The   \textit{de Sitter horocyclic parallel} $\widetilde{\mathcal{DP}}^\pm_{\widetilde{\lambda}_\pm} (\bgamma_T,\bgamma_T^d):I\rightarrow S_1^2$ of   $(\bgamma_T,\bgamma_T^d)$  is   defined as
$$
\widetilde{\mathcal{DP}}^\pm_{\widetilde{\lambda}_\pm} (\bgamma_T,\bgamma_T^d)(t)=\left(1-\frac{\widetilde{\lambda}^2 _\pm(t)}{2} \right) \bgamma_T(t)+\widetilde{\lambda} _\pm(t)\bgamma_T^d(t)\pm \frac{\widetilde{\lambda}^2_\pm(t)}{2} \bgamma_T^h(t),
$$
where 
$ \widetilde{\lambda}_\pm(t)  $ is a solution of the Bernoulli equation \begin{equation}\label{Bernoulli2}
\frac{{\rm d}\widetilde{\lambda} _\pm}{{\rm d}t}(t)=\pm \left(\frac{ \widetilde{\lambda} _\pm^2(t) n_T(t)}{2}+\widetilde{\lambda} _\pm(t)m_T(t)\right).
\end{equation}
\end{definition}
%%%%%

%%%%%
\begin{proposition}  
Let  $\widetilde{\mathcal{DP}}^\pm_{\widetilde{\lambda}_\pm} (\bgamma_T,\bgamma_T^d)$    be  a de Sitter horocyclic parallel    of $(\bgamma_T,\bgamma_T^d)$. Then    $(\widetilde{\mathcal{DP}}^\pm_{\widetilde{\lambda}_\pm} (\bgamma_T,\bgamma_T^d),\widetilde{\bnu}_\mathcal{DP}^\pm(\bgamma_T,\bgamma_T^d)):I\rightarrow\Delta_5$ is a timelike Legendre curve with curvature 
$$   (\widetilde{m}_\mathcal{DP}^\pm(t),\widetilde{n}_\mathcal{DP}^\pm(t) )= (\widetilde{\lambda} _\pm(t)n_T(t)+m_T(t),n_T(t) ), $$   
where $\widetilde{\bnu}_\mathcal{DP}^\pm(\bgamma_T,\bgamma_T^d)(t)=-\widetilde{\lambda} _\pm(t)\bgamma_T(t)+\bgamma_T^d(t)\pm\widetilde{\lambda} _\pm(t)\bgamma_T^h(t)$.
\end{proposition}
%%%% 
\begin{corollary}
The image of the de Sitter horocyclic evolute $\widetilde{\mathcal{DE}}^\pm(\bgamma_T,\bgamma_T^d)$ is contained in the set of singular values of the de Sitter horocyclic parallels $\widetilde{\mathcal{DP}}^\pm_{\widetilde{\lambda}_\pm} (\bgamma_T,\bgamma_T^d)$ of $(\bgamma_T,\bgamma_T^d)$. 	 
\end{corollary}
%%%%

%%%%
\begin{proposition} \label{prop11} 
  Let  $\widetilde{\bmu}_\mathcal{DP}^\pm(\bgamma_T,\bgamma_T^d)=\widetilde{\mathcal{DP}}^\pm_{\widetilde{\lambda}_\pm} (\bgamma_T,\bgamma_T^d)\wedge \widetilde{\bnu}_\mathcal{DP}^\pm(\bgamma_T,\bgamma_T^d)$ and denote the moving frame of $\widetilde{\mathcal{DP}}^\pm_{\widetilde{\lambda}_\pm} (\bgamma_T,\bgamma_T^d)$ by $ (	\widetilde{\mathcal{DP}}^\pm_{\widetilde{\lambda}_\pm} , \widetilde{\bnu}_\mathcal{DP}^\pm,	 \widetilde{\bmu}_\mathcal{DP}^\pm   ) (\bgamma_T,\bgamma_T^d)$. Then   
$$(	\widetilde{\mathcal{DP}}^\pm_{\widetilde{\lambda}_\pm} , \widetilde{\bnu}_\mathcal{DP}^\pm,	 \widetilde{\bmu}_\mathcal{DP}^\pm   )^T(\bgamma_T,\bgamma_T^d)
=\overline{N}_1^\pm( \widetilde{\lambda}_\pm )(	\bgamma_T, \bgamma_T^d,\bgamma_T^h  ) ^T.$$
That is,	$ ( \widetilde{\mathcal{DP}}^\pm_{\widetilde{\lambda}_\pm} (\bgamma_T,\bgamma_T^d)(t), \widetilde{\bnu}_\mathcal{DP}^\pm(\bgamma_T,\bgamma_T^d)(t),	 \widetilde{\bmu}_\mathcal{DP}^\pm(\bgamma_T,\bgamma_T^d)(t)   )$
is obtained by applying a null rotation  around  the $(\bgamma_T\pm \bgamma_T^h)$-direction with   angle $  \widetilde{\lambda}_\pm(t)$ to $ ( 	\bgamma_T(t), \bgamma_T^d(t),\bgamma_T^h(t)   )$ for all $t\in I$.  
\end{proposition}
%%%%

%%%%%
\begin{definition}[De Sitter  parallel]  \rm
Let $(\bgamma_T,\bgamma_T^d):I\rightarrow\Delta_5$ be a timelike Legendre curve  with curvature $(m_T,n_T)$ and $\theta\in\R$. 
The \textit{de Sitter parallel}  $\widetilde{DP }_\theta (\bgamma_T, \bgamma_T^d):I\rightarrow S_1^2$ of $(\bgamma_T,\bgamma_T^d)$ is   defined as
$$
\widetilde{DP }_\theta (\bgamma_T, \bgamma_T^d)(t)=\cos \theta\bgamma_T(t)+\sin\theta\bgamma_T^d(t).$$
\end{definition}
%%%%%

%%%%
\begin{proposition} 
Let $\widetilde{DP }_\theta (\bgamma_T, \bgamma_T^d)$   be a  de Sitter  parallel  of $(\bgamma_T,\bgamma_T^d)$.  Then \\
$ ( \widetilde{DP }_\theta (\bgamma_T,\bgamma_T^d), \widetilde{\bnu}_{DP} (\bgamma_T,\bgamma_T^d) ):I\rightarrow \Delta_5$ is a timelike Legendre curve with curvature 
$$
\left(	\widetilde{m} _{DP}(t), 	\widetilde{n}_{DP}(t)\right)= \left(      m_T(t)\cos  \theta+ n_T(t)\sin \theta,   n_T(t)\cos  \theta+ m_T(t)\sin \theta\right),
$$ 
where $\widetilde{\bnu}_{DP}(\bgamma_T,\bgamma_T^d)(t)=\cos \theta\bgamma_T^d(t)-\sin\theta\bgamma_T(t)$.	
\end{proposition}
%%%%

%%%% 
\begin{corollary}
The image of the de Sitter   evolute $ \widetilde{DE} (\bgamma_T, \bgamma_T^d)$ is contained in the set of singular values of the de Sitter  parallels $\widetilde{DP}_\theta(\bgamma_T, \bgamma_T^d	)$ of $(\bgamma_T,\bgamma_T^d)$. 	 
\end{corollary}
%%%%

\begin{proposition} \label{prop12} 
 Let  $\widetilde{\bmu}_{DP}(\bgamma_T,\bgamma_T^d)=\widetilde{DP }_\theta (\bgamma_T, \bgamma_T^d)\wedge \widetilde{\bnu}_{DP}(\bgamma_T,\bgamma_T^d)$ and denote the moving frame of $\widetilde{DP }_\theta (\bgamma_T, \bgamma_T^d)$ by $ ( 	\widetilde{DP }_\theta ,\widetilde{\bnu}_{DP},  \widetilde{\bmu}_{DP}    ) (\bgamma_T,\bgamma_T^d)$. Then      
$$( 	\widetilde{DP }_\theta ,\widetilde{\bnu}_{DP},  \widetilde{\bmu}_{DP}   )^T(\bgamma_T,\bgamma_T^d)
=\overline{R} (-\theta) (	\bgamma_T, \bgamma_T^d,\bgamma_T^h  )^T .$$
That is,	$ ( \widetilde{DP }_\theta (\bgamma_T, \bgamma_T^d)(t),\widetilde{\bnu}_{DP}(\bgamma_T,\bgamma_T^d)(t),  \widetilde{\bmu}_{DP}(\bgamma_T,\bgamma_T^d)(t)   )    $ is obtained by applying a  rotation  with angle   $-\theta$  to $ ( 	\bgamma_T(t), \bgamma_T^d(t),\bgamma_T^h(t)  )  $ for all $t\in I$.   
\end{proposition}
%%%%

Combining Propositions \ref{prop11} and  \ref{prop12}, we have the relations between  the de Sitter horocyclic and  de Sitter parallels of $(\bgamma_T,\bgamma_T^d)$.
%%%%
\begin{theorem} 
	Let $(\bgamma_T,\bgamma_T^d):I\rightarrow\Delta_5$ be a timelike Legendre curve  with curvature $(m_T,n_T)$.    $\widetilde{\mathcal{DP}}^\pm_{\widetilde{\lambda}_\pm} (\bgamma_T,\bgamma_T^d)$  and
	$\widetilde{DP }_\theta (\bgamma_T, \bgamma_T^d)$   are  de Sitter  horocyclic and  de Sitter parallels of $(\bgamma_T,\bgamma_T^d)$, respectively. Then 
	$$\begin{aligned}
		(	\widetilde{\mathcal{DP}}^\pm_{\widetilde{\lambda}_\pm} , \widetilde{\bnu}_\mathcal{DP}^\pm,	 \widetilde{\bmu}_\mathcal{DP}^\pm  )^T (\bgamma_T,\bgamma_T^d)  &=\overline{N}_1^\pm( \widetilde{ \lambda}_\pm )\overline{R}( \theta)( 	\widetilde{DP }_\theta  ,\widetilde{\bnu}_{DP},  \widetilde{\bmu}_{DP}   )  ^T  (\bgamma_T,\bgamma_T^d),\\
		( \widetilde{DP }_\theta  ,\widetilde{\bnu}_{DP},  \widetilde{\bmu}_{DP}  )^T (\bgamma_T,\bgamma_T^d)  &=\overline{R}( -\theta)\overline{N}_1^\pm( - \widetilde{\lambda}_\pm ) (	\widetilde{\mathcal{DP}}^\pm_{\widetilde{\lambda}_\pm}  , \widetilde{\bnu}_\mathcal{DP}^\pm,	 \widetilde{\bmu}_\mathcal{DP}^\pm  )^T (\bgamma_T,\bgamma_T^d)    
		.\end{aligned}$$ 
\end{theorem}
%%%%
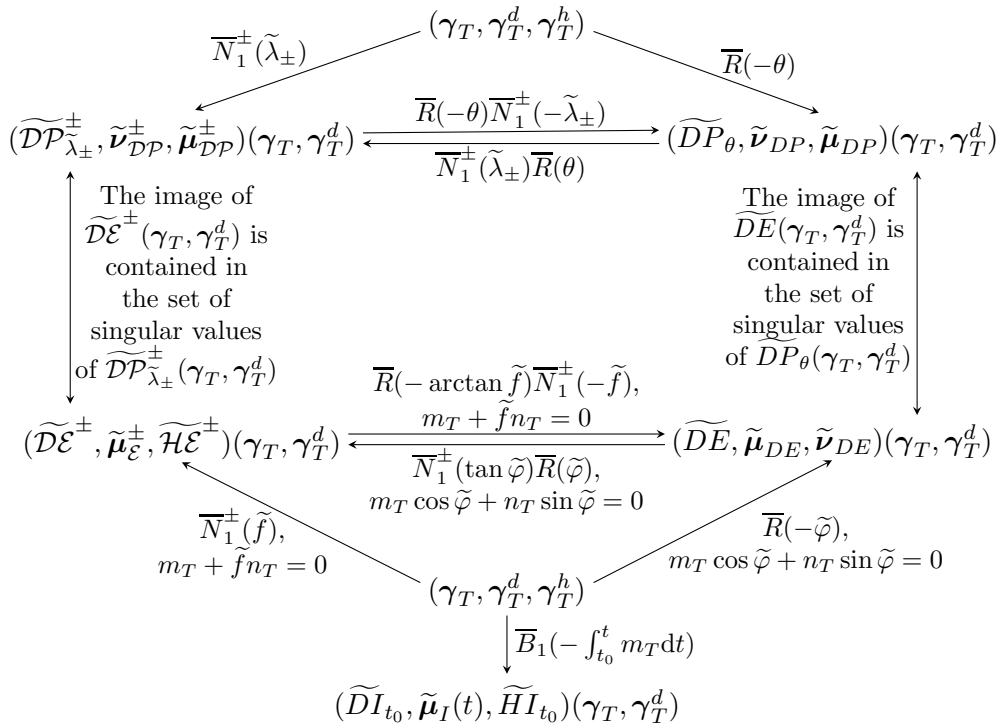
\begin{figure}[h]  
\centering
\begin{tikzpicture}[
>=stealth,  
node distance=3cm,
vertex/.style={circle, minimum size=6mm, inner sep=0pt, font= \small}
]

\node[vertex] (D1) at (-4.3, 3) {$ (\widetilde{\mathcal{DE}}^\pm,   \widetilde{\bmu}^\pm_\mathcal{E},\widetilde{\mathcal{HE}}^\pm     )(\bgamma_T,\bgamma_T^d)$};
\node[vertex] (D2) at (4.3, 3) {$  ( \widetilde{DE}, \widetilde{\bmu}_{DE},\widetilde{\bnu}_{DE}     ) (\bgamma_T,\bgamma_T^d)$};
\node[vertex] (D3) at (-4.3, 7) {$ (  \widetilde{\mathcal{DP}}^\pm_{\widetilde{\lambda}_\pm} , \widetilde{\bnu}_\mathcal{DP}^\pm,	 \widetilde{\bmu}_\mathcal{DP}^\pm    )(\bgamma_T,\bgamma_T^d)$};
\node[vertex] (D4) at (4.3, 7) {$  (\widetilde{DP }_\theta ,\widetilde{\bnu}_{DP},  \widetilde{\bmu}_{DP}    )(\bgamma_T, \bgamma_T^d) $};
\node[vertex] (D5) at (0,1) {$ ( 	\bgamma_T, \bgamma_T^d,\bgamma_T^h    )$};
\node[vertex] (D0) at (0,8.5) {$ ( 	\bgamma_T, \bgamma_T^d,\bgamma_T^h     )$};
\node[vertex] (D6) at (0,-0.5) {$ (   \widetilde{DI}_{t_0}, \widetilde{\bmu}_{I}(t),\widetilde{HI}_{t_0} )(\bgamma_T,\bgamma_T^d)   $};

% 0 -> 5  
\draw[->] ([yshift=-40pt]D5.north) -- ([yshift=77pt]D6.south) 
node[midway, xshift=37pt, font=\footnotesize, ] {$\overline{B}_1(- \int_{t_0}^t{m_T }{\rm d}t)$ };

% 1 <-  2  
\draw[->] ([xshift=-2pt, yshift=-2 pt]D2.west) -- ([xshift=2pt, yshift=-2 pt]D1.east) 
node[midway,   yshift=-14pt,font=\footnotesize,align=center] { $\overline{N}_1^\pm (\tan\widetilde{\varphi} )  \overline{R} (\widetilde{\varphi} )$,\\$m_T \cos\widetilde{\varphi} + n_T \sin\widetilde{\varphi} =0$ };
\draw[<-] ([xshift=-2pt, yshift=2 pt]D2.west) -- ([xshift=2pt, yshift=2 pt]D1.east) 
node[midway,   yshift=14pt,font=\footnotesize,align=center] {$ \overline{R} ( -\arctan\widetilde{f}) \overline{N}_1^\pm( - \widetilde{f} ) $,\\ $m_T +\widetilde{f} n_T =0$ };

% 3 <-> 4  
\draw[->] ([xshift=2pt, yshift=1.5pt]D3.east) -- ([xshift=-2pt, yshift=2.5pt]D4.west) 
node[midway,   yshift=8pt,font=\footnotesize,align=center] {$ \overline{R}( -\theta)\overline{N}_1^\pm( - \widetilde{\lambda}_\pm)$};
\draw[->] ([xshift=-2pt, yshift=-2.5pt]D4.west) -- ([xshift=2pt, yshift=-2.5pt]D3.east) 
node[midway,   yshift=-8pt,font=\footnotesize,align=center] { $\overline{N}_1^\pm( \widetilde{ \lambda}_\pm)\overline{R}( \theta)$};

% 5  -> 1
\draw[->] ([xshift=-2pt, yshift=3pt]D5.west) -- ([xshift=0pt, yshift=52pt]D1.south) 
node[midway,font=\footnotesize,align=center, xshift=-22pt,yshift=-8pt] {$\overline{N}_1^{\pm}( \widetilde{f}) $,\\$m_T +\widetilde{f} n_T =0$};

% 5  -> 2
\draw[->] ([xshift=2pt, yshift=3pt]D5.east) --([xshift=0pt, yshift=53pt]D2.south) 
node[midway,font=\footnotesize,align=center,  xshift=35pt,yshift=-8pt] {$\overline{R}(-\widetilde{\varphi})  $,\\$m_T \cos\widetilde{\varphi} + n_T \sin\widetilde{\varphi} =0$};

%	0  -> 3 
\draw[->] ([xshift=-2pt, yshift=-3pt]D0.west) -- ([xshift=5pt, yshift=-56pt]D3.north) 
node[midway,font=\footnotesize,align=center, xshift=-18pt,yshift=8pt] {$\overline{N}_1^\pm(  \widetilde{\lambda}_\pm) $};
%	0  -> 4 
\draw[->] ([xshift=2pt, yshift=-3pt]D0.east) -- ([xshift=-5pt, yshift=-54pt]D4.north) 
node[midway,font=\footnotesize,align=center, xshift=20pt,yshift=3pt] {$\overline{R}( -\theta)$};

% 1 <-> 3 
\draw[<- >] ([xshift=-42pt, yshift=-79pt]D3.north) -- ([xshift=-42pt, yshift=74pt]D1.south) 
node[midway,font=\footnotesize,align=center, xshift=40pt,yshift=0pt] {The image of\\$\widetilde{\mathcal{DE}}^\pm(\bgamma_T,\bgamma_T^d)$ is\\ contained in \\the set of \\singular values \\of   $\widetilde{\mathcal{DP}}_{\widetilde{\lambda}_\pm}^\pm (\bgamma_T,\bgamma_T^d)$ };

\draw[<-> ] ([xshift=32pt, yshift=-76pt]D4.north) -- ([xshift=32pt, yshift=70pt]D2.south) 
node[midway,font=\footnotesize,align=center, xshift= -37pt,yshift=5pt] {The image of\\$\widetilde{ {DE}} (\bgamma_T,\bgamma_T^d)$ is\\ contained in \\the set of \\singular values \\of   $\widetilde{ {DP}}_{\theta}(\bgamma_T,\bgamma_T^d)$ };

\pgfresetboundingbox
\path (-6, -1 ) rectangle (6, 9.2); 

\end{tikzpicture}
\caption{The relations among the evolutes, involutes and parallels of $(\bgamma_T,\bgamma_T^d)$.}  
\label{figure5} 
\end{figure}
%%%%

Based on the above conclusions, we present the relations among the evolutes, involutes and parallels of $(\bgamma_T,\bgamma_T^d)$, as illustrated in Figure \ref{figure5}.

We give  relations between the evolutes and   involutes of $(\bgamma_h,\bgamma_h^d)$, $(\bgamma_d,\bgamma_d^h)$ and $(\bgamma_T,\bgamma_T^d)$.
\begin{theorem} \label{th1}
Using the notations introduced above, we obtain the following results.
\begin{enumerate} 
\item[$(1)$]  Suppose that there exists   $f:I\rightarrow\R  $ such that $m_h(t)+f(t)n_h(t)=0$ for all $t\in I$ and $\mathrm{Reg}(\bnu_{\mathcal{HI}}^\pm )=\left\{t\in I\mid n_h(t)\pm m_h(t)\neq0\right\}$  is  dense in $I$. Then    $$
\begin{aligned}
&(\mathcal{HE}^\pm, \mathcal{DE}^\pm ,\bmu^\pm_\mathcal{E}   ) (\mathcal{HI}_{s_\pm}^\pm(\bgamma_h, \bgamma_h^d ),\bnu_{\mathcal{HI}}^\pm(\bgamma_h, \bgamma_h^d ) ) =  (\bgamma_h, \bgamma_h^d,\bgamma_h^s   )   	,\\
&(\mathcal{HE}^\pm, \mathcal{DE}^\pm,\bmu^\pm_\mathcal{E} ) (\mathcal{HI}_{s_\mp}^\mp(\bgamma_h, \bgamma_h^d ),-\bnu_{\mathcal{HI}}^\mp(\bgamma_h, \bgamma_h^d ))   =  (\bgamma_h, -\bgamma_h^d,-\bgamma_h^s   )   	,\\
&(\mathcal{HI}_{S ^\pm_\mathcal{HE}}^\pm , \bnu_{\mathcal{HI}}^\pm ,\bmu^\pm_\mathcal{HI}  )(\mathcal{HE}^\pm(\bgamma_h, \bgamma_h^d ),\mathcal{DE}^\pm(\bgamma_h, \bgamma_h^d ))      =   (\mathcal{HP}_{f-S ^\pm_\mathcal{HE}}^\pm , \bnu_{\mathcal{HP}}^\pm  ,	 \bmu_\mathcal{HP}^\pm   )(\bgamma_h, \bgamma_h^d)   	, \\
&(\mathcal{HI}_{ \overline{S}  ^\pm_\mathcal{HE}}^\pm , \bnu_{\mathcal{HI}}^\pm ,\bmu^\pm_\mathcal{HI}  )(\mathcal{HE}^\mp(\bgamma_h, \bgamma_h^d),-\mathcal{DE}^\mp(\bgamma_h, \bgamma_h^d))   =    (\mathcal{HP}_{f+\overline{S} ^\pm_\mathcal{HE}}^\mp   , -\bnu_{\mathcal{HP}}^\mp   ,	 -\bmu_\mathcal{HP}^\mp ) (\bgamma_h, \bgamma_h^d)  	, 
\end{aligned}
$$
where $  \overline{S}  ^\pm_\mathcal{HE}(t) =-S ^\mp_\mathcal{HE}(t)$ and $S ^\pm_\mathcal{HE}(t)$ is a solution of the Riccati equation $$
\begin{aligned}\frac{{\rm d} {S}_\mathcal{HE}^\pm(t)}{{\rm d}t} = \frac{m_\mathcal{HE}^\pm(t)\pm n_\mathcal{HE}^\pm(t)}{2}(S_\mathcal{HE}^\pm)^2(t)-m_\mathcal{HE}^\pm(t)  .
\end{aligned}$$
\item[$(2)$] Suppose that there exists   $\overline{f}:I\rightarrow\R  $ such that $m_d(t)+\overline{f}(t)n_d(t)=0$ for all $t\in I$ and $\mathrm{Reg}(\overline{\bnu}_\mathcal{DI}^\pm )=\left\{t\in I\mid n_d(t)\mp m_d(t)\neq0\right\}$  is  dense in $I$. Then 
$$\begin{aligned}
&(\widetilde{\mathcal{DE}}^\pm , \widetilde{\mathcal{HE}}^\pm ,\widetilde{\bmu}^\pm_\mathcal{E}  )(\overline{\mathcal{DI}}_{\overline{s}_\pm}^\pm(\bgamma_d, \bgamma_d^h ),\overline{\bnu}^\pm_\mathcal{DI}(\bgamma_d, \bgamma_d^h ))   =   (\bgamma_d, \bgamma_d^h,\bgamma_d^s   )   	,\\
&(\widetilde{\mathcal{DE}}^\pm , \widetilde{\mathcal{HE}}^\pm ,\widetilde{\bmu}^\pm_\mathcal{E}  )(\overline{\mathcal{DI}}_{\overline{s}_\mp}^\mp (\bgamma_d, \bgamma_d^h ),-\overline{\bnu}^\mp_\mathcal{DI}(\bgamma_d, \bgamma_d^h ))   =  (\bgamma_d, -\bgamma_d^h,-\bgamma_d^s   )  	,\\
&(\overline{\mathcal{DI}}_{\overline{S} ^\pm_\mathcal{DE}}^\pm , \overline{\bnu}^\pm_\mathcal{DI} ,\overline{\bmu}^\pm_\mathcal{DI}   )(\widetilde{\mathcal{DE}}^\pm(\bgamma_T,\bgamma_T^d),\widetilde{\mathcal{HE}}^\pm(\bgamma_T,\bgamma_T^d))  =     (\widetilde{\mathcal{DP}}_{\overline{f}-\overline{S} ^\pm_\mathcal{DE}}^\pm  , \widetilde{\bnu}_\mathcal{DP}^\pm  ,	 \widetilde{\bmu}_\mathcal{DP}^\pm   ) (\bgamma_T,\bgamma_T^d)  	,\\
&(\overline{\mathcal{DI}}_{\overline{\overline{S}} ^\pm_\mathcal{DE}}^\pm  , \overline{\bnu}^\pm_\mathcal{DI} ,\overline{\bmu}^\pm_\mathcal{DI}   )(\widetilde{\mathcal{DE}}^\mp(\bgamma_T,\bgamma_T^d),-\widetilde{\mathcal{HE}}^\mp(\bgamma_T,\bgamma_T^d)) =   (\widetilde{\mathcal{DP}}_{\overline{f}+\overline{\overline{S}} ^\pm_\mathcal{DE}}^\mp , -\widetilde{\bnu}_\mathcal{DP}^\mp ,	 -\widetilde{\bmu}_\mathcal{DP}^\mp  )(\bgamma_T,\bgamma_T^d)   , 
\end{aligned}$$ 
where $ \overline{\overline{S}} ^\pm_\mathcal{DE}(t) =-\overline{S} ^\mp_\mathcal{DE}(t)$ and $\overline{S} ^\pm_\mathcal{DE}(t)$ is a solution of the Riccati equation
\begin{equation}\notag
\frac{{\rm d}\overline{S} ^\pm_\mathcal{DE}}{{\rm d}t}(t)= \frac{\widetilde{m}^\pm_\mathcal{DE}(t)\pm \widetilde{n}^\pm_\mathcal{DE}(t)}{2}(\overline{S} ^\pm_\mathcal{DE})^2(t)-\widetilde{m}^\pm_\mathcal{DE}(t).\end{equation}
\item[$(3)$] Suppose	that   there exists   $\varphi:I\rightarrow\R$ such that $m_h(t)\cosh\varphi(t)+ n_h(t)\sinh\varphi(t)=0$ for all $t\in I$. Then  
$$\begin{aligned}& (HE, \bnu_{HE},\bmu_{HE}    ) (HI_{t_0}(\bgamma_h,\bgamma_h^d),DI_{t_0}(\bgamma_h,\bgamma_h^d) )  =  ( \bgamma_h,\bgamma_h^d,\bgamma_h^s   ),\\ 
& ( HI_{t_0}, DI_{t_0},\bmu_I    ) (HE(\bgamma_h,\bgamma_h^d),\bnu_{HE}(\bgamma_h,\bgamma_h^d) )=  (  	HP_{\varphi(t_0)}, DP_{\varphi(t_0)},	\bmu_{P}  )(\bgamma_h,\bgamma_h^d).
\end{aligned}$$

\item[$(4)$] Suppose	that there exists  $\overline{\varphi}:I\rightarrow\R$ such that $m_d(t)\sinh\overline{\varphi}(t)+ n_d(t)\cosh\overline{\varphi}(t)=0$ for all $t\in I$. Then 
$$\begin{aligned} &(\overline{HE} , \overline{\bnu}_{HE} ,\overline{\bmu}_{HE}   ) (DI_{t_0}(\bgamma_h,\bgamma_h^d),HI_{t_0}(\bgamma_h,\bgamma_h^d))= ( \bgamma_h,\bgamma_h^d,\bgamma_h^s   ),\\ 
&( DI_{t_0}, HI_{t_0} ,\bmu_I    ) (\overline{HE}(\bgamma_d,\bgamma_d^h),\overline{\bnu}_{HE}(\bgamma_d,\bgamma_d^h))= ( -\overline{DP}_{\overline{\varphi}(t_0)}, 	\overline{HP}_{\overline{\varphi}(t_0)}, 	-\overline{\bmu}_{P} )(\bgamma_d,\bgamma_d^h).
\end{aligned}$$ 
\item[$(5)$] Suppose that  there exists   $\overline{\psi}:I\rightarrow\R$ such that $m_d(t)\cosh\overline{\psi}(t)+ n_d(t)\sinh\overline{\psi}(t)=0$ for all $t\in I$.  Then 
$$\begin{aligned}&(\overline{DE}, \overline{\bnu}_{HE},\overline{\bmu}_{HE}  ) (\widetilde{DI}_{t_0}(\bgamma_T,\bgamma_T^d ),\widetilde{HI}_{t_0}(\bgamma_T,\bgamma_T^d )) = ( \bgamma_T,\bgamma_T^d,\bgamma_T^h   ),\\ 
&( \widetilde{DI}_{t_0}, \widetilde{HI}_{t_0} ,\widetilde{\bmu}_{I}   ) (\overline{DE}(\bgamma_d,\bgamma_d^h),\overline{\bnu}_{DE}(\bgamma_d,\bgamma_d^h)) = (  \overline{DP}_{\overline{\psi}(t_0)}, 	\overline{HP}_{\overline{\psi}(t_0)}, 	 \overline{\bmu}_{P} )(\bgamma_d,\bgamma_d^h).
\end{aligned}$$   
\item[$(6)$] Suppose	that  there exists   $\psi:I\rightarrow\R$ such that $m_h(t)\sinh\psi(t)+ n_h(t)\cosh\psi(t)=0$ for all $t\in I$.   Then 
$$\begin{aligned}&( {DE},  {\bnu}_{HE}, {\bmu}_{HE}   ) (\widetilde{HI}_{t_0}(\bgamma_T,\bgamma_T^d ),\widetilde{DI}_{t_0}(\bgamma_T,\bgamma_T^d )) = ( \bgamma_T,\bgamma_T^d,\bgamma_T^h   ),\\ 
&( \widetilde{HI}_{t_0} , \widetilde{DI}_{t_0},\widetilde{\bmu}_{I}  ) ( {DE}(\bgamma_h,\bgamma_h^d), {\bnu}_{DE}(\bgamma_h,\bgamma_h^d))= (  - {HP}_{ {\psi}(t_0)}, 	 {DP}_{ {\psi}(t_0)} , 	-  {\bmu}_{P} )(\bgamma_h,\bgamma_h^d).
\end{aligned}$$    
%$$DE(\widetilde{HI}_{t_0}(\bgamma_T,\bgamma_T^d),\widetilde{DI}_{t_0}(\bgamma_T,\bgamma_T^d))=\bgamma_T, \; \widetilde{HI}_{t_0}(DE(\bgamma_h,\bgamma_h^d),\bnu_{DE} )=-HP_{\psi(t_0)}(\bgamma_h,\bgamma_h^d). $$
\item[$(7)$] Suppose that
there exists  $\widetilde{\varphi}:I\rightarrow\R$ such that $m_T(t)\cos\widetilde{\varphi}(t)+ n_T(t)\sin\widetilde{\varphi}(t)=0$ for all $t\in I$.  Then
$$\begin{aligned}( \widetilde{DE},  \widetilde{\bnu}_{DE}, \widetilde{\bmu}_{DE}   ) (\overline{DI}_{t_0}(\bgamma_d,\bgamma_d^h),\overline{\bnu}_{DI}(\bgamma_d,\bgamma_d^h))& = ( \pm\bgamma_d,\bgamma_d^h,\pm\bgamma_d^s   ),\\ 
( \overline{DI}_{t_0} , \overline{\bnu}_{DI} ,\overline{\bmu}_{DI} ) ( \widetilde{DE}(\bgamma_T,\bgamma_T^d), \widetilde{\bnu}_{DE}(\bgamma_T,\bgamma_T^d) ) &= ( \widetilde{DP}_{\widetilde{\varphi}(t_0)},   \widetilde{\bnu}_{DP} , 	  	  \widetilde{\bmu}_{DP} )(\bgamma_T,\bgamma_T^d).
\end{aligned}$$     
\end{enumerate}
\end{theorem}
\demo {\rm(1)} Since the set $\mathrm{Reg}(\bnu_{\mathcal{HI}}^\pm )=\left\{t\in I\mid n_h(t)\pm m_h(t)\neq0\right\}$  is  dense in $I$,    there exists a unique smooth function $ {f}_\mathcal{HI}^\pm(t)= {s}_\pm(t)$ satisfying $ {f}_\mathcal{HI}^\pm(t) {n}_\mathcal{HI}^\pm(t)+ {m}_\mathcal{HI}^\pm(t)=0$ for all $t\in I$. 
Then by Propositions \ref{prop1} and \ref{prop3}, we have
$$\begin{aligned}
	&(\mathcal{HE}^\pm, \mathcal{DE}^\pm,\bmu^\pm_\mathcal{E}  )^T(\mathcal{HI}_{s_\pm}^\pm(\bgamma_h,\bgamma_h^d),\bnu_{\mathcal{HI}}^\pm(\bgamma_h,\bgamma_h^d)) \\ &= N^\mp_1(-f^\pm_{\mathcal{HI}})  R(-\pi/2)( \mathcal{HI}_{s_\pm}^\pm,\bnu_{\mathcal{HI}}^\pm ,   \bmu^\pm_\mathcal{I}  )^T(\bgamma_h,\bgamma_h^d)\\& = N^\mp_1(-f^\pm_{\mathcal{HI}})  R(-\pi/2)N^\mp_2(-s_\pm)  R( \pi/2)(\bgamma_h, \bgamma_h^d,\bgamma_h^s   )^T   \\&= N^\mp_1(-s_\pm) N^\mp_1( s_\pm)   (\bgamma_h, \bgamma_h^d,\bgamma_h^s   )^T = (\bgamma_h, \bgamma_h^d,\bgamma_h^s  ) ^T  	.\\\end{aligned}$$
%%%%
\begin{table}[h]
	\centering
	\caption{The relations between the evolutes and involutes of $(\bgamma_h,\bgamma_h^d)$, $(\bgamma_d,\bgamma_d^h)$ and $(\bgamma_T,\bgamma_T^d)$. Double-line arrows denote correspondence, and single-line arrows labeled $\Delta_1$ and $\Delta_5$ denote $\Delta_1$-duality and $\Delta_5$-duality, respectively. The wavy arrow denotes identity under the condition $(\bgamma_h,\bgamma_h^d)=(\bgamma_d^h,\bgamma_d)$.}
	\vspace{0.3cm}
	\renewcommand{\arraystretch}{1.8}
	\begin{tabular}{ |>{\centering\arraybackslash}p{1.7cm} | >{\centering\arraybackslash}m{1.8cm} | >{\centering\arraybackslash}m{1.8cm} | >{\centering\arraybackslash}m{1.7cm} | >{\centering\arraybackslash}m{1.7cm} | >{\centering\arraybackslash}m{1.7cm} | >{\centering\arraybackslash}m{1.7cm} | }
		\hline
		& \multicolumn{2}{c|}{ $(\bgamma_h,\bgamma_h^d)$} & \multicolumn{2}{c|}{ $(\bgamma_d,\bgamma_d^h)$} &  \multicolumn{2}{c|}{ $(\bgamma_T,\bgamma_T^d)$} \\
		\hline
		Horocyclic evolute 
		& \multicolumn{1}{>{\centering\arraybackslash}m{1.8cm}}{{\scriptsize$\mathcal{HE}^\pm(\bgamma_h,\bgamma_h^d)$:}  {\scriptsize hyperbolic horocyclic  evolute}\tikznode{C1}{}}
		& \multicolumn{1}{>{\centering\arraybackslash}m{1.8cm}|}{{\scriptsize $\mathcal{DE}^\pm(\bgamma_h,\bgamma_h^d)$:}  {\scriptsize de Sitter horocyclic evolute}\tikznode{C2}{}}
		&\multicolumn{2}{c|}{Not exist} 
		& \multicolumn{1}{>{\centering\arraybackslash}m{1.7cm}}{{\scriptsize $\widetilde{\mathcal{DE}}^\pm(\bgamma_T,\bgamma_T^d)$:}  {\scriptsize de Sitter horocyclic evolute}\tikznode{C5}{}}
		& \multicolumn{1}{>{\centering\arraybackslash}m{1.7cm}|}{{\scriptsize $\widetilde{\mathcal{HE}}^\pm(\bgamma_T,\bgamma_T^d)$:} {\scriptsize hyperbolic horocyclic evolute}\tikznode{C6}{}} \\
		\hline
		Horocyclic involute 
		& \multicolumn{1}{>{\centering\arraybackslash}m{1.8cm}}{\tikznode{D1}{}{\scriptsize$\mathcal{HI}_{t_0}^\pm(\bgamma_h,\bgamma_h^d)$:}  {\scriptsize  hyperbolic horocyclic involute}}
		& \multicolumn{1}{>{\centering\arraybackslash}m{1.8cm}|}{\tikznode{D2}{}{\scriptsize $\bnu_{\mathcal{HI}}^\pm(\bgamma_h,\bgamma_h^d)$}  }
		& \multicolumn{1}{>{\centering\arraybackslash}m{1.7cm}}{\tikznode{D3}{}{\scriptsize$\overline{\bnu}^\pm_\mathcal{DI}(\bgamma_d,\bgamma_d^h)$}}
		& \multicolumn{1}{>{\centering\arraybackslash}m{1.7cm}|}{\tikznode{D4}{}{\scriptsize $\overline{\mathcal{DI}}_{t_0}^\pm(\bgamma_d,\bgamma_h^d)$:}  {\scriptsize de Sitter horocyclic involute}}
		& \multicolumn{2}{c|}{Not exist} \\
		\hline
		Evolute 
		& \multicolumn{1}{>{\centering\arraybackslash}m{1.8cm}}{{\scriptsize$HE(\bgamma_h,\bgamma_h^d)$:}   {\scriptsize hyperbolic evolute} \tikznode{A1}{}}
		& \multicolumn{1}{>{\centering\arraybackslash}m{1.8cm}|}{{\scriptsize $DE(\bgamma_h,\bgamma_h^d)$:}  {\scriptsize de Sitter evolute} \tikznode{A2}{}}
		& \multicolumn{1}{>{\centering\arraybackslash}m{1.7cm}}{{\scriptsize$\overline{HE}(\bgamma_d,\bgamma_d^h)$:}  {\scriptsize  hyperbolic evolute}\tikznode{A3}{}}
		& \multicolumn{1}{>{\centering\arraybackslash}m{1.7cm}|}{{\scriptsize $\overline{DE}(\bgamma_d,\bgamma_d^h)$:}  {\scriptsize de Sitter evolute} \tikznode{A4}{}}
		& \multicolumn{1}{>{\centering\arraybackslash}m{1.7cm}}{{\scriptsize $\widetilde{DE}(\bgamma_T,\bgamma_T^d)$:}  {\scriptsize de Sitter evolute} \tikznode{A5}{}}
		& \multicolumn{1}{>{\centering\arraybackslash}m{1.7cm}|}{$\bgamma_T^h$\tikznode{A6}{}} \\
		\hline
		Involute 
		& \multicolumn{1}{>{\centering\arraybackslash}m{1.8cm}}{\tikznode{B1}{}{\scriptsize$HI_{t_0}(\bgamma_h,\bgamma_h^d)$:}  {\scriptsize  hyperbolic involute}}
		& \multicolumn{1}{>{\centering\arraybackslash}m{1.8cm}|}{\tikznode{B2}{}{\scriptsize $DI_{t_0}(\bgamma_h,\bgamma_h^d)$:}  {\scriptsize de Sitter involute}}
		& \multicolumn{1}{>{\centering\arraybackslash}m{1.7cm}}{\tikznode{B3}{}{\scriptsize$\overline{\bnu}_{DI}(\bgamma_d,\bgamma_d^h)$}}
		& \multicolumn{1}{>{\centering\arraybackslash}m{1.7cm}|}{\tikznode{B4}{}{\scriptsize $\overline{DI}_{t_0}(\bgamma_d,\bgamma_d^h)$:}  {\scriptsize de Sitter involute}}
		& \multicolumn{1}{>{\centering\arraybackslash}m{1.7cm}}{\tikznode{B5}{}{\scriptsize $\widetilde{DI}_{t_0}(\bgamma_T,\bgamma_T^d)$:}  {\scriptsize de Sitter involute}}
		& \multicolumn{1}{>{\centering\arraybackslash}m{1.7cm}|}{\tikznode{B6}{}{\scriptsize $\widetilde{HI}_{t_0}(\bgamma_T,\bgamma_T^d)$:}  {\scriptsize  hyperbolic involute}} \\
		\hline
	\end{tabular}
\begin{tikzpicture}[remember picture, overlay, >=Stealth, thin, black]
	
	\draw[<->,black,decorate,decoration={snake,amplitude=0.5mm,segment length=2mm}] ($(A1.north) + (-17pt,36pt)$)   to[bend left=20] ($(A3.north) + (-22pt,38pt)$);
	\draw[<->,black,decorate,decoration={snake,amplitude=0.5mm,segment length=2mm}] ($(A2.north) + (-17pt,36pt)$)   to[bend left=20] ($(A4.north) + (-22pt,38pt)$);
	
	\draw[<->, black, double] ($(A1.south) + (-17pt,0pt)$) -- ($(B1.north) + (23pt,6pt)$);
	\draw[<->, black, double] ($(A2.south) + (-3pt,2pt)$) -- ($(B6.north) + (2pt,6pt)$);
	\draw[<->, black, double] ($(A3.south) + (-27pt,10pt)$) -- ($(B2.north) + (32pt,6pt)$);
	\draw[<->, black, double] ($(A4.south) + (-2pt,10pt)$) -- ($(B5.north) + (18pt,6pt)$);
	\draw[<->, black, double] ($(A5.south) + (-30pt,10pt)$) -- ($(B4.north) + (32pt,5pt)$);
	\draw[<->, black ] ($(B1.south) + (46pt,-11pt)$) -- ($(B2.north) + (5pt,-12pt)$) node[above, xshift=-10pt,yshift=-2pt, font=\footnotesize] {$\Delta_1$};
	\draw[<->, black ] ($(B5.south) + (42pt,-11pt)$) -- ($(B6.north) + (4pt,-12pt)$) node[above, xshift=-9pt,yshift=-3pt, font=\footnotesize] {$\Delta_1$};
	\draw[<->, black ] ($(A5.south) + (6pt,18pt)$) -- ($(A6.north) + (-18pt, 3pt)$) node[above, xshift=-12pt,yshift=-3pt, font=\footnotesize] {$\Delta_1$};
	\draw[<->, black] ($(B3.south) + (47pt,1pt)$) -- ($(B4.north) + ( 6pt, -14.5pt)$)  node[above, xshift=-9pt,yshift=-3pt, font=\footnotesize] {$\Delta_5$};

	\draw[<->, black, double] ($(C1.south) + (-15pt,0pt)$) -- ($(D1.north) + (23pt,6pt)$);
	\draw[<->, black, double] ($(C5.south) + (-30pt,10pt)$) -- ($(D4.north) + (32pt,6pt)$);
	\draw[<->, black ] ($(C1.south) + (7pt,23pt)$) -- ($(C2.north) + (-32pt,22pt)$) node[above, xshift=-10pt,yshift=-3pt, font=\footnotesize] {$\Delta_1$};
	\draw[<->, black ] ($(C5.south) + (5pt,23pt)$) -- ($(C6.north) + (-32pt,22pt)$) node[above, xshift=-10pt,yshift=-3pt, font=\footnotesize] {$\Delta_1$};
	\draw[<->, black ] ($(D1.south) + (43pt,-21pt)$) -- ($(D2.north) + (0pt,0pt)$) node[above, xshift=-10pt,yshift=-3pt, font=\footnotesize] {$\Delta_1$};
	\draw[<->, black] ($(D3.south) + (46pt,1pt)$) -- ($(D4.north) + ( 6pt, -21.5pt)$) node[above, xshift=-9pt,yshift=-3pt, font=\footnotesize] {$\Delta_5$};
\end{tikzpicture}
\end{table}

Note that $(\mathcal{HI}_{s_\pm}^\pm(\bgamma_h,\bgamma_h^d),-\bnu_{\mathcal{HI}}^\pm (\bgamma_h,\bgamma_h^d)) $ is a spacelike Legendre curve with  curvature \\  
$(- {m}^\pm_\mathcal{HI}, {n}^\pm_\mathcal{HI} )$ and $ (-{f}_\mathcal{HI}^\pm(t)) {n}_\mathcal{HI}^\pm(t)- {m}_\mathcal{HI}^\pm(t)=0$ for all $t\in I$. It follows that
$$\begin{aligned}
&(\mathcal{HE}^\pm, \mathcal{DE}^\pm,\bmu^\pm_\mathcal{E}  )^T(\mathcal{HI}_{s_\mp}^\mp(\bgamma_h,\bgamma_h^d),-\bnu_{\mathcal{HI}}^\mp(\bgamma_h,\bgamma_h^d)) \\& =N^\mp_1( f^\mp_{\mathcal{HI}})  R(-\pi/2)( \mathcal{HI}_{s_\mp}^\mp,-\bnu_{\mathcal{HI}}^\mp ,   -\bmu^\mp_\mathcal{I} )^T(\bgamma_h,\bgamma_h^d)\\
&= N^\mp_1( f^\mp_{\mathcal{HI}})  R(-\pi/2) R( \pi )( \mathcal{HI}_{s_\mp}^\mp, \bnu_{\mathcal{HI}}^\mp ,    \bmu^\mp_\mathcal{I}  )^T(\bgamma_h,\bgamma_h^d)\\ &=N^\mp_1( s_\mp )  R(-\pi/2)R( \pi )N^\pm_2(-s_\mp)  R( \pi/2)(\bgamma_h, \bgamma_h^d,\bgamma_h^s   )^T   \\&= R (\pi)     (\bgamma_h, \bgamma_h^d,\bgamma_h^s   )^T   = (\bgamma_h, -\bgamma_h^d,-\bgamma_h^s  ) ^T  	.\\\end{aligned}$$
Consider the Riccati equation 
$$ 
\begin{aligned}\frac{{\rm d} {S}_\mathcal{HE}^\pm(t)}{{\rm d}t} = \frac{m_\mathcal{HE}^\pm(t)\pm n_\mathcal{HE}^\pm(t)}{2}(S_\mathcal{HE}^\pm)^2(t)-m_\mathcal{HE}^\pm(t)= \pm \frac{n_h(t)({S}_\mathcal{HE}^\pm)^2(t)}{2}  +\dot{f}(t)\mp\frac{f^2(t)n_h(t)}{2}.
\end{aligned}$$
Since $$ 
\begin{aligned}\frac{{\rm d}f(t)-{\rm d} {S}_\mathcal{HE}^\pm(t)}{{\rm d}t} &=  \mp \frac{n_h(t)({S}_\mathcal{HE}^\pm)^2(t)}{2}   \pm\frac{f^2(t)n_h(t)}{2}\\&=\mp\left(\frac{(f(t)-  {S}_\mathcal{HE}^\pm(t))^2n_h(t)}{2}+(f(t)-  {S}_\mathcal{HE}^\pm(t))m_h(t)\right),
\end{aligned}$$
then $  f(t)-  {S}_\mathcal{HE}^\pm(t) $ is a solution of equation \eqref{Bernoulli1}. Consequently, 
$$ 
\begin{aligned}
&(\mathcal{HI}_{S ^\pm_\mathcal{HE}}^\pm , \bnu_{\mathcal{HI}}^\pm ,\bmu^\pm_\mathcal{HI} ) ^T(\mathcal{HE}^\pm(\bgamma_h,\bgamma_h^d),\mathcal{DE}^\pm(\bgamma_h,\bgamma_h^d) )\\
&=    N_2^\mp(-S ^\pm_\mathcal{HE})R(\pi/2)(\mathcal{HE}^\pm , \mathcal{DE}^\pm,\bmu^\pm_\mathcal{E}  )^T(\bgamma_h,\bgamma_h^d) \\ & = N_2^\mp(-S ^\pm_\mathcal{HE})R(\pi/2)N_1^{\mp}(-f )R(- \pi/2) (\bgamma_h, \bgamma_h^d,\bgamma_h^s  ) ^T \\
&=  N_2^\mp(f-S ^\pm_\mathcal{HE})    (\bgamma_h, \bgamma_h^d,\bgamma_h^s  ) ^T \\
&=  (\mathcal{HP}_{f-S ^\pm_\mathcal{HE}}^\pm , \bnu_\mathcal{HP}^\pm  ,	 \bmu_\mathcal{HP}^\pm   )  ^T (\bgamma_h, \bgamma_h^d)	.
\end{aligned}$$ For the spacelike Legendre curve $ (\mathcal{HE}^\pm(\bgamma_h,\bgamma_h^d),-\mathcal{DE}^\pm(\bgamma_h,\bgamma_h^d) )$, we denote $\overline{S}_\mathcal{HE}^\pm = -S_\mathcal{HE}^\mp$. It can be  proved that $\overline{S}_\mathcal{HE}^\pm$ is a solution of the Riccati equation 
$$ 
\begin{aligned}\frac{{\rm d} \overline{S}_\mathcal{HE}^\pm(t)}{{\rm d}t} = \frac{-m_\mathcal{HE}^\pm(t)\pm n_\mathcal{HE}^\pm(t)}{2}(\overline{S}_\mathcal{HE}^\pm)^2(t)+m_\mathcal{HE}^\pm(t).
\end{aligned}$$
Therefore, $$ 
\begin{aligned}
& (\mathcal{HI}_{\overline{S} ^\pm_\mathcal{HE}}^\pm , \bnu_{\mathcal{HI}} ^\pm ,\bmu^\pm_\mathcal{HI}  ) ^T(\mathcal{HE}^\mp(\bgamma_h,\bgamma_h^d),-\mathcal{DE}^\mp(\bgamma_h,\bgamma_h^d))\\
&=    N_2^\mp(-\overline{S} ^\pm_\mathcal{HE})R(\pi/2)(\mathcal{HE}^\mp, -\mathcal{DE}^\mp,-\bmu^\mp_\mathcal{E}  )^T(\bgamma_h,\bgamma_h^d)\\ 
&=    N_2^\mp(-\overline{S} ^\pm_\mathcal{HE})R(\pi/2)R(\pi)(\mathcal{HE}^\mp ,  \mathcal{DE}^\mp , \bmu^\mp_\mathcal{E}  )^T(\bgamma_h,\bgamma_h^d)\\ 
&=  N_2^\mp(-\overline{S} ^\pm_\mathcal{HE})R(\pi/2)R(\pi)N_1^{\pm}(-f )R(- \pi/2) (\bgamma_h, \bgamma_h^d,\bgamma_h^s ) ^T \\
&=  R(\pi)N_2^\pm(f+\overline{S} ^\pm_\mathcal{HE})    (\bgamma_h, \bgamma_h^d,\bgamma_h^s ) ^T \\
&=   (\mathcal{HP}_{f+\overline{S}^\pm_\mathcal{HE}}^\mp  , -\bnu_\mathcal{HP}^\mp  ,	 -\bmu_\mathcal{HP}^\mp  )  ^T (\bgamma_h, \bgamma_h^d)	.
\end{aligned}$$

{\rm(2)}  Similar discussion yields assertion (2). 

{\rm(3)} Since 
$( HI_{t_0}(\bgamma_h,\bgamma_h^d), DI_{t_0}(\bgamma_h,\bgamma_h^d) ) $ is a spacelike Legendre curve with curvature 
$$
\left(m _{HI}(t), n_{HI}(t)\right)= \left(     n_h(t)\sinh\left(\int_{t_0}^t{m_h(t)}{\rm d}t\right), n_h(t)\cosh\left(\int_{t_0}^t{m_h(t)}{\rm d}t\right)\right) 
$$ and the set $\mathrm{Reg}(\bgamma_h^d)= \{t\in I\mid n_h(t)\neq0 \}$  is   dense  in $ I $, there exists a unique smooth function $\varphi_{HI}(t)=-\int_{t_0}^t{m_h(t)}{\rm d}t$ such that  $m_{HI}(t)\cosh\varphi_{HI}(t)+ n_{HI}(t)\sinh\varphi_{HI}(t)=0$  for all $t\in I$. Then  
$$\begin{aligned}
&(HE , \bnu_{HE}  ,\bmu_{HE}  )^T (HI_{t_0}(\bgamma_h,\bgamma_h^d),DI_{t_0}(\bgamma_h,\bgamma_h^d)) \\
&= B_2(-\varphi_{HI} )R( -\pi/2)( HI_{t_0} , DI_{t_0} ,\bmu_I   )^T(\bgamma_h,\bgamma_h^d)\\&= 	B_2(-\varphi_{HI} )R( -\pi/2)B_1\left( \int_{t_0}^t{m_h }{\rm d}t\right)R(\pi/2) \left(   \bgamma_h,\bgamma_h^d, \bgamma_h^s  \right)^T \\
&= 	B_2\left(\int_{t_0}^t{m_h }{\rm d}t\right) B_2\left(- \int_{t_0}^t{m_h }{\rm d}t\right) (   \bgamma_h,\bgamma_h^d, \bgamma_h^s  )^T 
=  (   \bgamma_h,\bgamma_h^d, \bgamma_h^s  )^T.
\end{aligned}$$
Moreover, 
$$\begin{aligned}
&( HI_{t_0} , DI_{t_0},\bmu_I  )^T(HE(\bgamma_h,\bgamma_h^d),\bnu_{HE}(\bgamma_h,\bgamma_h^d) )\\&=  B_1\left( \int_{t_0}^t{m_{HE} }{\rm d}t\right)R(\pi/2) \left(  HE,\bnu_{HE},\bmu_{HE}\right)^T(\bgamma_h,\bgamma_h^d) \\
&=  B_1\left( \int_{t_0}^t-\dot{\varphi}  {\rm d}t\right)R(\pi/2)B_2(-\varphi  )R( -\pi/2) (   \bgamma_h,\bgamma_h^d, \bgamma_h^s )^T\\ 
&= B_1\left( \int_{t_0}^t-\dot{\varphi}  {\rm d}t\right)B_1( \varphi  ) (   \bgamma_h,\bgamma_h^d, \bgamma_h^s )^T \\ 
&=  B_1(  {\varphi}(t_0)   )  (   \bgamma_h,\bgamma_h^d, \bgamma_h^s  )^T= (  	HP_{\varphi(t_0)}, DP_{\varphi(t_0)},	\bmu_{P} )^T(\bgamma_h,\bgamma_h^d).
\end{aligned}$$

Since we can prove from (4) to (7) by the similar calculations, we omit the details here.
\enD

%%%%%%%%%%%%%%%%%%%%%%%%%%%%%%%%%%%%%%%%%%%%%%%%%%%%%%%%%%

%%%%%%%
Nozomi Nakatsuyama, 
\\
Muroran Institute of Technology, Muroran 050-8585, Japan,
\\
E-mail address: 25096009b@muroran-it.ac.jp
\\
\\
Masatomo Takahashi, 
\\
Muroran Institute of Technology, Muroran 050-8585, Japan,
\\
E-mail address: masatomo@muroran-it.ac.jp
\\
\\
Anjie Zhou,
\\
School of Mathematics and Statistics, Northeast Normal University, Changchun 130024, P. R. China
\\
E-mail address: zhouaj882@nenu.edu.cn

\end{document}